\documentclass{elsarticle}
\usepackage[utf8]{inputenc}

\usepackage{tikz-cd}
\usepackage{url}

\usepackage{graphicx}
\usepackage{caption}
\usepackage{subcaption}

\usepackage{xspace}
\usepackage{xcolor}
\definecolor{darkorange}{rgb}{0.8, 0.4, 0.0}
\definecolor{darkmagenta}{rgb}{0.55, 0.0, 0.55}

\usepackage{epstopdf}
\usepackage{comment}

\usepackage{amsmath,amssymb}
\usepackage{dsfont}
\usepackage{bm}
\usepackage{mathrsfs}
\usepackage{amsthm}
\usepackage{a4wide}
\usepackage{stmaryrd}

\usepackage{multirow}
\usepackage{derivative}
\usepackage{arydshln,mathtools}
\usepackage[most]{tcolorbox}
\usepackage{soul}

\usepackage[numbers]{natbib}
\newtheorem{theorem}{Theorem}

\newtheorem{remark}{Remark}
\newtheorem{proposition}{Proposition}

\renewcommand\d{\ensuremath{\mathrm{d}}}

\newcommand{\bbR}{\mathbb{R}}

\DeclareMathOperator{\tr}{tr}

\newcommand*{\dual}[1]{\ensuremath{\widehat{#1}}}

\newcommand{\inprDom}[3][M]{\ensuremath{( #2, \, #3 )_{#1}}}
\newcommand{\inprBd}[3][\partial M]{\ensuremath{(#2 , #3 )_{#1}}}

\newcommand{\dualprDom}[3][M]{\ensuremath{\langle #2 \; \vert #3 \rangle_{#1}}}
\newcommand{\dualprBd}[3][\partial M]{\ensuremath{\langle #2 \; \vert #3 \rangle_{#1}}}

\DeclareMathOperator*{\grad}{grad}

\renewcommand{\div}{\operatorname{div}}

\DeclareMathOperator*{\curl}{curl}

\newcommand{\firedrake}{\textsc{Firedrake}\xspace}

\newcommand{\revone}[1]{{\color{black}#1}}
\newcommand{\revtwo}[1]{{\color{black}#1}}

\journal{Applied Mathematics and Computation}

\graphicspath{{./images_revision}}

\begin{document}

\begin{frontmatter}	
		
		\title{Finite element hybridization of port-Hamiltonian systems}	
        
            \author[ICA]{Andrea Brugnoli\corref{cor1}}
		\ead{andrea.brugnoli@isae-supaero.fr}
		\cortext[cor1]{Corresponding author}
  
	    \author[FAHD]{Ramy Rashad}
	    \ead{ramy.rashad@kfupm.edu.sa}
     
	    \author[GUI]{Yi Zhang}
	    \ead{zhangyi_aero@hotmail.com}
            \author[UT]{Stefano Stramigioli}
            \ead{s.stramigioli@utwente.nl}

            \address[ICA]{ICA, Universit\'e de Toulouse, ISAE–SUPAERO, INSA, CNRS, MINES ALBI, UPS, Toulouse, France}
            \address[FAHD]{Control and Instrumentation Engineering Department, King Fahd University of Petroleum and Minerals, Saudi Arabia}
            \address[GUI]{School of Mathematics and Computing Science, Guilin University of Electronic Technology, Guilin, China}
            \address[UT]{Robotics and Mechatronics Department, University of Twente, The Netherlands}

		\begin{abstract}
            In this contribution, we extend the hybridization framework for the Hodge Laplacian [Awanou et  al., \textit{Hybridization and postprocessing in finite element exterior calculus, 2023}] to  port-Hamiltonian systems describing linear wave propagation phenomena. To this aim, a dual field mixed Galerkin discretization is introduced, in which one variable is approximated via conforming finite element spaces, whereas the second is completely local. \revone{The mixed formulation is then hybridized to obtain an equivalent formulation that can be more efficiently solved using a static condensation procedure in discrete time. The size reduction achieved thanks to the hybridization is greater than the one obtained for the Hodge Laplacian as the final system only contains the globally coupled traces of one variable.}  Numerical experiments on the 3D wave and Maxwell equations illustrate the convergence of the method and the size reduction achieved by the hybridization.
		\end{abstract}
		
		\begin{keyword}
		Port-Hamiltonian systems \sep Finite element exterior calculus \sep  Hybridization \sep Dual field  
		\end{keyword}

	\end{frontmatter}
	
\section{Introduction}

Distributed port-Hamiltonian systems were first introduced in \cite{vanderSchaft2002} (see \cite{rashad2020review} for a comprehensive review) as an extension of Hamiltonian PDEs undergoing boundary interaction. The framework relies on a particular geometrical structure, that of a Dirac manifold~\cite{courant1990}.  This geometrical structure \revone{does not specify the actual boundary conditions of the problem at hand but captures all admissible boundary flows. Therefore, it defines a property of families of solutions.}  Because of the fact that Dirac manifolds are composable, port-Hamiltonian systems are closed under interconnection, and thus well suited for modular modelling of complex systems. The Dirac structure underlying port-Hamiltonian systems is characterized by the topological integration by parts formula (i.e. the combination of Leibniz rule and Stokes theorem), based on the Hodge duality of differential forms. 
\revtwo{The construction a non-degenerate duality product of forms is intimately related to that of a discrete Hodge star, since this operator is precisely the Riesz map. Obtaining an isometric discrete Hodge star constitutes a major difficulty in a finite element context \cite{hiptmair2001,hirani2003discrete}. If only computational mesh is used, a weak interpretation of the Hodge star leads to an $L^2$ projection between finite element spaces \cite{hiptmair2001}. This projection entails the loss of the isometric nature of the Hodge star and the duality pairing of forms becomes degenerate. To circumvent this problem, a primal dual formulation, either based on primal dual meshes or on dual spaces of finite elements, can be used.} Several recent contributions have explored the latter strategy to obtain a discretization of the Hodge operator. In \cite{kumar2022port} the authors show that the duality between he complete and trimmed polynomial families of finite element exterior calculus can be used to construct a non degenerate duality product. The higher regularity of splines has been used in \cite{kapidani2022} to construct Hodge dual isomorphic sequences. Our previous work \cite{brugnoli2022df} relies on finite element exterior calculus \cite{arnold2006acta} to discretize linear hyperbolic port-Hamiltonian  models (that represent a special case of the more general Hodge wave equation \cite{wu2021error}). \revtwo{The discretization strategy relies on a primal dual representation, in which the Hodge star arising from the constitutive equations is embedded in the codifferential}. This idea is at the core of finite element exterior calculus, where the codifferential operator has a weak interpretation by means of integration by parts. Even though the primal dual structure of wave propagation problems has been already analyzed in \cite{joly2003}, the connection with exterior calculus and the Hodge duality is the main novelty of our previous work. \\

\revone{Hybrid finite element methods are variations of continuous finite element methods that introduce additional variables, such as Lagrange multipliers or inter-element boundary terms, to enforce continuity or other constraints more flexibly \cite{cockburn2009unified}. This leads to several practical and theoretical benefits. Concerning the computational efficiency, an important advantage is given by the so called static condensation \cite{guyan1965reduction}. Degrees of freedom for discontinuous function spaces can be locally eliminated using a Schur complement. The resulting system of equations only involves degrees of freedom for the globally coupled boundary Lagrange multipliers. This leads to a smaller system to be solved than
the original problem. Furthermore, via local postprocessing the quality of the approximated solution can be
improved. Static condensation can be used in both static and dynamic problems. However in the latter case, instead of condensing internal degrees of freedom, Lagrange multipliers may be eliminated by regularizing the differential algebraic system arising from the hybridization. The resulting differential system may be solved in parallel as illustrated in \cite{park2023partitioned}. This strategy leads to other advantages like parallel solutions of quasi static problems (where a suitable preconditioner of the stiffness matrix plays the role of the mass matrix) and fluid structure interaction problems \cite{gonzales2023three}. On a theoretical level the Lagrange multiplier functions correspond to weak boundary traces of solution components, thus providing additional information about the
solution. This information may be crucial in a port-Hamiltonian context, as the trace variables correspond to the boundary port of the Dirac structure and represent the power flow between different elements.  The hybridization of the Hodge Laplacian has been recently presented in \cite{stern2023hyb} and the present contribution is built upon the framework presented therein. \\}

\revtwo{In this work, we consider the hybridization of a slightly different primal-dual formulation with respect to \cite{brugnoli2022df}. Linear hyperbolic port-Hamiltonian systems are special cases of the wave Hodge equation and lack one part of the Hodge Laplacian is missing. Therefore, the variational mixed formulation does not require both variables to be conforming and one variable (that does not undergo differentiation) can be taken to be less regular, i.e. in a broken finite element space (this choice is common to other previous mixed formulations \cite{cohen2002higher}). To ensure discrete conservation properties, this discontinuous space is chosen to verify local discrete subcomplex property and the resulting formulation is equivalent to the one presented \cite{brugnoli2022df}. Following the framework described in \cite{stern2023hyb}, the conforming variable in the mixed Galerkin scheme is then hybridized by introducing a broken (and therefore local) multiplier, capturing the information related to the normal trace, and an unbroken global multiplier, representing the tangential trace of the regular variable.} Explicit time integration schemes like the leapfrog scheme (also called St\"ormer-Verlet in the context of Hamiltonian dynamics) are typically preferred in the context of wave propagation. This is due to the fact that by introducing mass lumping strategies, the resulting linear system is block diagonal and can be efficiently solved \cite{egger2021mass}. However, we are here concerned with exact preservation of the power balance as it implies preservation of the Dirac structure. For this reason, the time discretization is obtained by using the implicit midpoint method, that provides an exact discrete power balance \cite{kirby2015wave,kotyczka2019discrete}. This implicit method can be directly applied to the index two differential-algebraic system arising from the hybridization procedure. The resulting system is a saddle point problem that can be efficiently solved via static condensation. This leads to a considerable reduction of the number of unknowns as the final linear only contains the global trace variable of the continuous field. Local variables can be subsequently computed in parallel, as they are completely uncoupled. \revone{Because of the peculiar structure of port-Hamiltonian systems, the size reduction obtained by means of the hybridization is even more important than the one obtained for the Hodge Laplacian in \cite{stern2023hyb}, as the continuity of only one variable needs to be enforced via a global trace multiplier.}\\

The paper is organized as follows: in Sec. \ref{sec:prel} a brief discussion of the $L^2$ theory of differential forms and the notation is presented. In Sec. \ref{sec:contGal} we introduce the modified mixed Galerkin formulation by making use of broken finite elements for differential forms. The hybridization of this formulation is presented in Sec. \ref{sec:hybridization}, where it is shown that the hybrid version is completely equivalent to the mixed formulation. In Sec. \ref{sec:realization}, the algebraic realization of the weak formulation is detailed. The reinterpretation of finite element assembly as interconnection of port-Hamiltonian descriptor system is also given. A time discrete system is obtained using the implicit midpoint method. Via static condensation the local variables are then eliminated. Section \ref{sec:num_exp} presents numerical experiments for the wave and Maxwell equations in 3D.

\section{Preliminaries}\label{sec:prel}

\subsection{Smooth theory of differential forms}
For sake of conciseness, we assume that the reader is familiar with the basic operators and results of exterior calculus (wedge product, exterior derivative, trace operator, Stokes theorem). \\

In the following $M$ is a Riemannian manifold of dimension $n$. In this paper we distinguish true forms and pseudo forms. The Hodge star maps inner-oriented (or true) forms to outer-oriented (or pseudo) forms \cite{kreeft2011mimetic,frankel2011geometry} and vice versa. This distinction allows defining integral quantities that are orientation independent (like mass, energy, etc.). In this paper, \revone{the space of differential form of degree $k$ is denoted by $\Omega^k(M)$, instead the space of outer-oriented forms carry a hat as its elements, i.e. $\dual{\alpha}^k \in \dual{\Omega}^k(M)$}. 
The Hodge-$\star$ operator $\star : \Omega^k(M) \rightarrow \dual{\Omega}^{n-k}(M)$ is such that
\begin{equation*}
    \alpha^k \wedge {\star \beta^k} = \inprDom{\alpha^k}{\beta^k} {\mathrm{vol}}, \qquad \alpha^k, \beta^k \in \Omega^k(M),
\end{equation*}
where $\inprDom{\alpha^k}{\beta^k}$ is the pointwise inner product of forms and ${\mathrm{vol}}$ the standard volume form.

Differential forms of complementary degree with respect to the manifold dimension possess a natural duality product
\begin{equation}\label{eq:dual_pr}
    \dualprDom[M]{\alpha^k}{\dual{\beta}^{n-k}} := \int_M \alpha^k \wedge \dual{\beta}^{n-k}, \qquad \alpha^{k} \in \Omega^{k}(M), \quad \dual{\beta}^{n-k} \in \dual{\Omega}^{n-k}(M), \quad k=0, \dots, n.
\end{equation}
An inner and an outer oriented forms are paired so that the resulting scalar quantity is orientation independent. The duality product is also defined on the boundary $\partial M$ whose orientation is inherited from the initial oriented manifold\footnote{If $M$ is not oriented $\partial M$ has always a transverse orientation pointing outwards.} $M$ ($k=0, \dots, n-1$)
\begin{equation}
    \dualprBd[\partial M]{\tr \alpha^k}{\tr \dual{\beta}^{n-k-1}} := \int_{\partial M} \tr \alpha^k \wedge \tr \beta^{n-k-1}, \qquad \alpha^k \in \Omega^{k}(M), \quad  \dual{\beta}^{n-k-1} \in \dual{\Omega}^{n-k-1}(M).
\end{equation}
Combing the Leibniz rule and the Stokes theorem, one has the integration by parts formula 
\begin{equation}\label{eq:int_byparts_d}
    \dualprDom[M]{\d\alpha}{\dual{\beta}} + (-1)^k \dualprDom[M]{\alpha}{\d\dual{\beta}} = \dualprBd[\partial M]{\tr \alpha}{\tr \dual{\beta}}, \qquad \alpha \in \Omega^{k}(M), \; \beta \in \dual{\Omega}^{n-k-1}(M), \quad k=0, \dots, n-1.
\end{equation}

\revtwo{
\begin{remark}
This duality pairing is important in port-Hamiltonian systems as it defines the power flow. Indeed the time derivative of the Hamiltonian is given by a duality pairing of forms where the variational derivative of the Hamiltonian is paired with the time derivative of a state variable, as follows    
$$
\dot{H}(\dual{\alpha}^p, \beta^q) = \dualprDom{\delta_{\dual{\alpha}} H}{\dot \alpha^p} + \dualprDom{\delta_\beta H}{\dot \beta^q},\qquad p+q=n+1
$$
where $\delta_{\dual{\alpha}} H \in \Omega^{n-p}(M), \; \delta_{\beta} H \in \Omega^{n-q}(M)$. A more detailed discussion on distributed port-Hamiltonian systems is found in \cite{vanderSchaft2002,brugnoli2022df}.
\end{remark}
}

\subsection{$L^2$ theory of differential forms}
As in vector calculus, the $L^2$ Hilbert space of differential forms is  the completion of the space of smooth forms $\Omega^k(M)$ in the norm induced by the $L^2$ inner product $\inprDom[M]{\cdot}{\cdot}$. Taking $\d$ in the sense of distributions allows it to be extended to a closed, densely defined operator with domain 
$$
H\Omega^k(M):= \{\omega^k \in L^2\Omega^k(M) | \; \d\omega^k \in L^2\Omega^k(M)  \},
$$
the space of $k$ forms belonging to $L^2$ with weak derivative in $L^2$. These spaces, connected by the operator $\d$, form the de Rham domain complex. The formal  adjoint of the exterior derivative with respect to the $L^2$ inner product is the codifferential operator  $\d^* : \Omega^k(M) \xrightarrow{} \Omega^{k-1}(M)$ that satisfies the following integration by parts formula
\begin{equation}\label{eq:intbyparts_codif}
    \inprDom[M]{\alpha^k}{\d^* \beta^{k+1}} = \inprDom[M]{\mathrm{d} \alpha^{k}}{\beta^{k+1}} - \dualprBd[\partial M]{\tr \alpha^{k}}{{\tr \star \beta^{k+1}}}, \qquad \alpha \in \Omega^k(M), \; \beta \in \Omega^{k+1}(M).
\end{equation}
Analogously to the exterior derivative, the codifferential $\d^*$ may be also extended to a closed, densely defined operator with domain $H^* \Omega^k(M)$. The integration by parts formula \eqref{eq:intbyparts_codif} can be rewritten using the normal trace \cite{stern2023hyb}. Given $\omega^k \in \Omega^k(M)$, its normal trace is a $k-1$ form on the boundary $\partial M$ defined by
\begin{equation*}
    \tr_{\bm{n}} \omega^k = \star^{-1}_\partial \tr \star \omega^k,
\end{equation*}
where ${\star}_\partial$ denotes the Hodge star at the boundary (that is constructed using the pullback of the metric at the boundary and the associated volume form). The integration by parts formula \eqref{eq:intbyparts_codif} is then rewritten using the inner product on the boundary, denoted by $\inprBd[\partial M]{\cdot}{\cdot}$. \revtwo{In \cite{weck2004,schulz2022boundary}, it is shown that the integration by parts formula \eqref{eq:intbyparts_nortr} can be extended to $\alpha^k \in H\Omega^k(M)$ and $\beta^{k+1} \in H^*\Omega^{k+1}(M)$ and to manifolds with Lipschitz boundary} 
\begin{equation}\label{eq:intbyparts_nortr}
    \inprBd[\partial M]{\tr \alpha^{k}}{{\tr_{\bm{n}} \beta^{k+1}}} = \inprDom[M]{\mathrm{d} \alpha^{k}}{\beta^{k+1}} - \inprDom[M]{\alpha^k}{\d^* \beta^{k+1}}, \qquad \alpha^k \in H\Omega^k(M), \; \beta^{k+1} \in H^*\Omega^{k+1}(M).
\end{equation}

\revtwo{The spaces where the trace variables $\tr \alpha^k, \; \tr_{\bm{n}} \beta^{k+1}$ live are denoted by $H\Omega^{k, \bm{t}}(\partial M)$ and $H^*\Omega^{k, \bm{n}}(\partial M)$ respectively. These are subspaces of the fractional Sobolev space\footnote{The space $H^{-1/2}{\Omega}^{k}(\partial M)$ arise from the fact that the trace operator on $H^1\Omega^k(M)$ (the space  of $k$-forms whose coefficients are in $H^1(M)$) extends to a bounded linear operator $\tr: H\Omega^k(M) \rightarrow H^{-1/2}\Omega^{k}(\partial M)$ \cite[Theorem 6.3]{arnold2018finite}.} $H^{-1/2}\dual{\Omega}^k(\partial M)$.
Thus, $\inprBd[\partial M]{\tr \alpha^{k}}{{\tr_{\bm{n}} \beta^{k+1}}}$  represents a duality pairing\footnote{It is important to underline the difference between $\dualprBd[\partial M]{\cdot}{\cdot}$ and $\inprBd[\partial M]{\cdot}{\cdot}$. The first is given by the duality of forms and does not need to be symmetric. The second is the extension of the $L^2$ inner product on the boundary to  distributional spaces and it is symmetric.} extending the $L^2$ inner product on the boundary. The exact definitions of $H\Omega^{k, \bm{t}}(\partial M), \; H^*\Omega^{k, \bm{n}}(\partial M)$ are involved. However, these spaces are isomorphic to the following quotient spaces (cf. \cite[Thms. 5 and 7]{Weck2004traces})
\begin{equation*}
\begin{aligned}
    H\Omega^{k, \bm{t}}(\partial M) &\cong H\Omega^{k}(M)/ \mathring{H}\Omega^{k}(M), \\
    H^*\Omega^{k, \bm{n}}(\partial M) &\cong H^*\Omega^{k+1}(M)/ \mathring{H}^*\Omega^{k+1}(M).
\end{aligned}
\end{equation*}
The subspaces $\mathring{H}\Omega^{k}(M) \subset {H}\Omega^{k}(M)$ and $\mathring{H}^*\Omega^{k}(M) \subset {H}^*\Omega^{k}(M)$ are defined as the closure of compactly supported smooth function in the ${H}\Omega^{k}(M)$ and ${H}^*\Omega^{k}(M)$ norm respectively. Therefore, $H\Omega^{k, \bm{t}}(\partial M)$ and $H^*\Omega^{k, \bm{n}}(\partial M)$ can be treated as quotient spaces equipped with the quotient norm. This equivalent characterization allows deducing that
 $\mathring{H}\Omega^{k}(M)$ and $\mathring{H}^*\Omega^{k}(M)$ are the spaces of function with vanishing tangential and normal trace respectively (see \cite[Lemma 2.3]{stern2023hyb} for details)
\begin{equation*}
\begin{aligned}
    \mathring{H}\Omega^{k}(M) &= \{{\omega}^{k} \in H{\Omega}^{k}(M) \vert \; \tr {\omega}^{k}\vert_{\partial M} = 0 \}, \\
    \mathring{H}^*\Omega^{k}(M) &= \{{\omega}^{k} \in H^*{\Omega}^{k}(M) \vert \; \tr_{\bm{n}} {\omega}^{k}\vert_{\partial M} = 0 \}.
\end{aligned} 
\end{equation*}}

For the reader convenience, Table \ref{tab:inner_dual} resumes the notation used for inner and duality products over the domain $M$ and the boundary $\partial M$.

\begin{table}[tbh]
    \centering
    \begin{tabular}{c|c|c}
        & Inner product  & Dual Product   \\
        \hline
        \rule{0pt}{.45cm}
        Domain $M$ & $\inprDom[M]{\alpha}{\beta} = \int_M \alpha \wedge \star \beta$ &  $\dualprDom[M]{\alpha}{\dual{\beta}} = \int_{M} \alpha \wedge \dual{\beta}$ \\
        \rule{0pt}{.4cm}
        Boundary $\partial M$ & $\inprBd[\partial M]{\alpha}{\beta} = \int_{\partial M} \alpha \wedge \star_\partial \beta$  & $\dualprBd[\partial M]{\alpha}{\dual{\beta}} = \int_{\partial M} \alpha \wedge \dual{\beta}$
    \end{tabular}
    \caption{\revtwo{Inner products are denoted using round brackets $(\;, )$, whereas duality products are denoted using angle brackets $\langle \;| \rangle$. The spatial location, i.e. the domain $M$ or the boundary $\partial M$ is indicated with a subscript.}}
    \label{tab:inner_dual}
\end{table}

\subsection{Broken Sobolev spaces and decomposition of Hilbert complexes}
We briefly recall here some concepts related to decomposition of Hilbert complex, following closely \cite{stern2023hyb}. Additional details over broken Sobolev spaces can be found in \cite{bacuta2021}. \\

\revtwo{From now on we will assume that $M \subset \mathbb{R}^n$ and $\mathcal{T}_h$ is a regular mesh corresponding to the tessellation of the domain $M$}. We denote single elements of $\mathcal{T}_h$ by $T$. On each element the $L^2\Omega^k(T)$ and $H\Omega^k(T)$ space can be defined. In particular $H\Omega^k(T) \subset L^2\Omega^k(T)$. The space $L^2\Omega^k(\mathcal{T}_h)$ is isomorphic to $L^2(M)$ due to the properties of Lebesgue integration. The broken Sobolev spaces are the product spaces on the disjoint union $\bigsqcup_{T \in\mathcal{T}_h} T$
\begin{equation}
    H\Omega^k(\mathcal{T}_h) := \prod_{T \in \mathcal{T}_h} H\Omega^k(T)
\end{equation}
with inner product inherited by the Cartesian product structure. \revone{So a broken Sobolev space is a Sobolev space on each element of the computational mesh but non on the overall domain.} Sobolev spaces are naturally included in their broken versions $H\Omega^k(M) \hookrightarrow H\Omega^k(\mathcal{T}_h)$ by restriction over each cell. Broken versions of the differential operators are defined element-wise, i.e. $$\d : H\Omega^k(\mathcal{T}_h) \rightarrow H\Omega^{k+1}(\mathcal{T}_h) : = \d|_{H\Omega^k(T)} \quad \forall  T\in \mathcal{T}_h.$$

 \revtwo{In the context of broken forms, the integration by parts formula is rewritten in a split form in which the broken nature of the variables. Indeed traces of broken differential forms are double valued, as no continuity is imposed at interfaces
\begin{equation*}
    \sum_{T \in \mathcal{T}_h} \inprDom[T]{\d \omega}{\sigma} - \inprDom[T]{\omega}{\d^* \sigma} = \sum_{T \in \mathcal{T}_h} \inprBd[\partial T]{\tr \omega}{\tr_{\bm{n}} \sigma}, \qquad \forall\; \omega \in H\Omega^{k}(\mathcal{T}_h), \quad \sigma \in H^*\Omega^{k+1}(\mathcal{T}_h).
\end{equation*}
This expression expresses the fact that tangential and normal traces of broken forms live on the mesh's disjoint facets $\partial \mathcal{T}_h:= \bigsqcup_{T \in\mathcal{T}_h} \partial T$. Given the following definitions 
$$\inprDom[\mathcal{T}_h]{\cdot}{\cdot} := \sum_{T \in \mathcal{T}_h} \inprDom[T]{\cdot}{\cdot}, \qquad \inprBd[\partial \mathcal{T}_h]{\cdot}{\cdot} := \sum_{T \in \mathcal{T}_h} \inprBd[\partial T]{\cdot}{\cdot},$$
a compact expression for the integration by parts formula is readily obtained 
\begin{equation*}
    \inprDom[\mathcal{T}_h]{\d \omega}{\sigma} - \inprDom[\mathcal{T}_h]{\omega}{\d^* \sigma} = \inprBd[\partial\mathcal{T}_h]{\tr \omega}{\tr_{\bm{n}} \sigma}, \qquad \forall\; \omega \in H\Omega^{k}(\mathcal{T}_h), \quad \sigma \in H^*\Omega^{k+1}(\mathcal{T}_h).
\end{equation*}
}
The following proposition, corresponding to Pr. 3.1. in  \cite{stern2023hyb}, gives a characterization of unbroken Sobolev spaces as single-valued broken spaces.

\begin{proposition}\label{pr:H_br_uniform}
Let $\mathcal{T}_h$ be a Lipschitz decomposition of $M$, \revone{a decomposition of the domain into Lipschitz regular partitions,} then \begin{equation*}
H\Omega^k(M) = \{\omega \in H\Omega^k(\mathcal{T}_h): \inprBd[\partial\mathcal{T}_h]{\tr \omega}{\tr_{\bm{n}} \tau}=0,\quad \forall \tau  \in \mathring{H}^* \Omega^{k+1}(M)\}. 
\end{equation*}

\end{proposition}

\subsection{The primal-dual structure of port-Hamiltonian systems}

By combining the canonical port-Hamiltonian system and its adjoint, two different formulations are deduced \cite{brugnoli2022df}. The variables in one system are the Hodge dual of those of the second one. \revtwo{In the following it is assumed that $p, q$ are two integers such that $p +q = n+ 1$ (cf. \cite{vanderSchaft2002}).}

\paragraph{Primal system of outer oriented forms}
Find $\dual{\alpha}^p: (0, T] \rightarrow \dual{\Omega}^p(M)$ and $\dual{\beta}^{p-1}: (0, T] \rightarrow  \dual{\Omega}^{p-1}(M)$ such that
 \begin{equation}\label{eq:strong_primalPH}
    \odv{}{t}\begin{pmatrix}
        \dual{\alpha}^p \\
        \dual{\beta}^{p-1}
    \end{pmatrix} = (-1)^{p}
    \begin{bmatrix}
        0 & \d \\
        -\d^* & 0\\
    \end{bmatrix}
    \begin{pmatrix}
        \dual{\alpha}^p \\
        \dual{\beta}^{p-1}
    \end{pmatrix}, \qquad
    \begin{aligned}
        \tr {\star \dual{\alpha}^p}|_{\Gamma_1} = u_1^{q-1}, \\
        (-1)^p\tr \dual{\beta}^{p-1}|_{\Gamma_2} = \dual{u}_2^{p-1}, \\
    \end{aligned}
 \end{equation}
 with initial condition $\dual{\alpha}^p(0) = \dual{\alpha}^p_{0}, \; \dual{\beta}^{p-1}(0) = \dual{\beta}^{p-1}_{0}$.
 The collocated outputs are given by
 \begin{equation}
      y^{q-1}_1 := \tr \star \dual{\alpha}^{p} \vert_{\Gamma_2}, \qquad
    \dual{y}^{p-1}_2 := (-1)^p \tr \dual{\beta}^{p-1} \vert_{\Gamma_1}. 
 \end{equation}

 \revtwo{\begin{remark}
     Collocation is a term borrowed from control theory. It refers to the fact that output is situated at the same location of the input and it is the conjugated variable for the energetic behavior \cite{vanderSchaft2002}, leading to the power balance $\dot{H} = \dualprBd[\partial M]{y}{u}$, where $H$ is the energy of the system.
 \end{remark}}

\paragraph{Dual system of inner oriented forms}
Find ${\alpha}^{q-1}: (0, T] \rightarrow {\Omega}^{q-1}(M)$ and ${\beta}^{q}: (0, T] \rightarrow  {\Omega}^{q}(M)$ such that
\begin{equation}\label{eq:strong_dualPH}
    \odv{}{t}\begin{pmatrix}
        {\alpha}^{q-1} \\
        {\beta}^{q}
    \end{pmatrix} =
    \begin{bmatrix}
        0 & \d^* \\
        -\d & 0\\
    \end{bmatrix}
    \begin{pmatrix}
        {\alpha}^{q-1} \\
        {\beta}^{q}
    \end{pmatrix}, \qquad
    \begin{aligned}
        \tr \alpha^{q-1}|_{\Gamma_1} = u_1^{q-1}, \\
        (-1)^{p+q(n-q)}\tr {\star \beta^{q}}|_{\Gamma_2} = \dual{u}_2^{p-1}, \\
    \end{aligned}
 \end{equation}
with initial condition ${\alpha}^{q-1}(0) = {\alpha}^{q-1}_{0}, \; {\beta}^{q}(0) = {\beta}^{q}_{0}$. The outputs are given by
 \begin{equation}\label{eq:y_dual}
      y^{q-1}_1 := \tr {\alpha}^{q-1} \vert_{\Gamma_2}, \qquad
    {y}^{p-1}_2 := (-1)^{p+q(n-q)} \tr {\star{\beta}^{q}} \vert_{\Gamma_1}. 
 \end{equation}

\section{Mixed Galerkin discretization of port-Hamiltonian systems}\label{sec:contGal}

\revtwo{A finite element exterior calculus discretization of the primal  and dual formulations \eqref{eq:strong_primalPH} and \eqref{eq:strong_dualPH} is obtained by interpreting the codifferential weakly via integration by parts.
One can notice that these two formulations represent a simplified Hodge wave Laplacian \cite{wu2021hodgewave}, one where either the term $\d\d^*$ or $\d^*\d$ is dropped. One of the two variables in the mixed formulation undergoes differentiation and is chosen in a suitable Sobolev space. The other variable is simply taken to be a square integrable function.
\subsection{Weak primal and dual port-Hamiltonian system} 
In this section we introduce the weak formulation to highlight the role played by the two variables. To accomodate for the essential boundary conditions, the test functions are taken in Sobolev spaces with homogeneous boundary conditions
\begin{equation*}
H_\Gamma\Omega^k(M) := \{{\omega}^{k} \in H\Omega^{k}(M) \vert \; \tr \omega^{k}\vert_{\Gamma} = 0 \}, \qquad \Gamma \subset \partial M.
\end{equation*}
In boundary controlled systems (port-Hamiltonian systems belong to this category), the data for the boundary conditions are taken to be free parameters. This is useful for control or modelling purposes (for instance to interconnect different systems together). In this work the space for the controls (i.e. boundary conditions) is taken to be the trace of Sobolev spaces, leading to fractional Sobolev space
\begin{equation*}
\begin{aligned}
    U_1 &:= H^{-1/2}\Omega^{q-1}(\Gamma_1) =  \{\tr {\omega}^k|_{\Gamma_1}, \; {\omega}^{q-1} \in H{\Omega}^{q-1}(M)\}, \\ 
    \dual{U}_2 &:= H^{-1/2}\dual{\Omega}^{p-1}(\Gamma_2) =  \{\tr \dual{\omega}^{p-1}|_{\Gamma_2}, \; \dual{\omega}^{p-1} \in H{\Omega}^{p-1}(M)\}.
\end{aligned}
\end{equation*}
\paragraph{Weak primal system}
The weak formulation for the primal system \eqref{eq:strong_primalPH} reads:\\
find $\dual{\alpha}^p \in L^2\dual{\Omega}^p(M)$,  $\dual{\beta}^{p-1} \in H\dual{\Omega}^{p-1}(M)$ such that $(-1)^p \tr \dual{\beta}^{p-1} \vert_{\Gamma_2} = \dual{u}^{p-1}_2 \in \dual{U}_2$ and
 \begin{subequations}\label{eq:weak_primalPH}
    \begin{align}
    \inprDom[M]{\dual{v}^p}{\partial_t \dual{\alpha}^p} &= (-1)^{p}\inprDom[M]{{v}^p}{\d \dual{\beta}^{p-1}}, \qquad 
    &&\forall \dual{v}^p \in L^2\dual{\Omega}^p(M), \label{eq:weak_primalPH_1}\\
     \inprDom[M]{\dual{v}^{p-1}}{\partial_t \dual{\beta}^{p-1}} &= (-1)^{p} \{- \inprDom[M]{\d\dual{v}^{p-1}}{\dual{\alpha}^p} + \dualprBd[\Gamma_1]{\tr \dual{v}^{p-1}}{u^{q-1}_1}\}, \qquad
    &&\forall \dual{v}^{p-1} \in H_{\Gamma_2}\dual{\Omega}^{p-1}(M). \label{eq:weak_primalPH_2}
    \end{align}
 \end{subequations}
By integration by parts, the control variable $u_1^{q-1} \in U_1$ is naturally included above.
\paragraph{Weak dual system}
For the dual system \eqref{eq:strong_dualPH}, the following weak formulation is obtained: find ${\alpha}^{q-1} \in H\Omega^{q-1}(M), \; {\beta}^{q}_2 \in L^2\Omega^q(M)$ such that $\tr {\alpha}^{q-1} \vert_{\Gamma_1} = u_1^{q-1} \in U_1, $ and
 \begin{subequations}\label{eq:weak_dualPH}
    \begin{align}
    \inprDom[M]{{v}^{q-1}}{\partial_t {\alpha}^{q-1}} &= \inprDom[M]{\d {v}^{q-1}}{{\beta}^{q}} +(-1)^{(p-1)(q-1)}\dualprBd[\Gamma_2]{\tr v^{q-1}}{\dual{u}_2^{p-1}}, \quad &&\forall {v}^{q-1} \in H_{\Gamma_1}\Omega^{q-1}(M), \\
     \inprDom[M]{{v}^{q}}{\partial_t {\beta}^{q}} &= -\inprDom[M]{v^{q}}{\d{\alpha}^{q-1}}, \quad &&\forall {v}^{q} \in L^2\Omega^q(M), 
    \end{align}
 \end{subequations}
Again, the control input $\dual{u}^{p-1}_2 \in U_2$ is naturally included in this formulation.
}

\revtwo{
\subsubsection{Dual field mixed conforming Galerkin discretization}
\revone{A subcomplex of de Rham complex is a sequence of subspaces $S^k \subset H\Omega^k(M)$ such that $\d{S}^k \subset S^{k+1}$.}
Conforming finite element differential forms that constitute a subcomplex of the de Rham complex are denoted by $V_h^k \subset H\Omega^k(M)$ (for instance the trimmed polynomial family $V_h^k = \mathcal{P}^-\Omega^k(\mathcal{T}_h)$ \cite{arnold2006acta} or the mimetic polynomial spaces \cite{kreeft2011mimetic,palha2014}). To accomodate for the essential boundary conditions, the test functions are taken in the corresponding discrete space with boundary conditions
\begin{equation*}
{V}^{k}_h(\Gamma) := \{\dual{\omega}^{k}_h \in {V}^{k}_{h} \vert \; \tr {\omega}^{k}_h\vert_{\Gamma} = 0 \}, \qquad \Gamma \subset \partial M.
\end{equation*}
In boundary controlled systems (port-Hamiltonian systems belong to this category \cite{vanderSchaft2002}), the boundary conditions data are taken to be free parameters one can act upon for control or modelling purposes. In this work the space for the controls (i.e. boundary conditions) is taken to be the restriction at the boundary of the $V_h^k$ spaces
\begin{equation}
    U_{1, h} = \tr V_h^{q-1}|_{\Gamma_1}, \qquad \dual{U}_{2, h} = \tr \dual{V}_h^{p-1}|_{\Gamma_2}.
\end{equation}
The following conforming discretization is the same as in \cite{brugnoli2022df}.
\paragraph{Discrete primal conforming system}
The discrete formulation for the primal system reads:\\
find $\dual{\alpha}^p_{h} \in \dual{V}^p_h, \; \dual{\beta}^{p-1}_{h} \in \dual{V}^{p-1}_h$ such that $(-1)^p \tr \dual{\beta}^{p-1}_{h} \vert_{\Gamma_2} = \dual{u}^{p-1}_{2, h} \in U_{2, h}$ and
 \begin{subequations}\label{eq:discrete_primal_conformingPH}
    \begin{align}
    \inprDom[M]{\dual{v}^p_h}{\partial_t\dual{\alpha}^p_{h}} &= (-1)^{p}\inprDom[M]{\dual{v}^p_h}{\d \dual{\beta}^{p-1}_{h}}, \qquad &&\forall \dual{v}^p_h \in \dual{V}_h^p, \label{eq:discrete_primal_conformingPH_1}\\
     \inprDom[M]{\dual{v}^{p-1}_h}{\partial_t \dual{\beta}^{p-1}_{h}} &= (-1)^{p} \{- \inprDom[M]{\d\dual{v}_h^{p-1}}{\dual{\alpha}^p_{h}} + \dualprBd[\Gamma_1]{\tr \dual{v}_h^{p-1}}{u^{q-1}_{1, h}}\} \qquad &&\forall \dual{v}^{p-1} \in \dual{V}_h^{p-1}(\Gamma_2), \label{eq:discrete_primal_conformingPH_2}
    \end{align}  
 \end{subequations}
where $u_{1, h}^{q-1} \in U_{1, h}$. 

\paragraph{Discrete dual conforming system}
The discrete dual port-Hamiltonian system is given by:\\
find ${\alpha}^{q-1}_{h} \in V^{q-1}_{h}, \; {\beta}^{q}_{h} \in V^q_h$ such that $\tr {\alpha}^{q-1}_{h} \vert_{\Gamma_1} = {u}^{q-1}_{1, h} \in U_{1, h}$ and
\begin{subequations}\label{eq:discrete_dual_conformingPH}
	\begin{align}
    	\inprDom[M]{v^{q-1}_h}{\partial_t {\alpha}^{q-1}_{h}} &= \inprDom[M]{\d{v}^{q-1}_h}{{\beta}^{q}_{h}} + (-1)^{(p-1)(q-1)} \dualprBd[\Gamma_2]{\tr{v}^{q-1}_h}{\dual{u}_{2, h}^{p-1}}, \qquad &&\forall v^{q-1}_h \in {V}^{q-1}_{h}(\Gamma_1), \label{eq:discrete_dual_conformingPH_1}\\
		\inprDom[M]{v^{q}_h}{\partial_t {\beta}^{q}_{h}} &= -\inprDom[M]{v^q_h}{\d {\alpha}^{q-1}_{h}}, \qquad &&\forall {v}^q_h \in V^{q}_{h}, \label{eq:discrete_dual_conformingPH_2}
	\end{align} 
\end{subequations}
where $\dual{u}_{2, h}^{p-1} \in \dual{U}_{2, h}$. 

\begin{remark}[Equivalence with the second order formulation]
Since $\d{}\, V_h^{p-1} \subset V_h^{p}$ and $\d{}\, V_h^{q-1} \subset V_h^{q}$, Eqs. \eqref{eq:discrete_primalPH_1} \eqref{eq:discrete_dualPH_2} hold in a strong sense. The mixed formulations coincides therefore with the second order formulation in time and space \cite{joly2003}.
\end{remark}

\subsubsection{Conforming discretization where the continuity on the $L^2$ variable is relaxed}\label{sec:conGal_dis}

From the weak formulations \eqref{eq:weak_primalPH}, \eqref{eq:weak_dualPH} the variables $\dual{\alpha}^p \in L^2\dual{\Omega}^p(M), \; \beta^{q} \in L^2\Omega^q(M)$ do not undergo any differential and can therefore be discretized by a discontinuous finite element space. In this work, this space is chosen to be the broken version of a finite element space for differential forms.  \\

 For each $T \in \mathcal{T}_h$, let $W_h^k(T) \subset H\Omega^k(T)$ be a finite subcomplex and 
\begin{equation}\label{eq:Wh_space}
    W_h^k := \prod_{T \in \mathcal{T}_h} W_h^k(T).
\end{equation}
Conforming spaces are naturally included in their broken counterpart, i.e. $V^k_h \hookrightarrow W_h^k$. If are the variables belonging to $L^2\Omega^k(M)$ (i.e. ${\alpha}^{p}_{h}, \; {\beta}^{q}_{h}$) are chosen to be in $W_h^k$ then 
$$\d V_h^{k-1} \subset V^k_h \hookrightarrow W_h^k$$ 
and one equation in the mixed formulation will still hold strongly. This property is important and leads to an equivalent formulation with respect to the conforming one presented in \cite{brugnoli2022df}.
}

\paragraph{Discrete primal system}
The discrete formulation for the primal system reads:\\
find $\dual{\alpha}^p_{h} \in \dual{W}^p_h, \; \dual{\beta}^{p-1}_{h} \in \dual{V}^{p-1}_h$ such that $(-1)^p \tr \dual{\beta}^{p-1}_{h} \vert_{\Gamma_2} = \dual{u}^{p-1}_{2, h} \in U_{2, h}$ and
 \begin{subequations}\label{eq:discrete_primalPH}
    \begin{align}
    \inprDom[M]{\dual{v}^p_h}{\partial_t\dual{\alpha}^p_{h}} &= (-1)^{p}\inprDom[M]{\dual{v}^p_h}{\d \dual{\beta}^{p-1}_{h}}, \qquad &&\forall \dual{v}^p_h \in \dual{W}_h^p, \label{eq:discrete_primalPH_1}\\
     \inprDom[M]{\dual{v}^{p-1}_h}{\partial_t \dual{\beta}^{p-1}_{h}} &= (-1)^{p} \{- \inprDom[M]{\d\dual{v}_h^{p-1}}{\dual{\alpha}^p_{h}} + \dualprBd[\Gamma_1]{\tr \dual{v}_h^{p-1}}{u^{q-1}_{1, h}}\} \qquad &&\forall \dual{v}^{p-1} \in \dual{V}_h^{p-1}(\Gamma_2), \label{eq:discrete_primalPH_2}
    \end{align}  
 \end{subequations}
where $u_{1, h}^{q-1} \in U_{1, h}$. 

\paragraph{Discrete dual system}
The discrete dual port-Hamiltonian system is given by:\\
find ${\alpha}^{q-1}_{h} \in V^{q-1}_{h}, \; {\beta}^{q}_{h} \in W^q_h$ such that $\tr {\alpha}^{q-1}_{h} \vert_{\Gamma_1} = {u}^{q-1}_{1, h} \in U_{1, h}$ and
\begin{subequations}\label{eq:discrete_dualPH}
	\begin{align}
    	\inprDom[M]{v^{q-1}_h}{\partial_t {\alpha}^{q-1}_{h}} &= \inprDom[M]{\d{v}^{q-1}_h}{{\beta}^{q}_{h}} + (-1)^{(p-1)(q-1)} \dualprBd[\Gamma_2]{\tr{v}^{q-1}_h}{\dual{u}_{2, h}^{p-1}}, \qquad &&\forall v^{q-1}_h \in {V}^{q-1}_{h}(\Gamma_1), \label{eq:discrete_dualPH_1}\\
		\inprDom[M]{v^{q}_h}{\partial_t {\beta}^{q}_{h}} &= -\inprDom[M]{v^q_h}{\d {\alpha}^{q-1}_{h}}, \qquad &&\forall {v}^q_h \in W^{q}_{h}, \label{eq:discrete_dualPH_2}
	\end{align} 
\end{subequations}
where $\dual{u}_{2, h}^{p-1} \in \dual{U}_{2, h}$. 

\revtwo{
\begin{remark}[Equivalence of the two conforming formulations]
    Since $\d V_h^{k-1} \subset V^k_h \hookrightarrow W_h^k$, applying the results of \cite{joly2003} one obtains the equivalence of \eqref{eq:discrete_primal_conformingPH} and \eqref{eq:discrete_primalPH}  and \eqref{eq:discrete_dual_conformingPH} and \eqref{eq:discrete_dualPH}.
\end{remark}
}

\subsection{Discrete power balance}\label{sec:discrete_powbal}

The dual field mixed Galerkin discretization is capable of retaining the following discrete power balance (that characterizes the Dirac structure underlying port-Hamiltonian systems).

\begin{proposition}\label{pr:powbal}
The primal-dual discrete port-Hamiltonian systems \eqref{eq:discrete_primalPH}, \eqref{eq:discrete_dualPH} encode the following discrete power balance
\begin{equation*}
    \begin{aligned}
    (-1)^{p(n-p)}\dualprDom[M]{\alpha^{q-1}_h}{\partial_t \dual{\alpha}^p_h} + \dualprDom[M]{\dual{\beta}^{p-1}_h}{\partial_t \beta^q_h} = \dualprBd[\Gamma_1]{\dual{y}^{p-1}_{2, h}}{u^{q-1}_{1, h}} + \dualprBd[\Gamma_2]{\dual{u}^{p-1}_{2, h}}{y^{q-1}_{1, h}}
    \end{aligned}
\end{equation*}
\end{proposition}

\begin{proof}
The proof generalizes Pr. 5 in \cite{brugnoli2022df}. Since the conforming discrete spaces $V_h^k$ form a de Rham subcomplex and conforming spaces are naturally included in their broken counterpart, two pointwise discrete conservation laws are verified, namely \eqref{eq:discrete_primalPH_1}, \eqref{eq:discrete_dualPH_2}
\begin{equation*}
    \begin{aligned}
    &\d \dual{V}^{p-1}_h \subset \dual{V}^p_h \hookrightarrow \dual{W}_h^p, \\ 
    &\d V^{q-1}_h \subset V^q_h \hookrightarrow W^q_h,
    \end{aligned} \quad \implies \quad 
    \begin{aligned}
     \partial_t \dual{\alpha}_h^p &= (-1)^p \d \dual{\beta}_h^{p-1}, \\
    \partial_t {\beta}_h^q &= -\d {\alpha}_h^{q-1}.
    \end{aligned}
\end{equation*}

Taking the duality product of the first equation with $(-1)^{p(n-p)}\alpha_h^{q-1}$ and with $\dual{\beta}^{p-1}$ for the second and summing over each cell leads 
\begin{equation*}
    \begin{aligned}
    (-1)^{p(n-p)}\dualprDom[M]{\alpha^{q-1}_h}{\partial_t \dual{\alpha}^p_h} + \dualprDom[M]{\dual{\beta}^{p-1}_h}{\partial_t \beta^q_h} &= \sum_{T\in \mathcal{T}_h} (-1)^p\dualprDom[T]{\d \dual{\beta}_h^{p-1}}{\alpha^{q-1}_h} - \dualprDom[T]{\dual{\beta}^{p-1}_h}{\d {\alpha}_h^{q-1}},\\
    &= \sum_{T\in \mathcal{T}_h} (-1)^p\dualprBd[\partial T]{\tr \dual{\beta}_h^{p-1}}{\tr \alpha^{q-1}_h},\\
    &=\dualprBd[\Gamma_1]{\dual{y}^{p-1}_{2, h}}{u^{q-1}_{1, h}} + \dualprBd[\Gamma_2]{\dual{u}^{p-1}_{2, h}}{y^{q-1}_{1, h}}.
    \end{aligned}
\end{equation*}
The final equality uses the continuity property of differential forms (cf. Appendix A in \cite{brugnoli2022df}). 
\end{proof}

\section{Hybridization of mixed Galerkin schemes}\label{sec:hybridization}

Let $W_h^k \subset H\Omega^k(\mathcal{T}_h)$ be defined as in \eqref{eq:Wh_space}. Conforming finite element spaces for the Sobolev spaces are given by $V_h^k = V^k \cap W_h^k$ where $V^k = H\Omega^k(M)$. We recall from \cite{stern2023hyb} the broken and unbroken space of tangential traces finite element differential forms  are then defined as follows
\begin{equation*}
\begin{aligned}
    W_h^{k, \bm{t}} &:= \{\tr \omega_h^k: \omega_h^k \in W_h^k\}, \\
    V_h^{k, \bm{t}} &:= \{\tr \omega_h^k: \omega_h^k \in V_h^k\} = V^{k, \bm{t}} \cap W_h^{k, \bm{t}}, \qquad \text{where} \qquad V^{k, \bm{t}} := H\Omega^{k, \bm{t}}(\partial \mathcal{T}_h), \\
    V_h^{k, \bm{t}}(\Gamma_1) &:= \{\tr \omega_h^k: \omega_h^k \in V_h^k(\Gamma_1)\} = V^{k, \bm{t}}(\Gamma_1) \cap W_h^{k, \bm{t}}.
\end{aligned}
\end{equation*}
For what concerns the space $W^{k, \bm{n}}$, its discrete version is taken to be the dual space of $W_h^{k, \bm{t}}$, i.e. $W^{k, \bm{n}}_h = (W_h^{k, \bm{t}})^*$. Since discrete tangential traces are piecewise polynomial, they are in $L^2(\partial\mathcal{T}_h)$ and so $W^{k, \bm{n}}_h = W_h^{k, \bm{t}}$ where $\inprBd[\partial\mathcal{T}_h]{\cdot}{\cdot}$ is the $L^2$ inner product.

\paragraph{Primal discrete system}
Consider the variational
problem: find 
\begin{itemize}
    \item Local variables: $\dual{\alpha}^p_{h} \in \dual{W}^p_h, \quad \dual{\beta}^{p-1}_{h} \in \dual{W}^{p-1}_h, \quad \dual{\alpha}^{p-1, \bm{n}}_{h} \in \dual{W}^{p-1, \bm{n}}_h$;
    \item Global variables $\dual{\beta}^{p-1, \bm{t}}_{h} \in \dual{V}^{p-1, \bm{t}}_h$;
\end{itemize}
such that $\dual{\beta}^{p-1, \bm{t}}_{h}|_{\Gamma_2} = \dual{u}^{p-1}_{2, h}$ and 

\begin{subequations}\label{eq:discrete_primalPHhybrid}
    \begin{align}
    \inprDom[\mathcal{T}_h]{\dual{v}^p_h}{\partial_t \dual{\alpha}^p_{h}} &= (-1)^{p}\inprDom[\mathcal{T}_h]{{v}^p_h}{\d \dual{\beta}^{p-1}_{h}}, \quad &&\forall \dual{v}^p_h \in \dual{W}^p_h, \label{eq:discrete_primalPHhybrid_1}\\
     \inprDom[\mathcal{T}_h]{\dual{v}^{p-1}_h}{\partial_t \dual{\beta}^{p-1}_{h}} &= (-1)^{p}  \{-\inprDom[\mathcal{T}_h]{\d\dual{v}^{p-1}_h}{\dual{\alpha}^p_{h}} + \inprBd[\partial\mathcal{T}_h]{\tr \dual{v}^{p-1}_h}{\dual{\alpha}^{p-1, \bm{n}}_{h}}\}, \quad &&\forall \dual{v}^{p-1}_h \in \dual{W}^{p-1}_h, \label{eq:discrete_primalPHhybrid_2}\\
    0 &= -(-1)^p \inprBd[\partial\mathcal{T}_h]{\dual{v}^{p-1, \bm{n}}_h}{\tr \dual{\beta}^{p-1}_{h} - \dual{\beta}^{p-1, \bm{t}}_{h}}, \quad &&\forall \dual{v}^{p-1, \bm{n}}_h \in \dual{W}^{p-1, \bm{n}}_h, \label{eq:discrete_primalPHhybrid_3}\\
    0 &= (-1)^p \{-\inprBd[\partial\mathcal{T}_h]{\dual{v}^{p-1, \bm{t}}_h}{\dual{\alpha}^{p-1, \bm{n}}_{h}} + \dualprBd[\Gamma_1]{\dual{v}^{p-1, \bm{t}}_h}{u_{1, h}^{q-1}}\}, \quad &&\forall \dual{v}^{p-1, \bm{t}}_h \in \dual{V}^{p-1, \bm{t}}_h(\Gamma_2). \label{eq:discrete_primalPHhybrid_4}
    \end{align}
 \end{subequations}

\paragraph{Dual discrete system}
 Consider the variational
problem: find 
\begin{itemize}
    \item Local variables: ${\alpha}^{q-1}_{h} \in {W}^{q-1}_h, \quad {\beta}^{q}_{h} \in {W}^{q}_h, \quad {\beta}^{q-1, \bm{n}}_{h} \in {W}^{q-1, \bm{n}}_h$;
    \item Global variables ${\alpha}^{q-1, \bm{t}}_{h} \in {V}^{q-1, \bm{t}}_h$;
\end{itemize}
such that ${\alpha}^{q-1, \bm{t}}_{h}|_{\Gamma_1} = u^{q-1}_{1, h}$ and 
\begin{subequations}\label{eq:discrete_dualPHhybrid}
    \begin{align}
    \inprDom[\mathcal{T}_h]{{v}^{q-1}_h}{\partial_t {\alpha}^{q-1}_{h}} &= \inprDom[\mathcal{T}_h]{\d {v}^{q-1}_h}{\beta^{q}_{h}} - \inprBd[\partial\mathcal{T}_h]{\tr {v}^{q-1}_h}{\beta^{q-1, \bm{n}}_{h}}, \quad &&\forall {v}^{q-1}_h \in W^{q-1}_h, \label{eq:discrete_dualPHhybrid_1}\\
     \inprDom[\mathcal{T}_h]{{v}^{q}_h}{\partial_t \beta^{q}_{h}} &= -\inprDom[\mathcal{T}_h]{v^{q}_h}{\d{\alpha}^{q-1}_{h}}, \quad &&\forall {v}^{q}_h \in W^{q}_h, \\
     0 &= \inprBd[\partial\mathcal{T}_h]{v^{q-1, \bm{n}}_h}{\tr {\alpha}^{q-1}_h - {\alpha}^{q-1, \bm{t}}_{h}}, \quad &&\forall v^{q-1, \bm{n}}_h \in W^{q-1, \bm{n}}_h, \\
     0 &= \inprBd[\partial\mathcal{T}_h]{v^{q-1, \bm{t}}_h}{\beta^{q-1, \bm{n}}_{h}} + (-1)^{(p-1)(q-1)}  \dualprBd[\Gamma_2]{{v}^{q-1, \bm{t}}_h}{\dual{u}_{2, h}^{p-1}} ,  &&\forall v^{q-1, \bm{t}}_h \in V^{q-1, \bm{t}}_h(\Gamma_1).
    \end{align}
 \end{subequations}

 \begin{figure}[h]
    \centering
    \begin{subfigure}[t]{0.3\textwidth}
        \includegraphics[width=\textwidth]{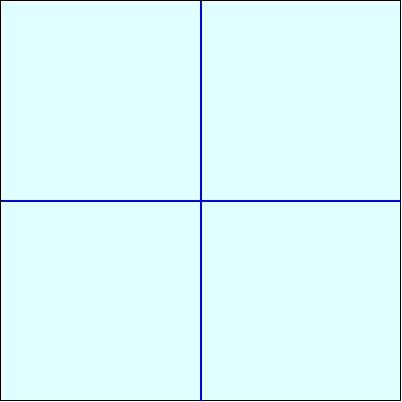}
        \caption{Computational mesh $\mathcal{T}_h$. The internal facets are highlighted in blue.}
        \label{fig:mesh}
    \end{subfigure}
    \hfill
    \begin{subfigure}[t]{0.3\textwidth}
        \includegraphics[width=\textwidth]{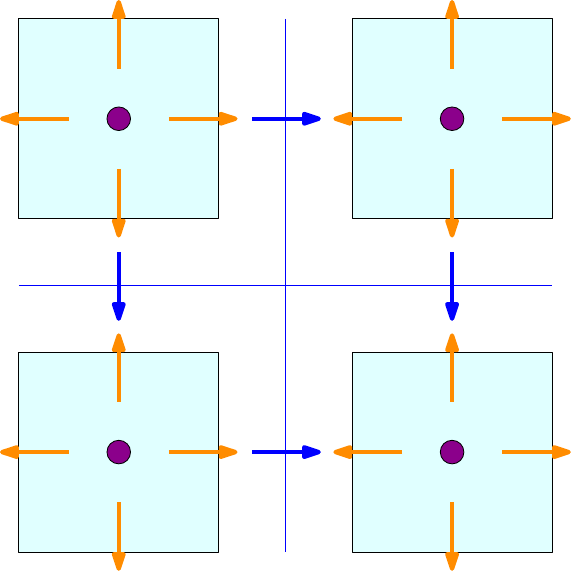}
        \caption{Primal system: the regular variable is an outer oriented 1-form. The globally coupled Lagrange multiplier enforces normal continuity.}
        \label{fig:image2}
    \end{subfigure}
    \hfill
    \begin{subfigure}[t]{0.3\textwidth}
        \includegraphics[width=\textwidth]{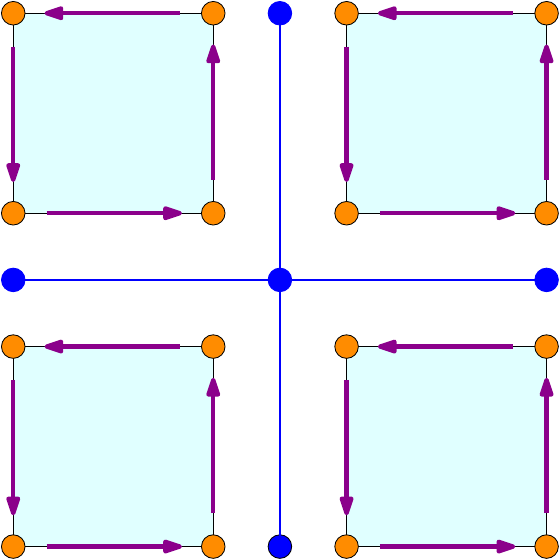}
        \caption{Dual system: the regular variable is a 0-form. The globally coupled Lagrange multiplier enforces continuity.}
        \label{fig:image3}
    \end{subfigure}
    \caption{\revone{Illustration of global and local variables for $n=2, p=2, q=1$ for the lowest order polynomial mesh. The global variable lives on the mesh's skeleton and its degrees of freedom are in blue \textcolor{blue}{\rule{.7em}{.7em}}. The degrees of freedom of the local variable that undergoes hybridization is depicted in orange \textcolor{darkorange}{\rule{.7em}{.7em}} and those of the variable that is discontinuous and does not require regularity is depicted in purple \textcolor{darkmagenta}{\rule{.7em}{.7em}}.}  }
    \label{fig:dofs_p2_q1}
\end{figure}

\revone{Fig. \ref{fig:dofs_p2_q1} depicts the degrees of freedom for the different variables for a 2 dimensional domain and $p=2, \; q=1$. In this case the primal system contains a discontinuous 2-form (discretized using a discontinuous Galerkin space) and a continuous 1-form (discretized using a Raviart-Thomas space). This latter variable undergoes hybridization, i.e. its continuity is first broken and reestablished using a Lagrange multiplier that enforces normal continuity (an outer oriented 1-form $\dual{\beta}^1 \in H\dual{\Omega}^1(M)$ in a two dimensional domain is normal continuous). The dual contains a discontinuous one-form (discretized using a broken N\'ed\'elec space) and a continuous 0-form (discretized using Lagrange finite elements). This latter variable is hybridized. Its continuity is broken and a Lagrange multiplier enforces continuity (a 0-form ${\alpha}^0 \in H\dual{\Omega}^0(M)$ is continuous).}

\paragraph{Equivalence of the discrete mixed and hybrid formulation}
The mixed and hybrid primal discrete formulations are equivalent. The only difference is that $\dual{\alpha}^{p-1, \bm{n}}_{h}$ (resp. ${\beta}^{q-1, \bm{n}}_{h}$) are not exactly equal to the normal trace of $\dual{\alpha}^p_{h}$ (resp. ${\beta}^q_{h}$), but only weakly \cite{stern2023hyb}.
\begin{theorem}\label{thm:equiv_discrete}
For the primal system, the following are equivalent:
\begin{itemize}
    \item $\dual{\alpha}^p_{h}, \dual{\beta}^{p-1}_{h}, \dual{\alpha}^{p-1, \bm{n}}_{h}, \dual{\beta}^{p-1, \bm{t}}_{h}$ is a solution to \eqref{eq:discrete_primalPHhybrid};
    \item $\dual{\alpha}^p_{h}, \dual{\beta}^{p-1}_{h}$ is a solution to \eqref{eq:discrete_primalPH}. Moreover, $\dual{\beta}^{p-1, \bm{t}}_{h} = \tr \dual{\beta}^{p-1}_{h}$ and $\dual{\alpha}^{p-1, \bm{n}}_{h}$ is uniquely determined by \eqref{eq:discrete_primalPHhybrid_2}.
\end{itemize}
For the dual system, the following are equivalent:
\begin{itemize}
    \item ${\alpha}^{q-1}_{h}, \; {\beta}^{q}_{h}, \; {\beta}^{q-1, \bm{n}}_{h}, \; {\alpha}^{q-1, \bm{t}}_{h}$ is a solution to \eqref{eq:discrete_dualPHhybrid};
    \item ${\alpha}^{q-1}_{h}, \; {\beta}^{q}_{h}$ is a solution to \eqref{eq:discrete_dualPH}. Moreover, ${\alpha}^{q-1, \bm{t}}_{h} = \tr {\alpha}^{q-1}_{h}$ and ${\beta}^{q-1, \bm{n}}_{h}$ is uniquely determined by \eqref{eq:discrete_dualPHhybrid_1}.
\end{itemize}
\end{theorem}
\begin{proof}
The proof will be given only for the primal system. The one for the dual system is completely analogous. Assume a solution to \eqref{eq:discrete_primalPHhybrid} is given. Then from the variational formulation, it follows that $\dual{\beta}_h^{p-1, \bm{t}} = \tr \dual{\beta}_h^{p-1}$. By Pr. \ref{pr:H_br_uniform} this means that $\dual{\beta}_h^{p-1} \in \dual{H}\Omega^{p-1}(M)$. This immediately implies \eqref{eq:discrete_primalPH_1}. Taking the test function $\dual{v}^{p-1} \in H\dual{\Omega}^{p-1}(M, \Gamma_2)$ in \eqref{eq:discrete_primalPHhybrid_2}, the boundary term on $\partial\mathcal{T}_h$ is replaced by the contribution on the $\Gamma_1$ boundary, thus leading to \eqref{eq:discrete_primalPH_2}. Uniqueness of $\dual{\alpha}^{p-1, \bm{n}}_{h}$  follows from the fact that the broken space of normal traces is in duality with the broken space of tangential traces \revtwo{(actually they are taken to be the same space). Therefore the resulting linear differential algebraic system is regular. Conversely assume that a solution $(\dual{\alpha}^p_{h}, \dual{\beta}^{p-1}_{h})$ to \eqref{eq:discrete_primalPH} is given. Then because of the fact that $\inprBd[\partial\mathcal{T}_h]{\cdot}{\cdot}$ is a duality pairing, equations \eqref{eq:discrete_primalPHhybrid_1}, \eqref{eq:discrete_primalPHhybrid_2}, \eqref{eq:discrete_primalPHhybrid_3} hold. For the last equation, combining \eqref{eq:discrete_primalPH_2} and \eqref{eq:discrete_primalPHhybrid_2} gives 
$$\inprBd[\partial\mathcal{T}_h]{\tr \dual{v}^{p-1}_h}{\dual{\alpha}^{p-1, \bm{n}}_{h}} = \dualprBd[\Gamma_1]{\tr \dual{v}^{p-1}_h}{u_{1, h}^{q-1}}, \qquad \text{for} \quad \dual{v}^{p-1}_h \in V_h^{p-1}(\Gamma_2).$$
This implies  \eqref{eq:discrete_primalPHhybrid_4}.
}
\end{proof}

\begin{remark}
Since the primal and dual discrete hybrid systems are equivalent to the mixed Galerkin formulations, they satisfy a discrete power balance.    
\end{remark}

\section{Algebraic realization of the discrete systems and static condensation}\label{sec:realization}
Given a finite element basis, the algebraic realization of the different terms in the weak formulation can be computed. Both the primal discrete system \eqref{eq:discrete_primalPHhybrid} and dual discrete system \eqref{eq:discrete_dualPHhybrid} have the following common structure
\begin{equation}\label{eq:alg_sys_common}
    \begin{bmatrix}
    \mathbf{E}_{l} & \mathbf{0} \\
    \mathbf{0} & \mathbf{0} \\
    \end{bmatrix}
    \begin{pmatrix}
    \dot{\mathbf{x}}_{l} \\
    \dot{\mathbf{x}}_{g} \\
    \end{pmatrix} = \begin{bmatrix}
    \mathbf{J}_l &  \mathbf{C}_{lg} \\
    -\mathbf{C}_{lg}^\top & \mathbf{0} \\
    \end{bmatrix}
    \begin{pmatrix}
    {\mathbf{x}}_{l} \\
    {\mathbf{x}}_{g} \\
    \end{pmatrix} + 
    \begin{bmatrix}
    \mathbf{B}_l &  \mathbf{0} \\
    \mathbf{0} & \mathbf{B}_g \\
    \end{bmatrix}
    \begin{pmatrix}
    {\mathbf{u}}_{l} \\
    {\mathbf{u}}_{g} \\
    \end{pmatrix},
\end{equation}
where the subscript $l$ denotes the local variables and the subscript $g$ denotes the global variable. Matrix 
$\mathbf{E}_l$ is symmetric and positive semi-definite, while $\mathbf{J}_l$ is skew-symmetric $\mathbf{J}_l = - \mathbf{J}_l^{\top}$. This structure is a particular instance of a port-Hamiltonian descriptor system of index 2~\cite{beattie2018phdae}. 
We present the detailed matrix expression of the primal and dual algebraic systems. We shall use in what follows the notation $[\mathbf{x}]_{R}$ to indicate the vector containing the rows of vector $\mathbf{x}$ associated with the index set $R$ only.
In particular, we will denote the index sets corresponding to the degrees of freedom of the interior of the domain by $I$, and the ones associated to the subpartitions ${\Gamma}_1$ and ${\Gamma}_2$ by the same letters $\Gamma_1,\, \Gamma_2$ (with a slight abuse of notation, since these symbols are also used for the corresponding continuous states). Local matrices have as subscript the disjoint union of cells or boundaries (i.e. $\mathcal{T}_h$ or $\partial \mathcal{T}_h$). Coupled matrices have a subscipt denoting the of facets $\mathcal{F}_h = \cup_{T \in \mathcal{T}_h} \partial T$. \revone{In Appendix A the matrix realization of the local and global weak
formulations is described. In particular the notation used just after is explained.}

\paragraph{The primal system}
For what concerns the primal system, the state and input variables are the following
\begin{equation*}
    \mathbf{x}_{l} = \begin{pmatrix}\dual{\bm{\alpha}}^p \\
    {\dual{\bm{\beta}}}^{p-1} \\
    {\dual{\bm{\alpha}}}^{p-1, \bm{n}}\end{pmatrix}, \qquad \mathbf{x}_{g} = [{\dual{\bm{\beta}}}^{p-1, \bm{t}}]_{I \cup \Gamma_1}, \qquad \begin{pmatrix}
    {\mathbf{u}}_{l} \\
    {\mathbf{u}}_{g} \\
    \end{pmatrix} = \begin{pmatrix}
    \dual{\mathbf{u}}_{2}^{p-1} \\
    {\mathbf{u}}_{1}^{q-1} \\
    \end{pmatrix}.
\end{equation*}
The specific structure of the matrices appearing in \eqref{eq:alg_sys_common} is given by

\begin{align*}
\mathbf{E}_l  &=
\begin{bmatrix}
\mathbf{M}^p_{\mathcal{T}_h} & \mathbf{0} & \mathbf{0}\\
\mathbf{0} & \mathbf{M}^{p-1}_{\mathcal{T}_h} & \mathbf{0} \\
\mathbf{0} & \mathbf{0} & \mathbf{0} \\
\end{bmatrix}, \qquad
\mathbf{J}_l = (-1)^p \begin{bmatrix}
\mathbf{0} & \mathbf{D}^{p-1}_{\mathcal{T}_h} & \mathbf{0} \\
-(\mathbf{D}^{p-1}_{\mathcal{T}_h})^\top & \mathbf{0} & (\mathbf{T}_{\partial\mathcal{T}_h}^{p-1})^\top \mathbf{M}_{\partial\mathcal{T}_h}^{p-1} \\
\mathbf{0} & -\mathbf{M}_{\partial\mathcal{T}_h}^{p-1}\mathbf{T}_{\partial\mathcal{T}_h}^{p-1} & \mathbf{0} \\
\end{bmatrix}, \\
\mathbf{C}_{lg} &= (-1)^p \begin{bmatrix}
\mathbf{0} \\
\mathbf{0} \\
 (\dual{\mathbf{\Xi}}^{p-1}_{\mathcal{F}_h/\Gamma_2})^\top \mathbf{M}_{\mathcal{F}_h/\Gamma_2}^{p-1}
\end{bmatrix}, \qquad
\mathbf{B}_l = (-1)^p 
\begin{bmatrix}
\mathbf{0}\\
\mathbf{0}\\
 ({\mathbf{T}}^{p-1}_{\Gamma_2})^\top \mathbf{M}_{\Gamma_2}^{p-1}
\end{bmatrix}, \\
\mathbf{B}_g &= (-1)^p (\mathbf{T}^{p-1, \bm{t}}_{\Gamma_1})^\top \mathbf{\Psi}_{\Gamma_1}^{q-1}.
\end{align*}

\paragraph{The dual system}
For the dual system, the state and input variables are the following
\begin{equation*}
    \mathbf{x}_{l} = \begin{pmatrix}{\bm{\alpha}}^{q-1} \\
    {{\bm{\beta}}}^{q} \\
    {{\bm{\beta}}}^{q-1, \bm{n}}\end{pmatrix}, \qquad \mathbf{x}_{g} = [{{\bm{\alpha}}}^{q-1, \bm{t}}]_{I \cup \Gamma_2}, \qquad \begin{pmatrix}
    {\mathbf{u}}_{l} \\
    {\mathbf{u}}_{g} \\
    \end{pmatrix} = \begin{pmatrix}
    {\mathbf{u}}_{1}^{q-1} \\
    \dual{\mathbf{u}}_{2}^{p-1} \\
    \end{pmatrix}.
\end{equation*}
The matrices exhibit the following structure
\begin{align*}
\mathbf{E}_l  &=
\begin{bmatrix}
\mathbf{M}^{q-1}_{\mathcal{T}_h} & \mathbf{0} & \mathbf{0}\\
\mathbf{0} & \mathbf{M}^{q}_{\mathcal{T}_h} & \mathbf{0} \\
\mathbf{0} & \mathbf{0} & \mathbf{0} \\
\end{bmatrix}, \qquad
\mathbf{J}_l = \begin{bmatrix}
\mathbf{0} & (\mathbf{D}^{q-1}_{\mathcal{T}_h})^\top & -(\mathbf{T}_{\partial\mathcal{T}_h}^{q-1})^\top  \mathbf{M}_{\partial\mathcal{T}_h}^{q-1}\\
-(\mathbf{D}^{q-1}_{\mathcal{T}_h}) & \mathbf{0} & \mathbf{0}  \\
\mathbf{M}_{\partial\mathcal{T}_h}^{q-1}\mathbf{T}_{\partial\mathcal{T}_h}^{q-1} & \mathbf{0} & \mathbf{0} \\
\end{bmatrix}, \\
\mathbf{C}_{lg} &= -\begin{bmatrix}
\mathbf{0} \\
\mathbf{0} \\
(\mathbf{\Xi}^{q-1}_{\mathcal{F}_h/\Gamma_1})^\top \mathbf{M}_{\mathcal{F}_h/\Gamma_1}^{q-1}
\end{bmatrix}, \qquad
\mathbf{B}_l = -
\begin{bmatrix}
\mathbf{0}\\
\mathbf{0}\\
(\mathbf{T}^{p-1}_{\Gamma_1})^\top \mathbf{M}_{\Gamma_1}^{p-1}
\end{bmatrix}, \\
\mathbf{B}_g &= (-1)^{(p-1)(q-1)} (\mathbf{T}^{q-1, \bm{t}}_{\Gamma_2})^\top \mathbf{\Psi}_{\Gamma_2}^{p-1}.
\end{align*}

\subsection{Time discretization and static condensation}
The time discretization leads to a linear saddle point system that can be efficiently solved via static condensation \cite{guyan1965reduction}. Since we are interested in Hamiltonian conservative systems, conservation of energy is of utmost importance. Implicit Runge-Kutta methods based on Gauss Legendre collocation points can be used to this aim \cite{sanzserna1992}. These methods are also the only collocation schemes
that lead to an exact discrete energy balance in the linear case \cite{kotyczka2019discrete}. The implicit midpoint method
is here used to illustrate the time discretization. \\\

Consider a total simulation time $T_{\mathrm{end}}$ and a equidistant splitting given by the time step $\Delta t =
T_{\mathrm{end}}/N_t$, where $N_t$ is the total number of simulation instants. The evaluation of a generic variable
$\mathbf{x}$ at the time instant $t^n = n\Delta t$ is denoted by $\mathbf{x}^n$. Consider once again system \eqref{eq:alg_sys_common}. The application of the implicit midpoint scheme leads to the following algebraic system
\begin{equation*}
    \begin{bmatrix}
    \mathbf{E}_l - \frac{\Delta t}{2}\mathbf{J}_l & -\frac{\Delta t}{2}\mathbf{C}_{lg}\\
    \frac{\Delta t}{2}\mathbf{C}_{lg}^\top & \mathbf{0}
    \end{bmatrix}
    \begin{pmatrix}
    \mathbf{x}_l^{n+1} \\
    \mathbf{x}_g^{n+1} 
    \end{pmatrix} = 
    \begin{bmatrix}
    \mathbf{E}_l + \frac{\Delta t}{2}\mathbf{J}_l & \frac{\Delta t}{2}\mathbf{C}_{lg}\\
    -\frac{\Delta t}{2}\mathbf{C}_{lg}^\top & \mathbf{0}
    \end{bmatrix}
    \begin{pmatrix}
    \mathbf{x}_l^{n} \\
    \mathbf{x}_g^{n} 
    \end{pmatrix} + \Delta t
    \begin{bmatrix}
    \mathbf{B}_l &  \mathbf{0} \\
    \mathbf{0} & \mathbf{B}_g \\
    \end{bmatrix}
    \begin{pmatrix}
    \mathbf{u}_{l}^{n+1/2} \\
    \mathbf{u}_{g}^{n+1/2}
    \end{pmatrix}.
\end{equation*}
This system can be rewritten as the following saddle point problem
\begin{equation}\label{eq:alg_saddlepoint}
\begin{bmatrix}
\mathbf{A} & -\mathbf{C} \\
\mathbf{C}^\top & \mathbf{0}
\end{bmatrix}
\begin{pmatrix}
\mathbf{x}_l \\
\mathbf{x}_g \\
\end{pmatrix} = 
\begin{pmatrix}
    \mathbf{b}_l \\
    \mathbf{b}_g \\
\end{pmatrix}.
\end{equation}

The application of a Schur complement leads to the following system for the global variables
\begin{equation*}
    \mathbf{C}^\top \mathbf{A}^{-1} \mathbf{C} \mathbf{x}_g = \mathbf{b}_g - \mathbf{C}^\top \mathbf{A}^{-1} \mathbf{b}_l.
\end{equation*}
The matrix $\mathbf{A}$ is block diagonal and contains local information. Each block can be inverted numerically or even analytically. Furthermore by Corollary 4.3. in \cite{guducu2022}, the Schur complement $\mathbf{C}^\top \mathbf{A}^{-1} \mathbf{C}$  has a positive semidefinite symmetric part, like the original matrix. This property can be exploited to design efficient preconditioning strategies \cite{guducu2022}.  Once the global variable has been computed, the local unknown is then computed by solving the following system
\begin{equation*}
    \mathbf{x}_l = \mathbf{A}^{-1}\mathbf{b}_l + \mathbf{A}^{-1}\mathbf{C}\mathbf{x}_g .
\end{equation*}
This system is block diagonal and can be solved in parallel. 

\begin{remark}[Well-posedness of the problem]
The matrix $\mathbf{A}$ is invertible as it is given by a saddle point matrix whose off diagonal block are full rank and $\mathbf{C}$ is full rank. These properties are a consequence of the fact that $\inprBd[\partial\mathcal{T}_h]{\cdot}{\cdot}$ is simply an $L^2$ inner product at the discrete level. \revtwo{To be more specific matrix $\mathbf{C}$ is of size $N_l \times N_g$ where $N_l$ is the number of local variable and $N_g$ is the number of global variable. Naturally $N_l<N_g$ and the last block of matrix $\mathbf{C}$ is $\mathbf{\Xi}_{\mathcal{F}_h}^\top \mathbf{M}_{\mathcal{F}_h}$ where matrix $\mathbf{\Xi}_{\mathcal{F}_h}$ is surjective as it represents an average matrix (see Appendix \ref{app:alg_mimetic}) and matrix $\mathbf{M}_{\mathcal{F}_h}$ is a positive definite mass matrix.} Therefore, the algebraic system \eqref{eq:alg_saddlepoint} is well posed.
\end{remark}

A major advantage when using hybridization strategies is the reduced size of the system of equations to be solved. For simplicity consider the case in which no essential boundary conditions are present (this depends on which of the two system, primal or dual, is considered). The  mixed Galerkin discretization leads to a total size of $\dim W_h^k + \dim V_h^{k-1}$ with $k=\{p,q\}$, whereas the hybridization strategy condenses it to $\dim V_h^{k-1, \bm{t}}$. The following result is a simple application of Th. 4.4 in \cite{stern2023hyb}.

\begin{theorem}
When only natural boundary conditions apply, the reduction in size from the conforming discretization to the hybridized system is given by
\begin{equation}
    (\dim W_h^k + \dim V_h^{k-1}) - \dim V_h^{k-1, \bm{t}} = \dim W_h^k + \sum_{T \in \mathcal{T}_h} \dim \mathring{W}_h^{k-1}(T),
\end{equation}
where $\mathring{W}_h(T) :=\{\omega_h \in W_h^{k-1}, \tr \omega_h|_{\partial T} = 0\}$ and $k=\{p,q\}$.
\end{theorem}

\begin{proof}
The space $V_h^{k-1, \bm{t}}$ is the image of $V_h^{k-1}$ under the trace map. By the rank nullity theorem their dimensions differ by the kernel
\begin{equation}
\begin{aligned}
    \dim V_h^{k-1} - \dim V_h^{k-1, \bm{t}} &= \dim\{ \omega_h \in V_h^{k-1}: \tr \omega_h|_{\partial \mathcal{T}_h} = 0 \}, \\
    &= \dim \prod_{T \in \mathcal{T}_h} \mathring{W}_h^{k-1}(T) = \sum_{T \in \mathcal{T}_h} \dim \mathring{W}_h^{k-1}(T).
\end{aligned}
\end{equation}
\end{proof}

\begin{remark}
    A different strategy for solving this problem consist in constructing a projector to eliminate the Lagrange multiplier. This idea is pursed in \cite{park2023partitioned} and offers several computational advantages, in particular for the parallel time domain simulation. 
\end{remark}

\section{Numerical experiments}\label{sec:num_exp}

In this section, the hybridization strategy  is tested for the wave and Maxwell equations. By introducing the musical isomorphism, given by the flat $\flat$ and the sharp operator $\sharp$, and the isomorphism $\Theta$ converting vector fields in $n-1$ forms, the commuting diagram in Fig. \ref{fig:cd_ext_vec}, that provides the link between the de Rham complex and the standard operators and Sobolev space from vector calculus, is obtained.
\begin{figure}[ht]
\centering
\begin{tikzcd}
H\Omega^0(M) \arrow[r, "\d"] \arrow[d, leftrightarrow, "Id"]
& H\Omega^{1}(M) \arrow[r, "\d"] \arrow[d, "\sharp", xshift=10pt] & H\Omega^2(M) \arrow[r, "\d"] \arrow[d, "\Theta^{-1}", xshift=10pt]
& H\Omega^{3}(M) \arrow[d, "\star", xshift=10pt]  \\
H^1(M) \arrow[r, "\grad"]
& H^{\curl}(M) \arrow[r, "\curl"] \arrow[u, "\flat"] & H^{\div}(M) \arrow[r, "\div"] \arrow[u, "\Theta"]
& L^2(M) \arrow[u, "\star^{-1}"]
\end{tikzcd}  
\caption{Equivalence of vector and exterior calculus Sobolev spaces.}
\label{fig:cd_ext_vec}
\end{figure}

For the discretization, the trimmed polynomial family $\mathcal{P}^{-}_s\Omega^k(\mathcal{T}_h)$ is used \revone{(we refer to \cite{arnold2006acta} for the specific definition of this family of finite elements)}. This family corresponds to the well known continuous Galerkin (or Lagrange) elements $\mathcal{P}^-_s\Omega^0(\mathcal{T}_h) \equiv \mathrm{CG}_s(\mathcal{T}_h)$, Nédélec of the first kind $\mathcal{P}^-_s\Omega^1(\mathcal{T}_h) \equiv \mathrm{NED}_s^1(\mathcal{T}_h)$, Raviart-Thomas $\mathcal{P}^-_s\Omega^2(\mathcal{T}_h) \equiv \mathrm{RT}_s(\mathcal{T}_h)$ and discontinuous Galerkin $\mathcal{P}^-_s\Omega^3(\mathcal{T}_h) \equiv \mathrm{DG}_{s-1}(\mathcal{T}_h)$, as illustrated in Figure \ref{fig:fe_ext_vec}.

\begin{figure}[htb]
\centering
\begin{tikzcd}
H^1(M) \arrow[r, "\grad"] \arrow[d, "\Pi_{s, h}^{-, 0}"]
& H^{\curl}(M) \arrow[r, "\curl"] \arrow[d, "\Pi_{s, h}^{-, 1}"] & H^{\div}(M) \arrow[r, "\div"] \arrow[d, "\Pi_{s, h}^{-, 2}"]
& L^2(M) \arrow[d, "\Pi_{s, h}^{-, 3}"]  \\
\mathrm{CG}_s(\mathcal{T}_h) \arrow[r, "\grad"] 
& \mathrm{NED}_s^1(\mathcal{T}_h) \arrow[r, "\curl"] & \mathrm{RT}_s(\mathcal{T}_h) \arrow[r, "\div"]
& \mathrm{DG}_{s-1}(\mathcal{T}_h) 
\end{tikzcd}  
\caption{Equivalence between finite element differential forms and classical elements.}
\label{fig:fe_ext_vec}
\end{figure}

The finite element library \firedrake \cite{rathgeber2017firedrake} is used for the numerical investigation. The firedrake version to reproduce the experiment is archived at  \cite{zenodo/Firedrake-20231027.1}. The \firedrake component Slate \cite{gibson2020slate} is used to implement the local solvers and static condensation.\\

For what concerns the computation of the error, care has to be taken when considering normal trace variable $\omega^{k, \bm{n}}_h \in W_h^{k, \bm{n}} \equiv W_h^{k, \bm{t}}$. Identification of $\omega^{k, \bm{n}}_h$  with an element of $L^2\Omega^k(\partial \mathcal{T}_h)$ is only unique up to the annihilator $(W_h^{k, \bm{t}})^\perp$ \cite{stern2023hyb}. Therefore, the $L^2$ error has to be computed after the annihilator has been removed, which is
equivalent to taking the following projection $P_h \omega^{k, \bm{n}}_{\mathrm{ex}}$ of the exact solution $\omega^{k+1}_{\mathrm{ex}}$
\begin{equation}
    \inprBd[\partial \mathcal{T}_h]{v_h^{k, \bm{t}}}{P_h \omega^{k, \bm{n}}_{\mathrm{ex}}} = \inprBd[\partial \mathcal{T}_h]{v_h^{k, \bm{t}}}{\omega^{k, \bm{n}}_{\mathrm{ex}}}, \qquad \forall v_h^{k, \bm{t}} \in W_h^{k, \bm{t}}.
\end{equation}
The scaled $L^2$ norm over a cell boundary is given by $|||\cdot|||_{\partial K}:= h_T ||\cdot ||_{\partial K}$ where $h_T$ denoted the diameter of the cell $T \in \mathcal{T}_h$. For the overall mesh, we use the notation $|||\cdot|||_{\partial \mathcal{T}_h} = \sum_{T \in \mathcal{T}_h}|||\cdot|||_{\partial K}$. The error for $\omega^{k, \bm{n}}_h$ will then measure as
\begin{equation*}
    \text{Error } \omega^{k, \bm{n}}_h = |||\omega^{k, \bm{n}}_h - P_h \omega^{k, \bm{n}}_{\mathrm{ex}}|||_{\partial \mathcal{T}_h}.
\end{equation*}
The norm $|||\cdot|||_{\partial \mathcal{T}_h}$ is also used to measure convergence for the tangential trace $\omega_h^{k, \bm{t}}$.

\subsection{The wave equation in 3D}

The acoustic wave equation corresponds to the case $p=3$ and $q=1$. The energy variables are the pressure top-form $\dual{p}^3:=\dual{\alpha}^3$ and the stress one-form $\sigma^1:=\beta^{1}$. The Hamiltonian is given by
\begin{equation}
    H(\dual{p}^3, u^1) = \frac{1}{2} \int_M \dual{p}^3 \wedge {\star \dual{p}^3} + \sigma^1 \wedge \star \sigma^1,
\end{equation}
with its variational derivatives given by
\begin{equation}
p^0:=\delta_{\dual{p}^3} H = \star \dual{p}^3, \qquad  \dual{\sigma}^2:= \delta_{\sigma^1} H = \star \sigma^1,  
\end{equation}
leading to the pH system
\begin{equation}\label{eq:wave_eq}
\begin{bmatrix}
    c^{-2} & 0 \\
    0 & 0 \\
    \end{bmatrix}
    \begin{pmatrix}
    \partial_t \dual{p}^3 \\
    \partial_t \sigma^1
    \end{pmatrix} =
    -
    \begin{bmatrix}
    0 & \d \\
    \d & 0 \\
    \end{bmatrix}
    \begin{pmatrix}
     p^0 \\
     \dual{\sigma}^2
    \end{pmatrix}, \qquad 
    \begin{aligned}
    \tr p^0 \vert_{\Gamma_1} &= u^{0}_1, \\
    -\tr \dual{\sigma}^2 \vert_{\Gamma_2} &= \dual{u}^2_2.
 \end{aligned}
\end{equation}
where $c$ is the speed of propagation. The employment of the dual field discretization leads to the resolution of two systems:
\begin{itemize}
    \item the primal system \eqref{eq:discrete_primalPHhybrid} of outer oriented variables $\dual{p}^3_h, \; \dual{\sigma}^2_h, \; \dual{p}_h^{2, \bm{n}}, \; \dual{\sigma}_h^{2, \bm{t}}$;
    \item the dual system \eqref{eq:discrete_dualPHhybrid} of inner oriented variables $p^0_h, \; \sigma^1_h, \; \sigma_h^{0, \bm{n}}, \; p^{0, \bm{t}}_h$.
\end{itemize}

\subsubsection{Convergence test}
in this section we assess the rate of convergence for the primal and the dual formulation. For this test, the domain is the unit cube $$M = \{ (x,y,z) \in [0, 1]^3\}.$$
The boundary sub-partitions are selected to be
\begin{equation*}
    \Gamma_1 = \{(x,y,z) \vert \; x=0 \cup y=0 \cup z=0\}, \qquad \Gamma_2 = \{(x, y, z) \vert \; x=1 \cup y=1 \cup z=1 \}.
\end{equation*}
A structured tetrahedral mesh $\mathcal{T}_h$ is formed by partitioning $M$ into $N_{\text{el}} \times N_{\text{el}} \times N_{\text{el}}$ cubes, each of which is divided into six tetrahedra. The total simulation time is 1 and the time step is taken to be $\Delta t = T_{\mathrm{end}}/500$. The speed propagation is taken to be $c=1$. Introducing the functions
\begin{equation}
    g(x, y, z) = \sin(x) \sin(y) \sin(z), \qquad f(t) = \frac{1}{2} t^2,
\end{equation}
the manufactured solution of \eqref{eq:wave_eq} is given by
\begin{equation}\label{eq:wave_exsol}
\begin{aligned}
\dual{p}^3_{\mathrm{ex}} &= \star g \odv{f}{t}, \\    
\sigma^1_{\mathrm{ex}} &= -\d{g} f, 
\end{aligned} \qquad 
\begin{aligned}
p^0_{\mathrm{ex}} &= g \odv{f}{t}, \\
\dual{\sigma}^2_{\mathrm{ex}} &= -\star \d{g} f,
\end{aligned}
\end{equation}
A quadratic polynomial in time is taken to ensure that the error is only due to the spatial integration \revone{(the implicit midpoint method is of order two and integrates quadratic polynomials exactly)}. For this to be a true solution a forcing has to be introduced in the pressure equation
\begin{equation}
    \xi_p = g \odv[order=2]{f}{t} - \div(\grad g)f.
\end{equation}
The exact solution provides the appropriate inputs to be fed into the system
\begin{equation}
    u^0_1 = \left. \tr p^0_{\mathrm{ex}} \right\vert_{\Gamma_1}, \qquad u^{2}_2 =  -\tr \dual{\sigma}^2_{\mathrm{ex}} \vert_{\Gamma_2}.
\end{equation}

\begin{figure}[tbhp]%
\centering
\subfloat[][$L^2$ error for $\dual{p}^3_h$]{%
	\label{fig:err_p3}%
\includegraphics[width=0.48\columnwidth]{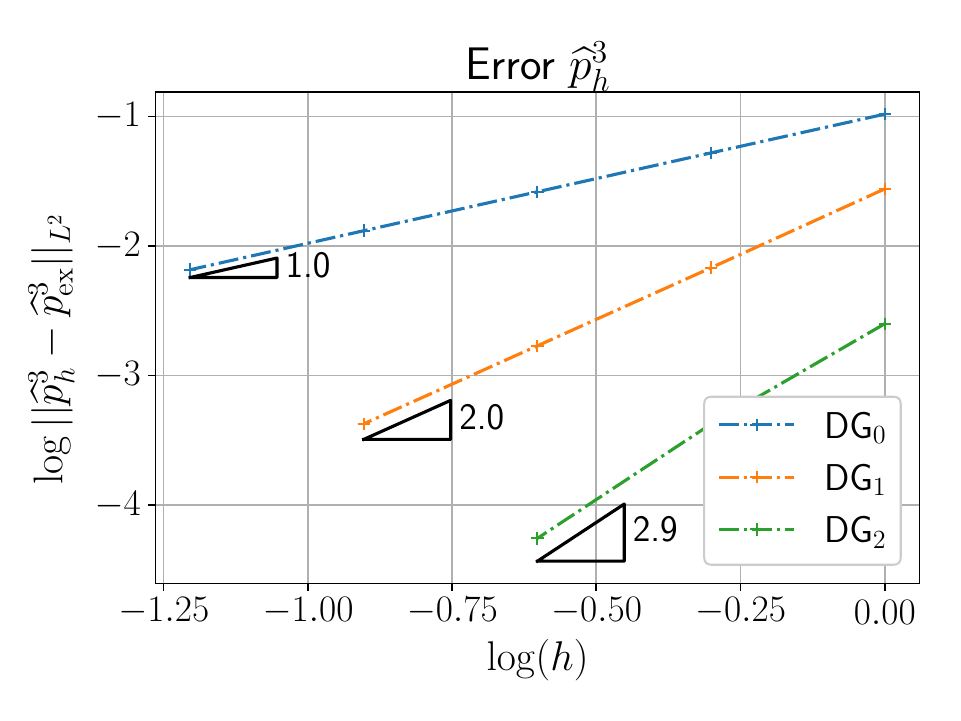}}%
 \hspace{8pt}%
\subfloat[][$H^{\div}$ error for $\dual{\sigma}^2_h$]{%
	\label{fig:err_u2}%
\includegraphics[width=0.48\columnwidth]{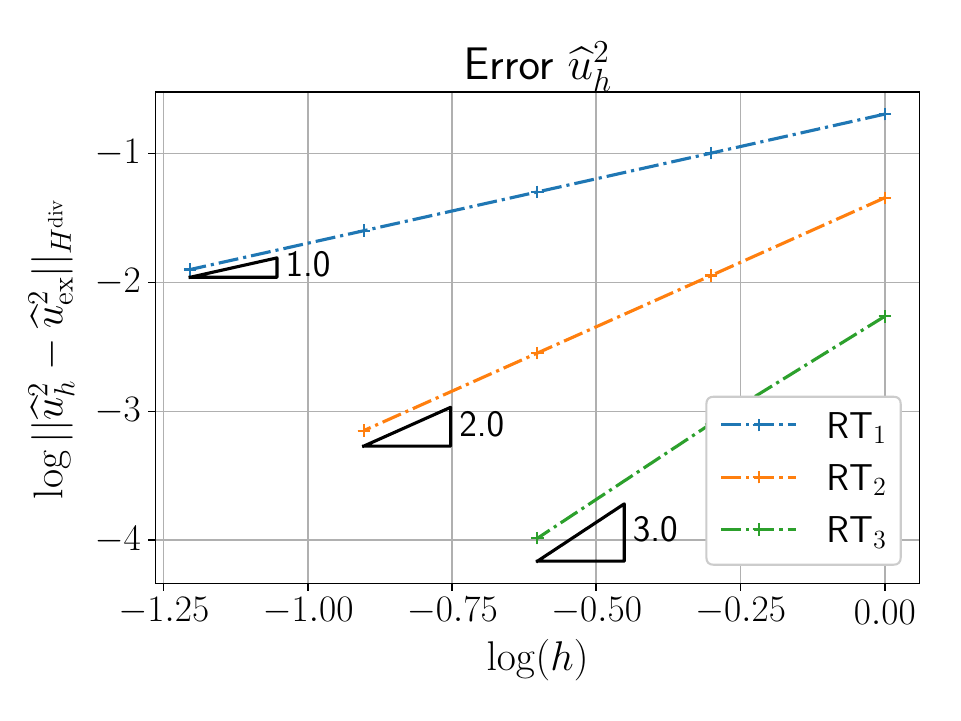}}
 \hspace{8pt}%
 \subfloat[][$L^2$ error for $\dual{p}^{2, \bm{n}}_h$]{%
	\label{fig:err_p2nor}%
\includegraphics[width=0.48\columnwidth]{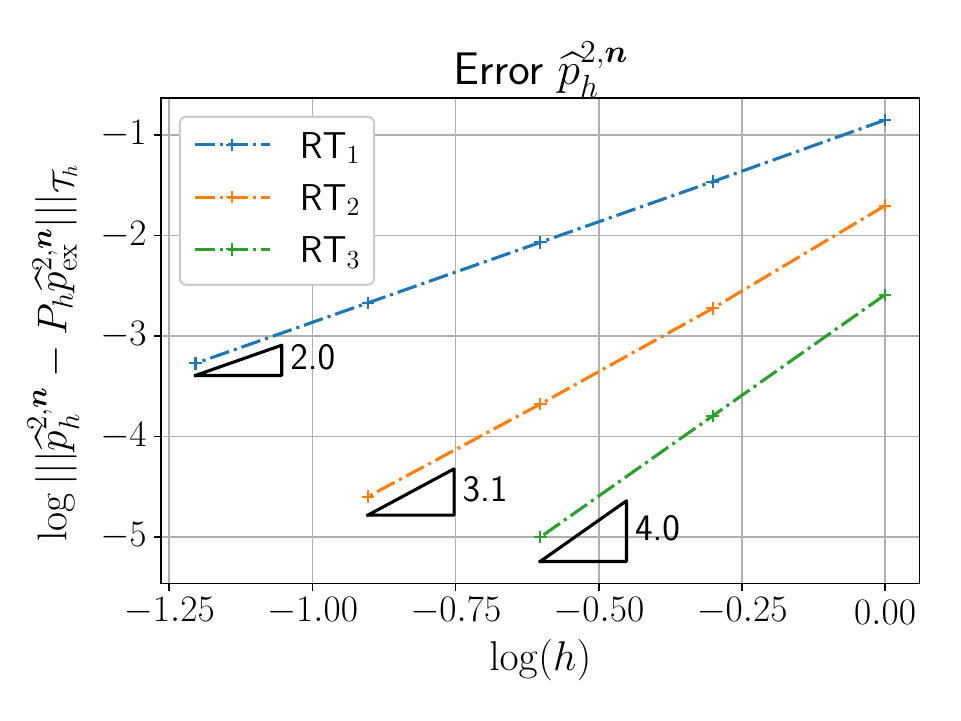}}%
\hspace{8pt}%
\subfloat[][$L^2$ error for $\dual{\sigma}^{2, \bm{t}}_h$]{%
	\label{fig:err_u2tan}%
\includegraphics[width=0.48\columnwidth]{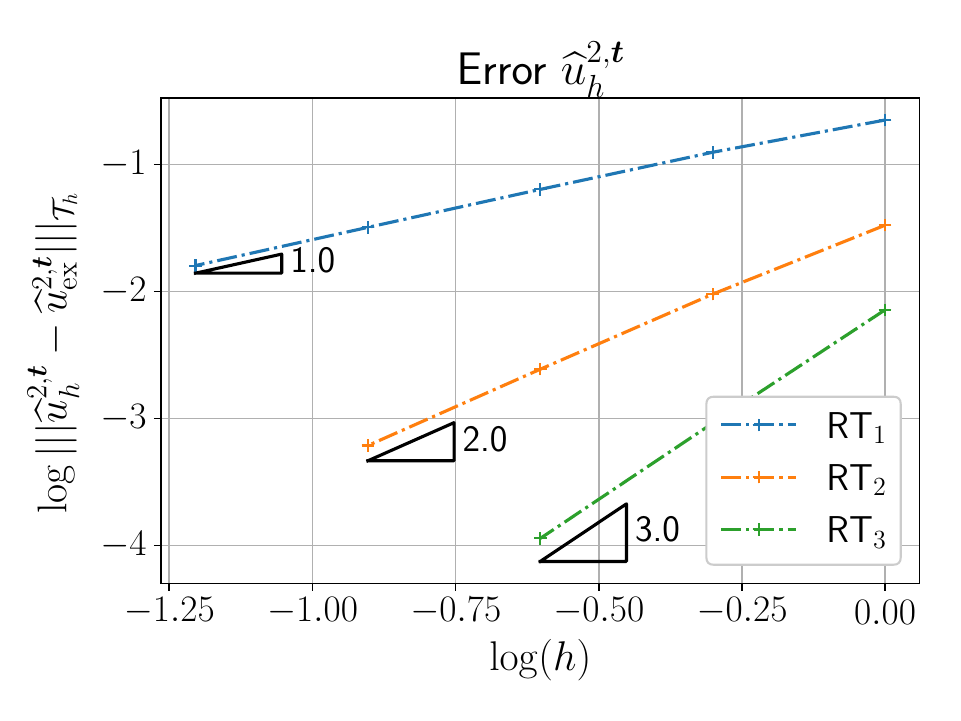}}
\caption{Convergence rate for the different variables in the primal formulation of the wave equation, measured at $T_{\text{end}}=1$ for $\Delta t = \frac{1}{500}$. }%
\label{fig:conv_var_wave_primal}%
\end{figure}

\begin{figure}[tbhp]%
\centering
\subfloat[][$H^1$ error for $p^0_h$]{%
\label{fig:err_p0}%
\includegraphics[width=0.48\columnwidth]{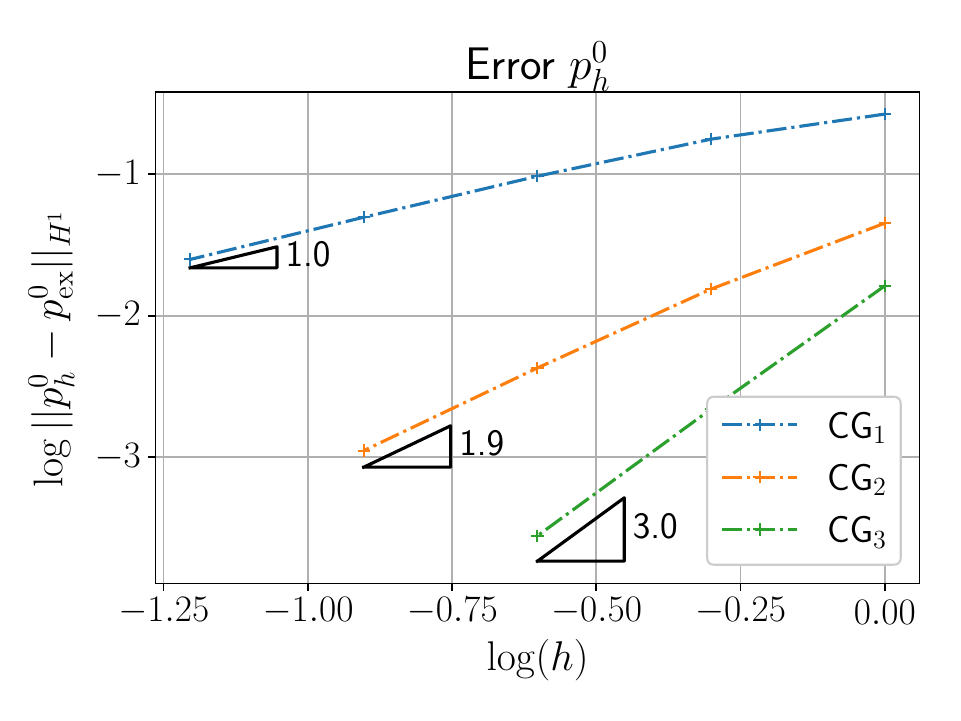}}%
\hspace{8pt}%
\subfloat[][$H^{\curl}$ error for $\sigma^1_h$]{%
	\label{fig:err_u1}%
\includegraphics[width=0.48\columnwidth]{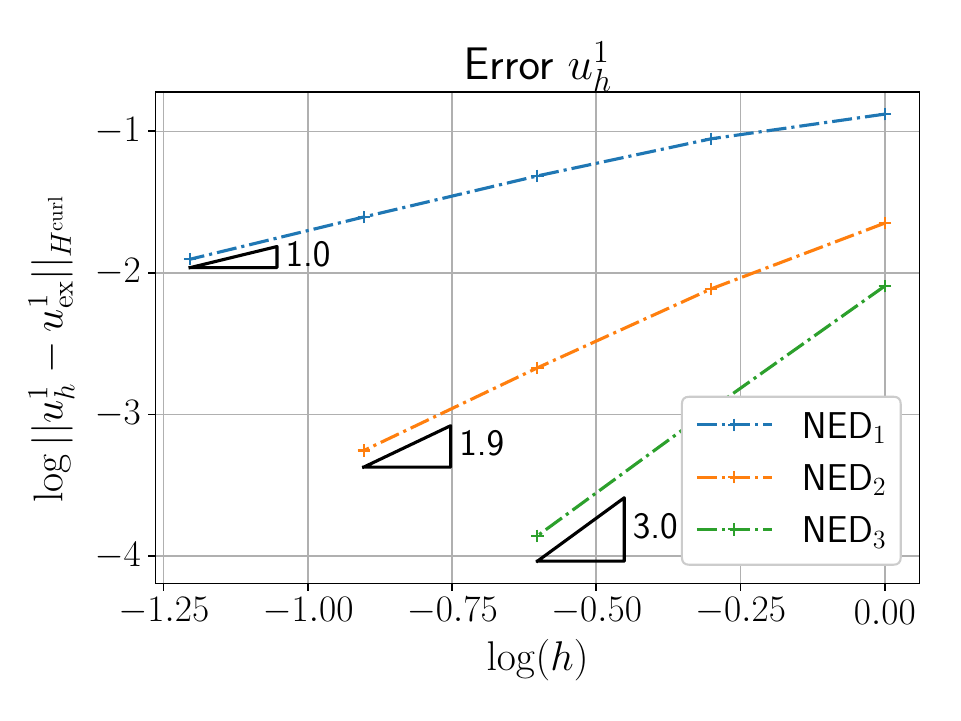}}
\hspace{8pt}%
 \subfloat[][$L^2$ error for ${\sigma}^{0, \bm{n}}_h$]{%
	\label{fig:err_u0nor}%
\includegraphics[width=0.48\columnwidth]{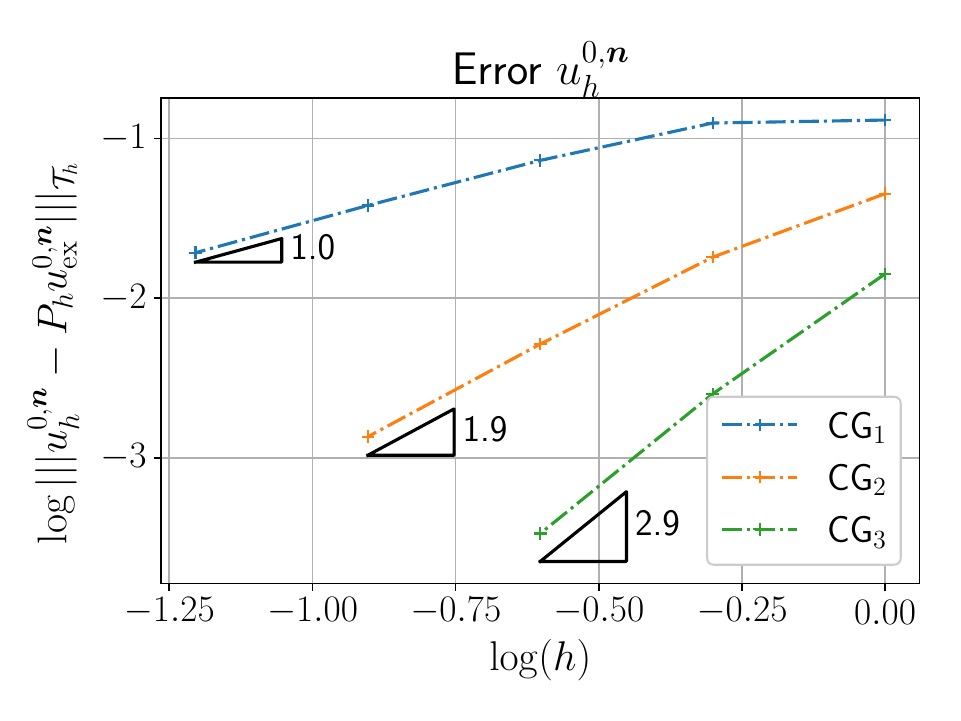}}%
\hspace{8pt}%
\subfloat[][$L^2$ error for ${p}^{0, \bm{t}}_h$]{%
	\label{fig:err_p0tan}%
\includegraphics[width=0.48\columnwidth]{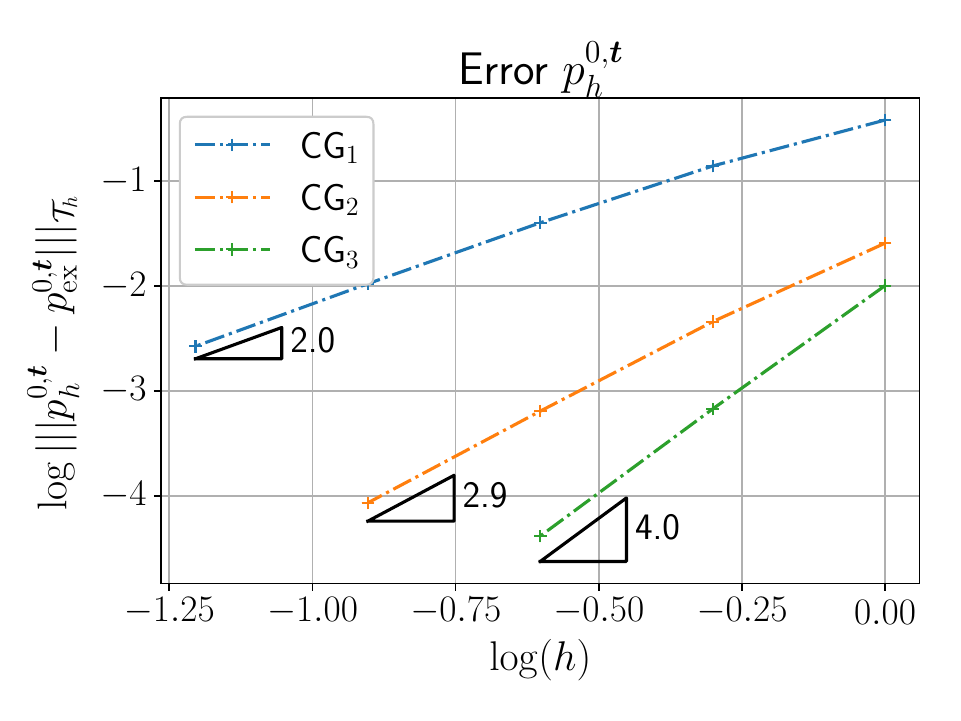}}
\caption{Convergence rate for the different variables in the dual formulation of the wave equation, measure at at $T_{\text{end}}=1$ for $\Delta t = \frac{1}{500}$.}%
\label{fig:conv_var_wave_dual}%
\end{figure}

\begin{figure}[p]%
\centering
\subfloat[][$L^2$ norm of the difference $\dual{p}^3_h - p^0_h$]{%
	\label{fig:diff_p30}%
\includegraphics[width=0.48\columnwidth]{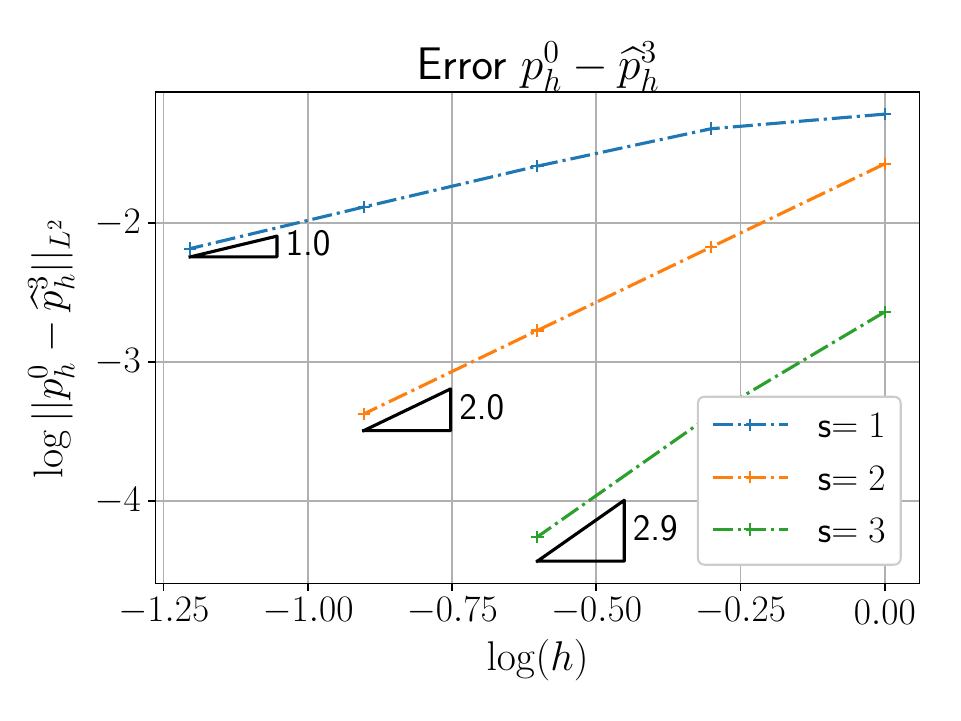}}%
\hspace{8pt}%
\subfloat[][$L^2$ norm of the difference $u^1_h - \dual{u}^2_h$]{%
	\label{fig:diff_u12}%
\includegraphics[width=0.48\columnwidth]{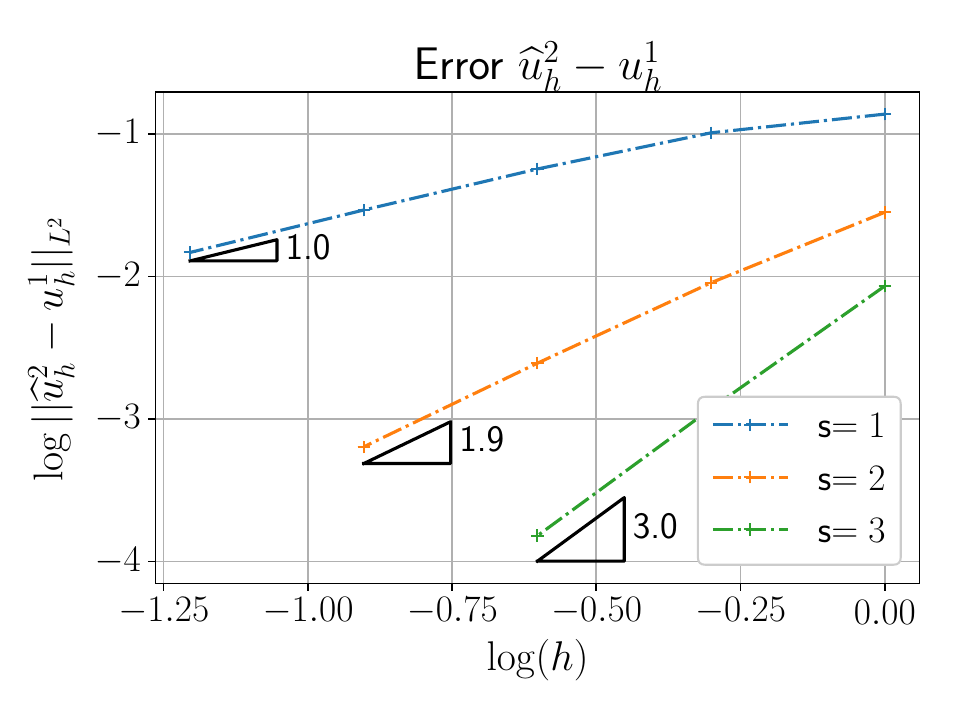}}%
\caption{$L^2$ difference of the dual representation of the solution for the wave equation at $T_{\text{end}}=1$ for $\Delta t = \frac{1}{500}$.}%
\label{fig:diff_dual_wave}%
\end{figure}

The error is measured in the Sobolev norm or the tangential norm for the facet variables at the final time $T_{\mathrm{end}}$. In Fig. \ref{fig:conv_var_wave_primal} the convergence rate of the variables is plotted against the mesh size. It can be immediately noticed that variable $\dual{p}^{2, \bm{n}}_h$ superconverges (cf. Fig. \ref{fig:err_p2nor}). This behaviour is well known for the RT and BDM elements in the static case \cite{douglas1985}. The other variables converge with the optimal order. For what concerns the dual system, the results are shown in Fig. \ref{fig:conv_var_wave_dual}. The $L^2$ tangential trace of the pressure superconverges (cf. Fig. \ref{fig:err_p0tan}) since the pressure converges with order $h^s$ in the $H^1$. All other variables converge with optimal order. In Fig. \ref{fig:diff_dual_wave}, the $L^2$ norm of the difference between the dual representation of variables. As in the dual field continuous Galerkin formulation, the dual representation of the variables converges under h-p refinement with order $h^s$. \\

The size reduction between the continuous and hybrid formulation is reported in Tables \ref{tab:dofs_cont_hyb_waveprimal},~\ref{tab:dofs_cont_hyb_wavedual}. For the primal formulation one only solves for $40\%$ of the degrees of freedom when third order polynomials are used. For the dual the size reduction is way more impressive (the hybrid formulation dimension is $20\%$ of the continuous formulation  in the worst case) as the broken Nédélec space is completely discarded when hybridization is used.

\begin{table}[tbhp]
    \centering
    \begin{tabular}{c|c|c|c|c}
        Pol. Degree $s$ & $N_{\text{elem}}$ & $N^\circ$ dofs. continuous & $N^\circ$ dofs. hybrid & $N_{\text{hyb}}/N_{\text{cont}}$ \\ \hline
        \multirow{5}{*}{1} & 1 & 24 & 18 & 75 \% \\
         & 2 & 168 & 120 & 71\% \\
         & 4 & 1248 & 864 & 69\% \\
         & 8 & 9600 & 6528 & 68\% \\
         & 16 & 75264 & 50688 & 67\% \\ \hline
         \multirow{4}{*}{2} & 1 & 96 & 54 & 56\% \\
          & 2 & 696 & 360 & 52\% \\
          & 4 & 5280 & 2592 & 49\% \\
          & 8 & 41088 & 19584 & 47\% \\ \hline
          \multirow{3}{*}{3} & 1 & 240 & 108 & 45\% \\    
           & 2 & 1776 & 720 & 41\% \\ 
           & 4 & 13632 & 5184 & 38\% \\ \hline
    \end{tabular}
    \caption{Size of the primal system $\dual{p}^3, \dual{u}^2$ for the wave equation: continuous and hybrid formulation.}
    \label{tab:dofs_cont_hyb_waveprimal}
\end{table}

\begin{table}[tbhp]
    \centering
    \begin{tabular}{c|c|c|c|c}
        Pol. Degree $s$ & $N_{\text{elem}}$ & $N^\circ$ dofs. continuous & $N^\circ$ dofs. hybrid & $N_{\text{hyb}}/N_{\text{cont}}$ \\ \hline
        \multirow{5}{*}{1} & 1 & 44 & 8 & 18\% \\
         & 2 & 315 & 27 & 9\% \\
         & 4 & 2429 & 125 & 5\% \\
         & 8 & 19161 & 729 & 4\% \\
         & 16 & 152369 & 4913 & 3\% \\ \hline
        \multirow{4}{*}{2} & 1 & 147 & 27 & 18\% \\
         & 2 & 1085 & 125 & 12\% \\
         & 4 & 8409 & 729 & 9\% \\
         & 8 & 66353 & 4913 & 7\% \\ \hline
         \multirow{3}{*}{3} & 1 & 334 & 64 & 19\% \\
         & 2 & 2503 & 343 & 13\% \\
         & 4 & 19477 & 2197 & 11\%
    \end{tabular}
    \caption{Size of the dual system ${p}^0, {u}^1$ for the wave equation: continuous and hybrid formulation.}
    \label{tab:dofs_cont_hyb_wavedual}
\end{table}

\begin{figure}[tbhp]%
\centering
\subfloat[][{$\dual{p}^3$ at $t=1 \; [\mathrm{s}]$}]{%
\label{fig:p3_1_discontinuous}%
\includegraphics[width=0.4\columnwidth]{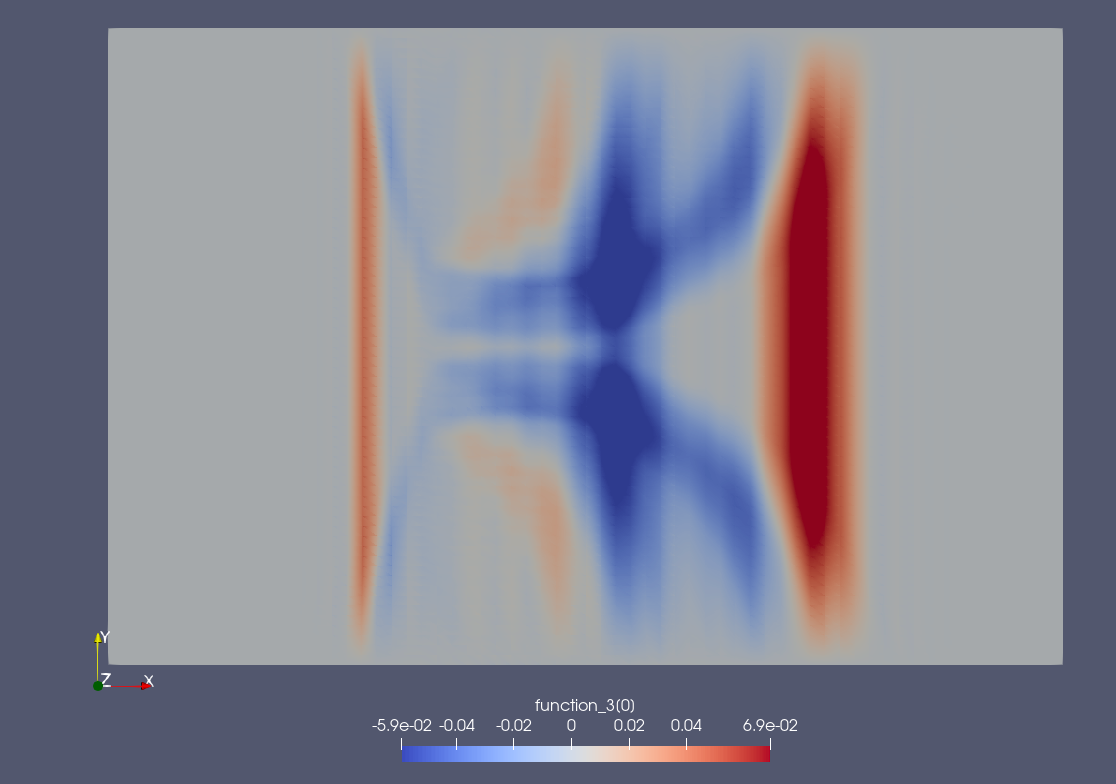}}%
\hspace{8pt}%
\subfloat[][{$p^0$ at $t=1 \; [\mathrm{s}]$}]{%
	\label{fig:p0_1_discontinuous}%
\includegraphics[width=0.4\columnwidth]{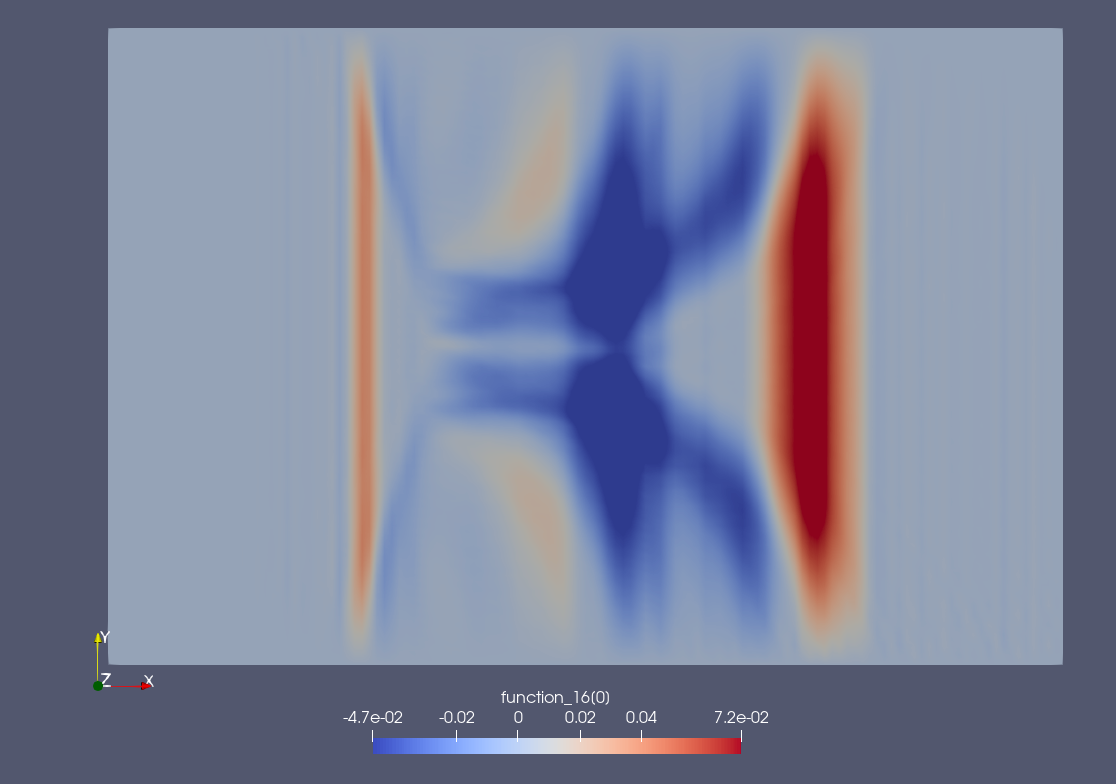}}\\
\subfloat[][{$\dual{p}^3$ at $t=2 \; [\mathrm{s}]$}]{%
	\label{fig:p3_2_discontinuous}%
\includegraphics[width=0.4\columnwidth]{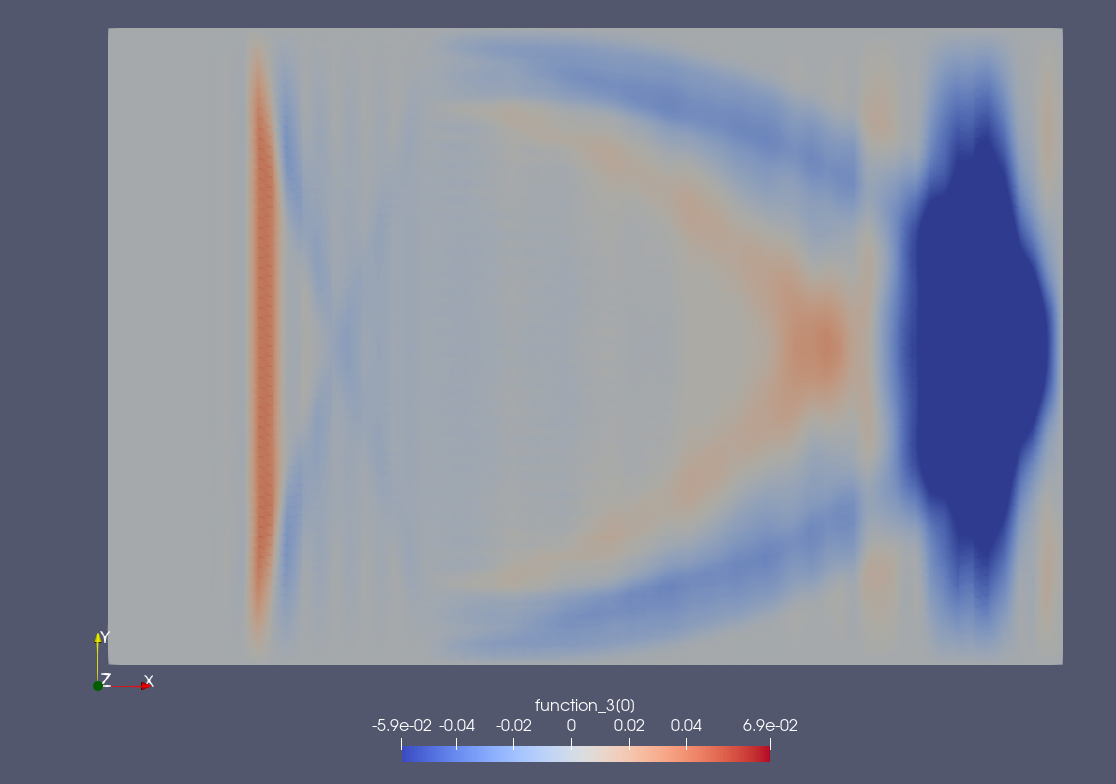}}%
\hspace{8pt}%
\subfloat[][{$p^0$ at $t=2 \; [\mathrm{s}]$}]{%
	\label{fig:p0_2_discontinuous}%
\includegraphics[width=0.4\columnwidth]{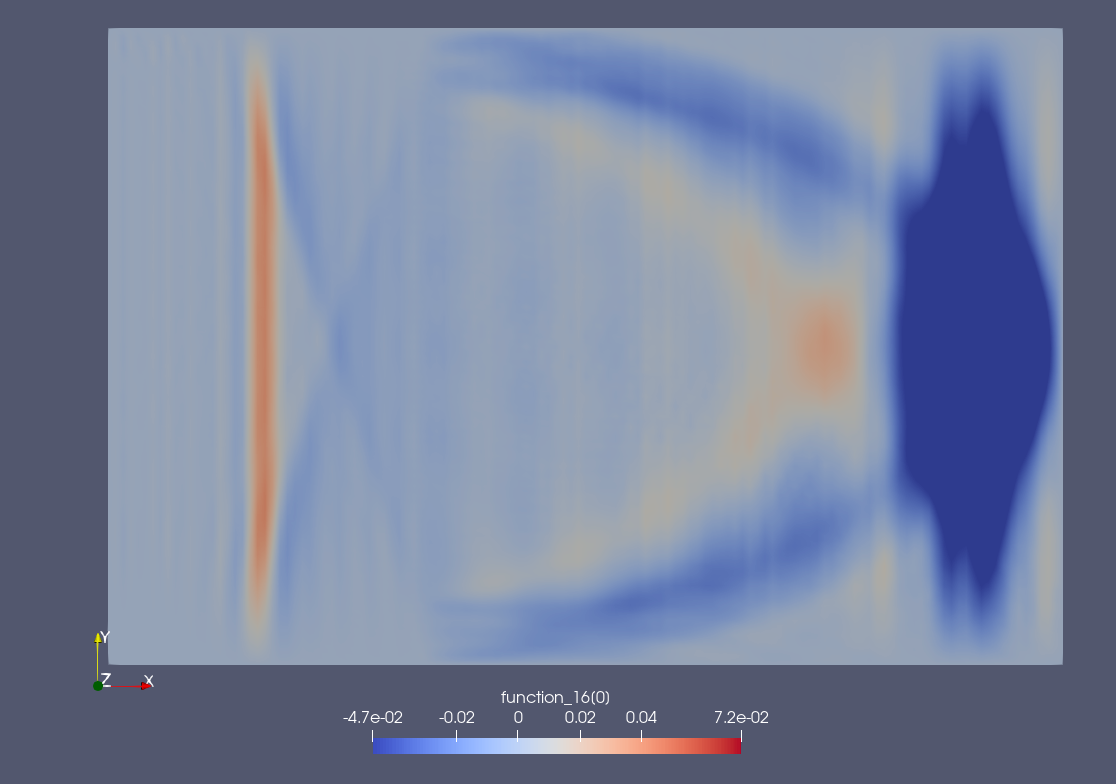}}\\
\subfloat[][{$\dual{p}^3$ at $t=3 \; [\mathrm{s}]$}]{%
	\label{fig:p3_3_discontinuous}%
\includegraphics[width=0.4\columnwidth]{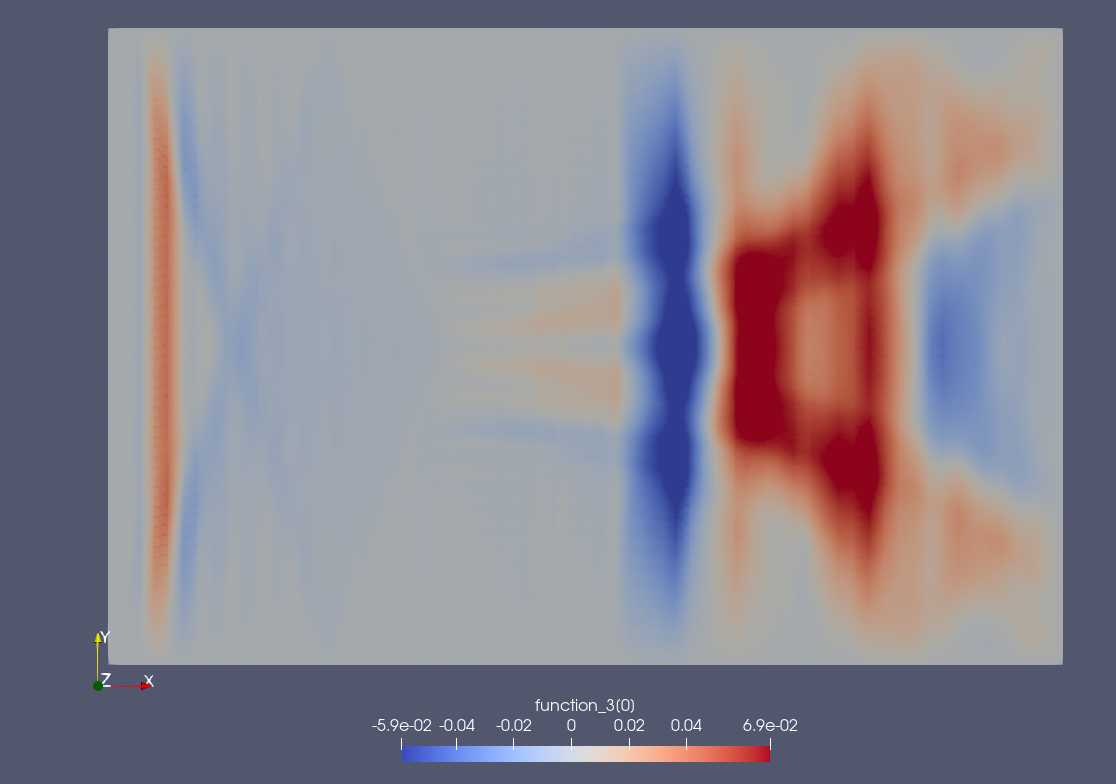}}%
\hspace{8pt}%
\subfloat[][{$p^0$ at $t=3 \; [\mathrm{s}]$}]{%
	\label{fig:p0_3_discontinuous}%
\includegraphics[width=0.4\columnwidth]{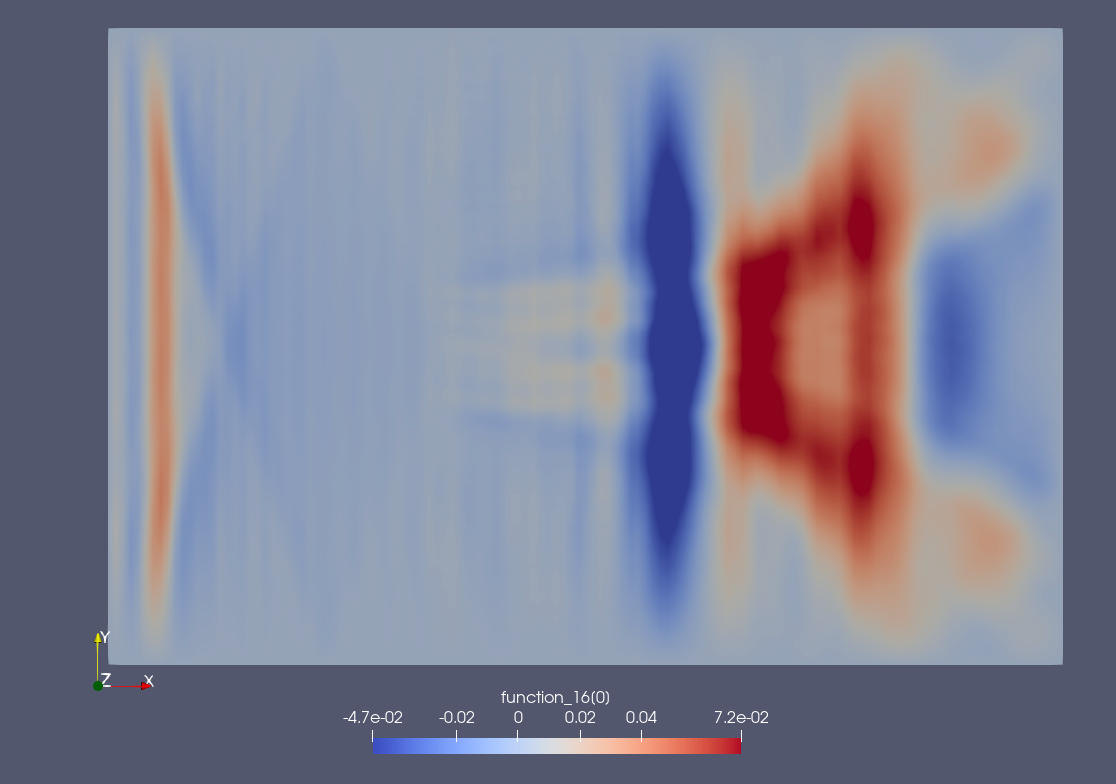}}\\
\subfloat[][{$\dual{p}^3$ at $t=4 \; [\mathrm{s}]$}]{%
	\label{fig:p3_4_discontinuous}%
\includegraphics[width=0.4\columnwidth]{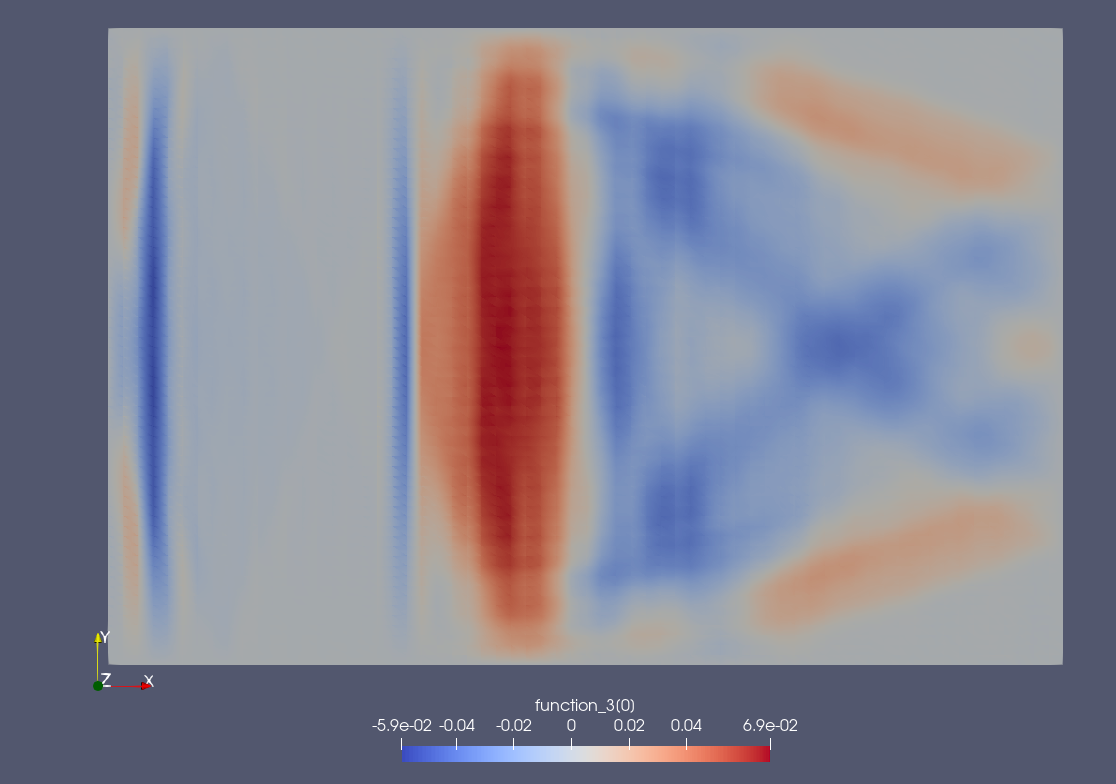}}%
\hspace{8pt}%
\subfloat[][{$p^0$ at $t=4 \; [\mathrm{s}]$}]{%
	\label{fig:p0_4_discontinuous}%
\includegraphics[width=0.4\columnwidth]{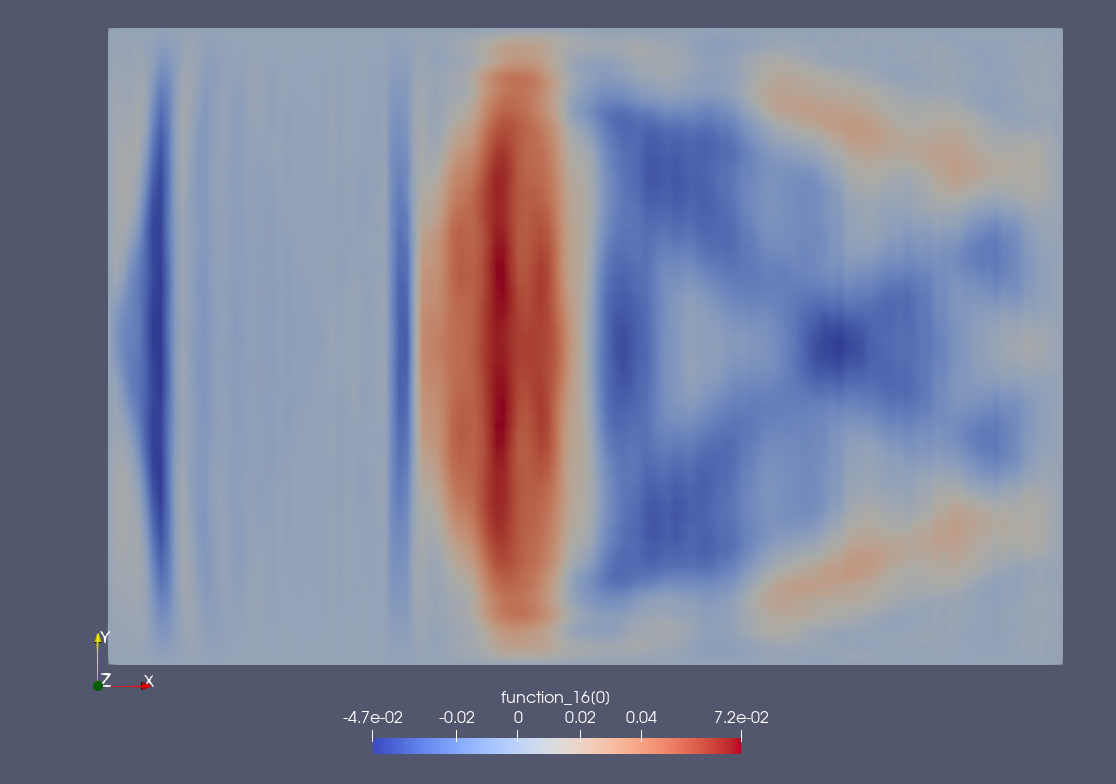}}%
\caption{Snapshots of the pressure for a domain with discontinuous coefficients calculated using the primal and dual system at different time instants.}%
\label{fig:discontinuous_coeff_wave}%
\end{figure}

{
\subsubsection{An example with discontinuous coefficients}
\revone{We consider the example  \cite{grote2006dg_wave} on the domain $M = [0, 3] \times [0, 1]$, where an heterogeneous material is considered}
\begin{equation*}
    c^{2} = \begin{cases}
        0.1 \\
        1 
    \end{cases} \qquad
    \begin{aligned}
    x \le 1, \\
    x > 1.
    \end{aligned}
\end{equation*}
The system is excited by the following forcing 
\begin{equation*}
    \xi_p = \begin{cases}
    1, \\
    0, \\
    \end{cases}
    \qquad
    \begin{aligned}
    1.2 < x < 1.4 \quad \text{and} \quad t \le 0.2, \\
    \text{ else.}
    \end{aligned}
\end{equation*}
Homogeneous Dirichlet conditions are considered. The system is simulated using 8192 elements and second order polynomial. The total time simulation time is $T_{\mathrm{end}} = 4[\mathrm{s}]$ and the time step is $\Delta t = 0.002[\mathrm{s}]$. The snapshots, reported in Fig. \ref{fig:discontinuous_coeff_wave}, show the agreement between the primal and dual formulation and the different velocity speed of the front wave in the different parts of the domain.
}

\subsection{The Maxwell equations in $3D$}
The Maxwell equations corresponds to the case $p=2, \; q=2$. The energy variables correspond to the electric displacement two form $\dual{D}^2 = \dual{\alpha}^2$ and the magnetic field $B^2=\beta^2$. The Hamiltonian reads
\begin{equation}
    H(\dual{D}^2, B^2) = \frac{1}{2} \int_M \dual{D}^2 \wedge \star \dual{D}^2 + B^2 \wedge \star B^2,
\end{equation}
The variational derivative of the Hamiltonian are given by 
\begin{equation}
    E^1 := \delta_{\dual{D}^2} H = \star \dual{D}^2, \qquad \dual{H}^1 := \delta_{B^2} H = \star B^2.
\end{equation}
Variables $E^1, \; \dual{H}^1$ are the electric field and the magnetizing field respectively. Since the reduction of the constitutive equation is such to keep only the efforts variables and their duals, the following dynamical system is obtained.
\begin{equation}\label{eq:maxwell_eq}
\begin{bmatrix}
    \varepsilon & 0 \\
    0 & \mu \\
    \end{bmatrix}
    \begin{pmatrix}
    \partial_t \dual{E}^2\\
    \partial_t H^2
    \end{pmatrix} = 
    \begin{bmatrix}
    0 & \d^1 \\
    -\d^1 & 0 \\
    \end{bmatrix}
    \begin{pmatrix}
    {E}^1\\
    \dual{H}^1
    \end{pmatrix},
\end{equation}
where $\dual{E}^2 = \star E^1, \; H^2= \star \dual{H}^1$ and $\varepsilon, \; \mu$ are the electric permittivity and magnetic permeability. 
The employment of the dual field discretization leads to the resolution of two systems:
\begin{itemize}
    \item the primal system \eqref{eq:discrete_primalPHhybrid} of outer oriented variables $\dual{E}^2_h, \dual{H}^1_h, \dual{E}^{1, \bm{n}}_h, \dual{H}^{1, \bm{t}}_h$;
    \item the dual system \eqref{eq:discrete_dualPHhybrid} of inner oriented variables ${E}^1_h, H^2_h, {H}_h^{1, \bm{n}}, {E}_h^{1, \bm{t}}$.
\end{itemize}

\subsubsection{Convergence results}
The same mesh and boundary partitions considered in the corresponding section for the wave equation are here considered. The electric and magnetic permeability are taken to be $\varepsilon=1, \; \mu=1$.

Given the functions
\begin{equation}
    \bm{g}(x, y, z) = \begin{pmatrix}
    -\cos(x)\sin(y)\sin(z) \\
    0 \\
    \sin(x)\sin(y)\cos(z)
    \end{pmatrix}, \qquad f(t) = \frac{1}{2} t^2
\end{equation}
The manufactured solution is given by
\begin{equation}\label{eq:maxwell_exsol}
\begin{aligned}
\dual{E}^2_{\mathrm{ex}} &= \star \bm{g}^\flat \odv{f}{t}, \\    
H^2_{\mathrm{ex}} &= -\d{\bm{g}^\flat} f, 
\end{aligned} \qquad
\begin{aligned}
E^1_{\mathrm{ex}} &= \bm{g}^\flat \odv{f}{t}, \\    
\dual{H}^1_{\mathrm{ex}} &= -\star \d{\bm{g}^\flat} f.
\end{aligned}
\end{equation}
A current has to be introduced in the electric field to make this a true solution of the problem
\begin{equation*}
    j = \bm{g} \odv[order=2]{f}{t} + \curl\curl \bm{g} f .
\end{equation*}

The exact solution provides the appropriate inputs to be fed into the system
\begin{equation}
    u^1_1 = \left. \tr E^1_{\mathrm{ex}}  \right\vert_{\Gamma_1}, \qquad \dual{u}^1_2 = \tr \dual{H}^1_{\mathrm{ex}} \vert_{\Gamma_2}.
\end{equation}

\begin{figure}[tbhp]%
\centering
\subfloat[][$L^2$ error for $\dual{E}^2_h$]{%
	\label{fig:err_L2_E2}%
\includegraphics[width=0.48\columnwidth]{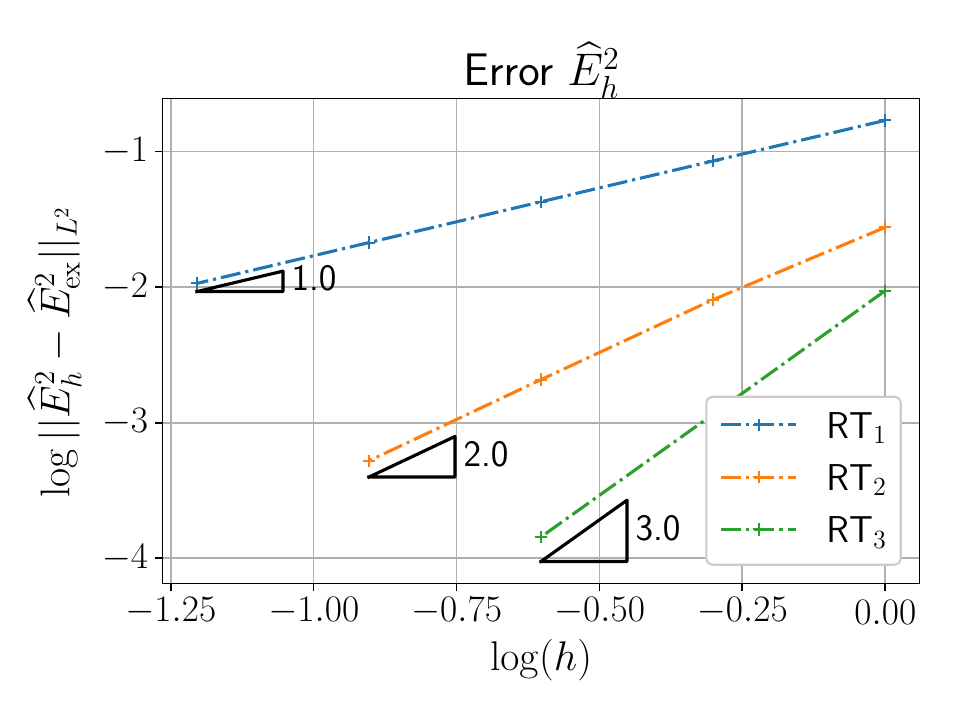}}%
 \hspace{8pt}%
\subfloat[][$H^{\div}$ error for $\dual{E}^2_h$]{%
	\label{fig:err_Hdiv_E2}%
\includegraphics[width=0.48\columnwidth]{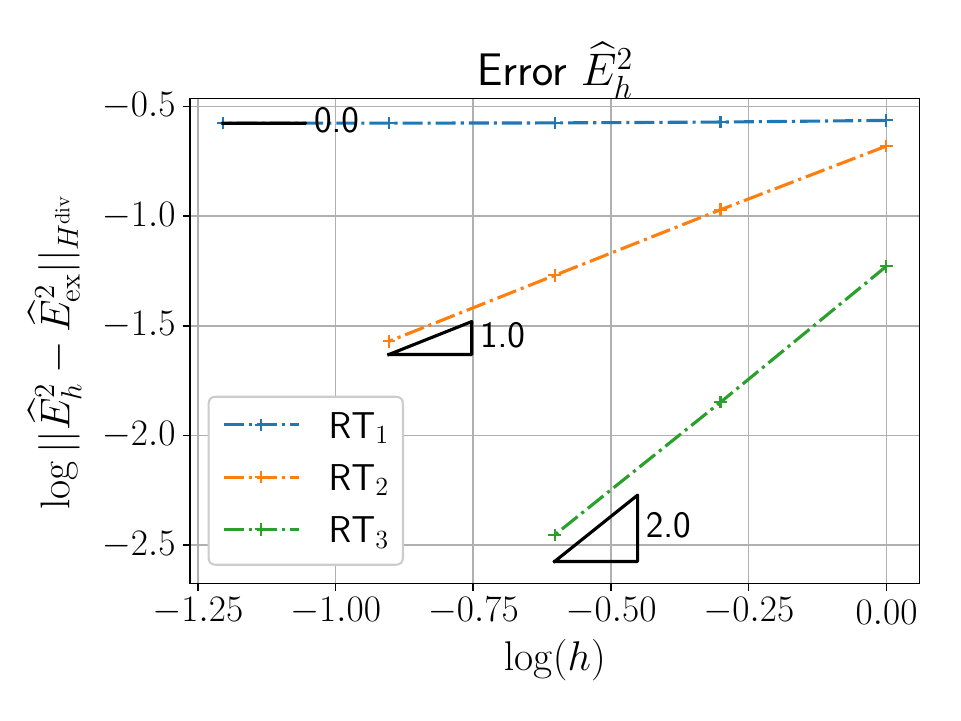}}%
 \hspace{8pt}%
\subfloat[][$H^{\curl}$ error for $\dual{H}^1_h$]{%
	\label{fig:err_H1}%
\includegraphics[width=0.48\columnwidth]{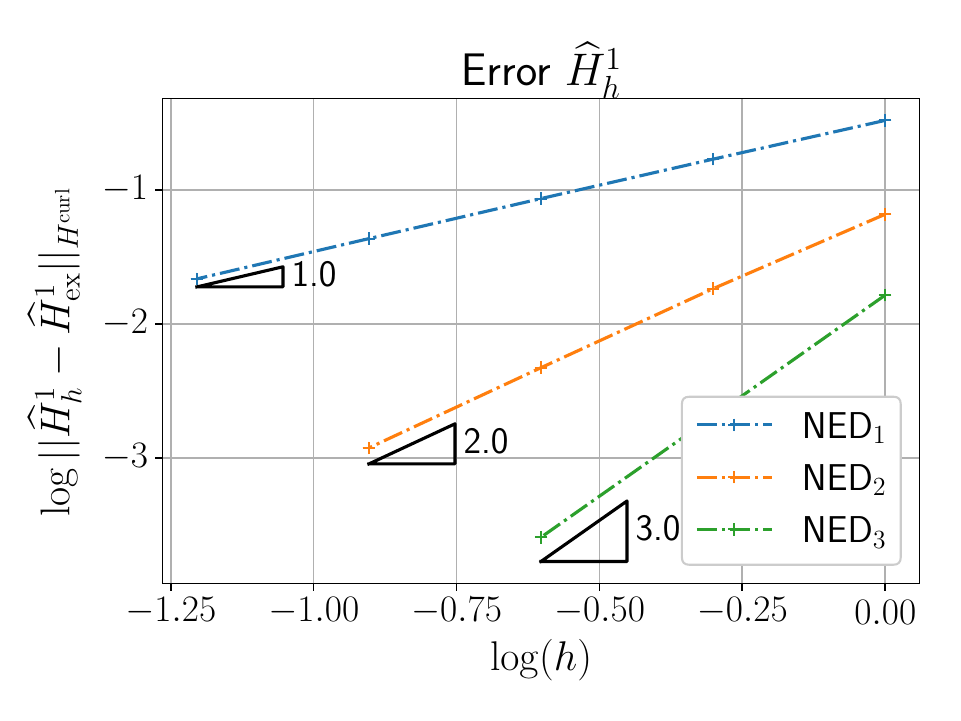}}
 \hspace{8pt}%
 \subfloat[][$L^2$ error for $\dual{E}^{1, \bm{n}}_h$]{%
	\label{fig:err_E1nor}%
	\includegraphics[width=0.48\columnwidth]{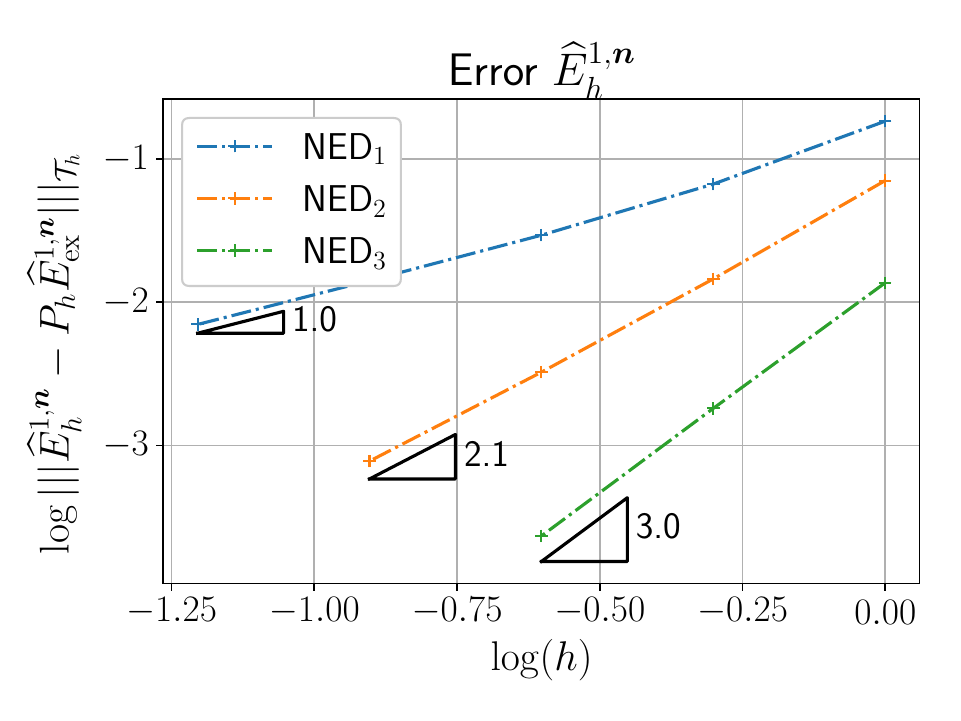}}%
\hspace{8pt}%
\subfloat[][$L^2$ error for $\dual{H}^{1, \bm{t}}_h$]{%
	\label{fig:err_H1tan}%
 \includegraphics[width=0.48\columnwidth]{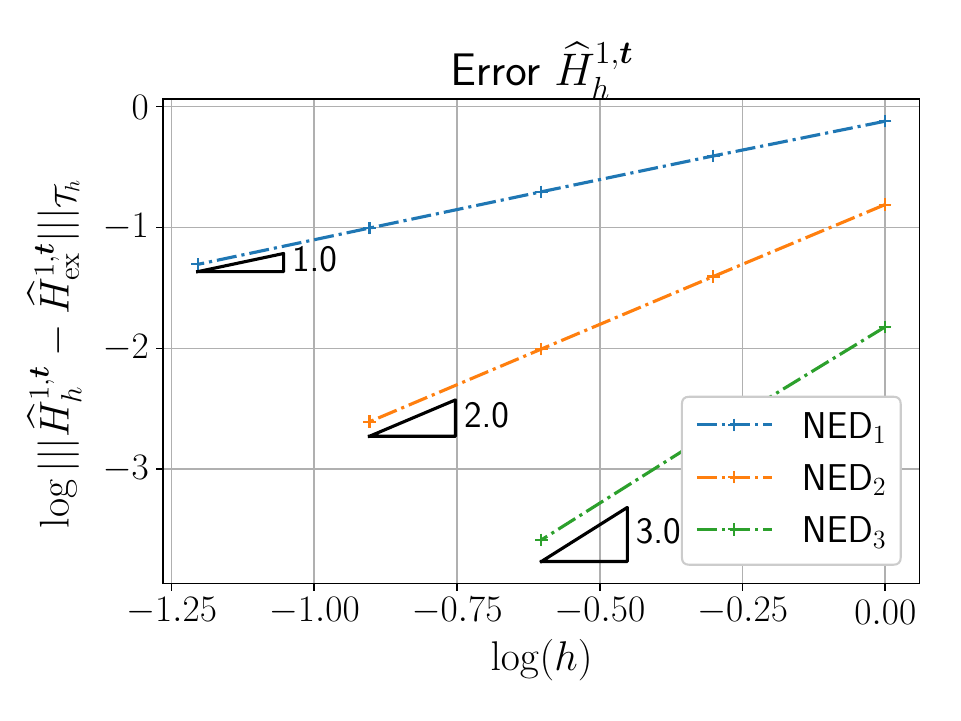}}
\caption{Convergence rate for the different variables in the primal formulation of the Maxwell equations, measured at $T_{\text{end}}=1$ for $\Delta t = \frac{1}{500}$. }%
\label{fig:conv_var_maxwell_primal}%
\end{figure}

\begin{figure}[tbhp]%
\centering
\subfloat[][$H^{\curl}$ error for $E^1_h$]{%
	\label{fig:err_E1}%
	\includegraphics[width=0.48\columnwidth]{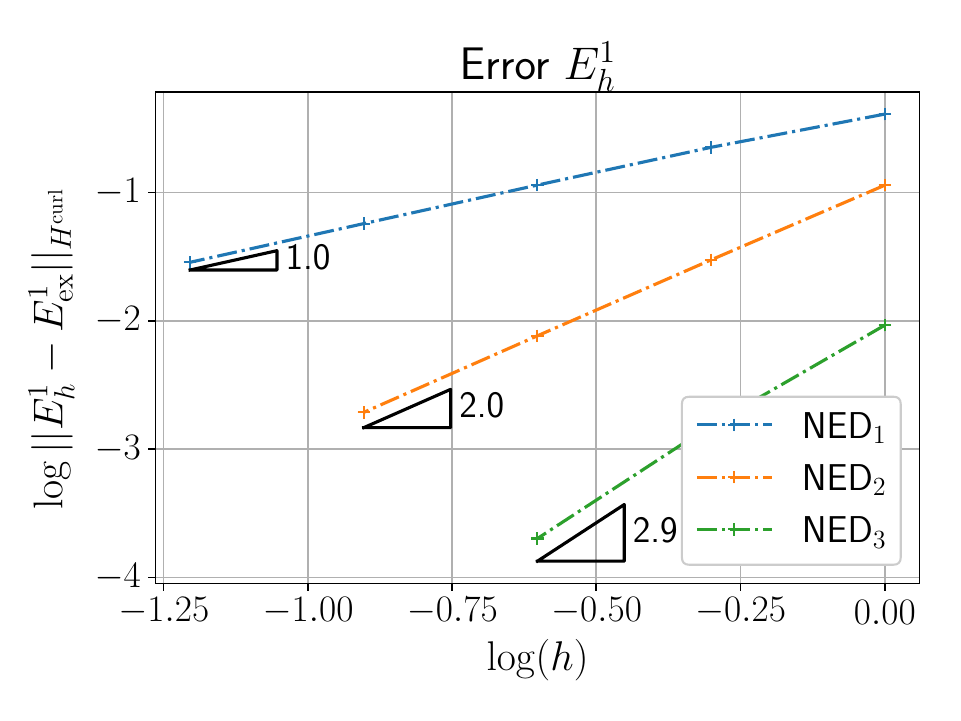}}%
\hspace{8pt}%
\subfloat[][$H^{\div}$ error for $H^2_h$]{%
	\label{fig:err_H2}%
\includegraphics[width=0.48\columnwidth]{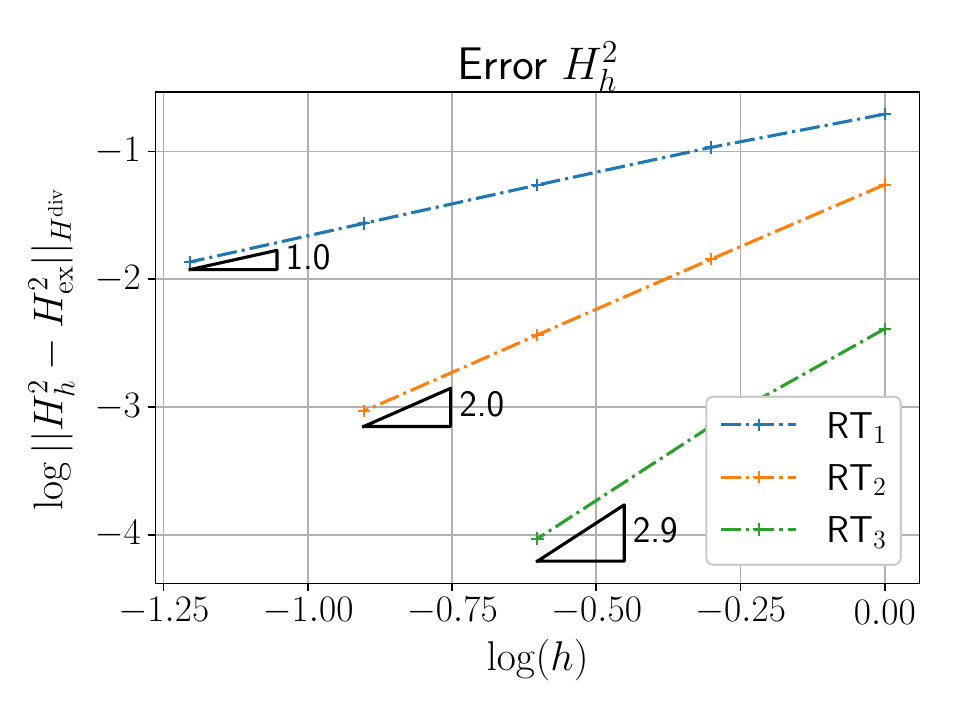}}
\hspace{8pt}%
 \subfloat[][$L^2$ error for ${H}^{1, \bm{n}}_h$]{%
	\label{fig:err_H1nor}%
\includegraphics[width=0.48\columnwidth]{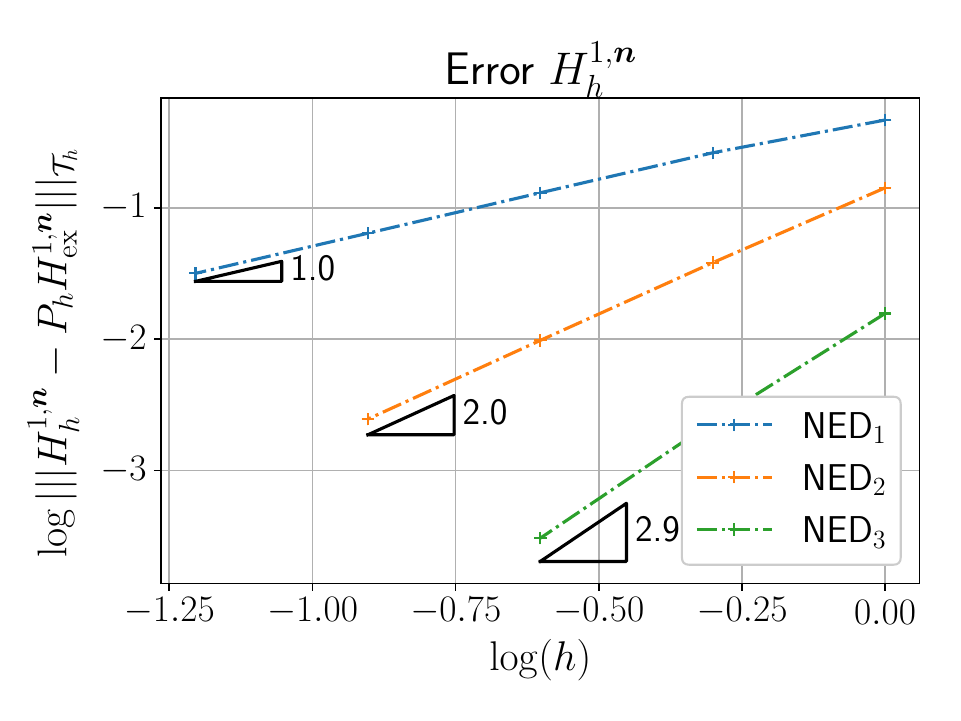}}%
\hspace{8pt}%
\subfloat[][$L^2$ error for ${E}^{1, \bm{t}}_h$]{%
	\label{fig:err_E1tan}%
\includegraphics[width=0.48\columnwidth]{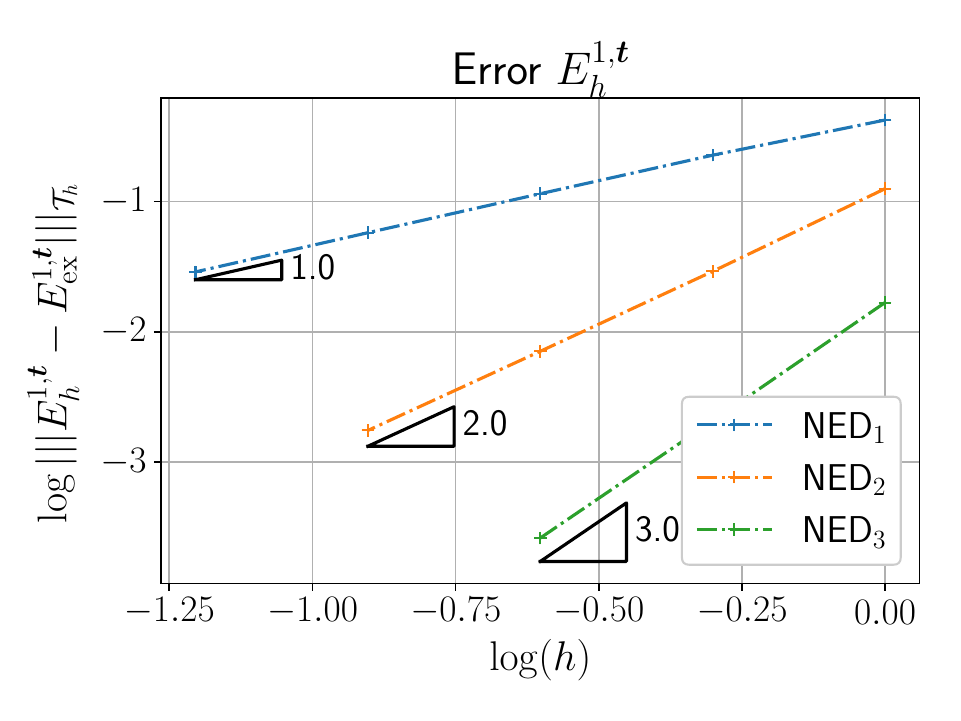}}
\caption{Convergence rate for the different variables in the dual formulation of the Maxwell equation, measure at at $T_{\text{end}}=1$ for $\Delta t = \frac{1}{100}$.}%
\label{fig:conv_var_maxwell_dual}%
\end{figure}

\begin{figure}[p]%
\centering
\subfloat[][$L^2$ norm of the difference $\dual{E}^2_h - E^1_h$]{%
	\label{fig:diff_E12}%
	\includegraphics[width=0.48\columnwidth]{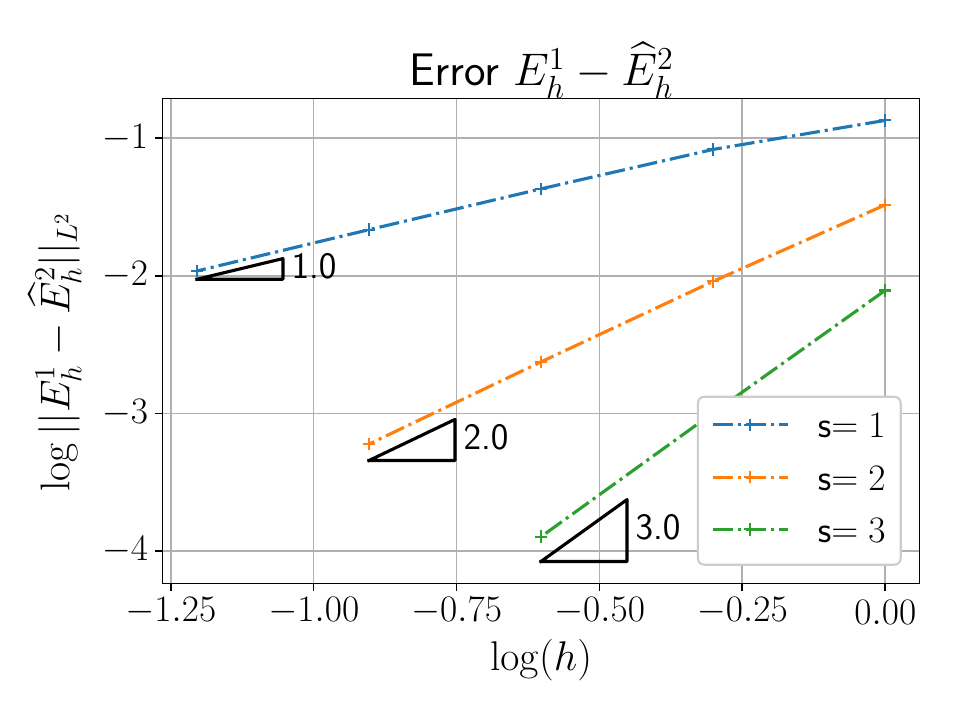}}%
\hspace{8pt}%
\subfloat[][$L^2$ norm of the difference $H^2_h - \dual{H}^1_h$]{%
	\label{fig:diff_H12}%
	\includegraphics[width=0.48\columnwidth]{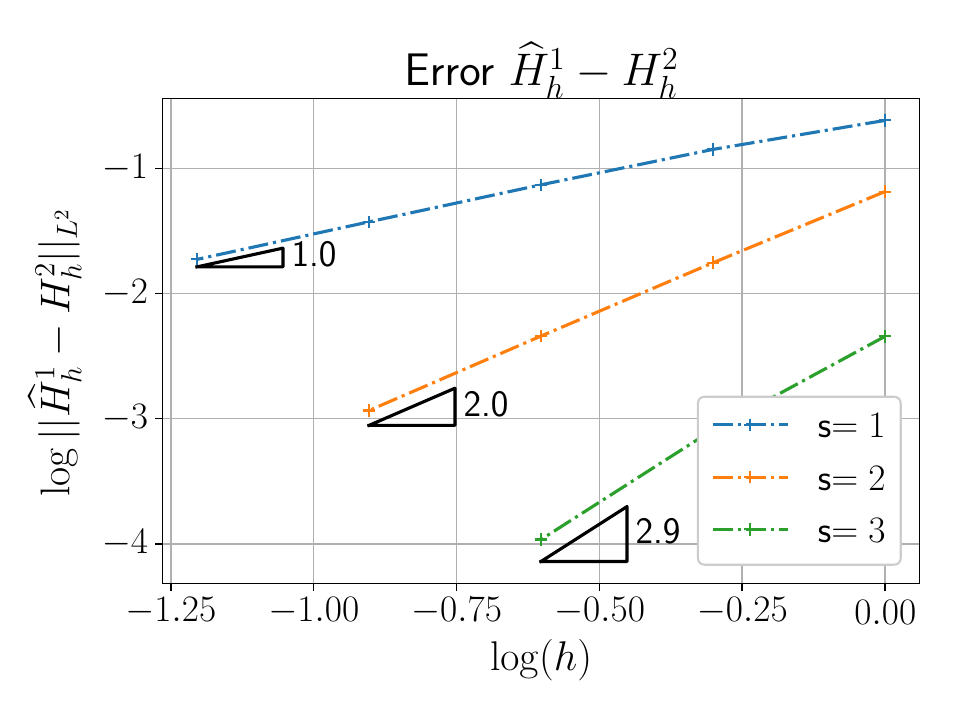}}%
\caption{$L^2$ difference of the dual representation of the solution for the Maxwell equation at $T_{\text{end}}=1$ for $\Delta t = \frac{1}{100}$.}%
\label{fig:diff_dual_maxwell}%
\end{figure}

In Fig. \ref{fig:conv_var_maxwell_primal} the convergence rate of the variables is plotted against the mesh size. The error is again measured in the Sobolev norm or the tangential trace norm for the facet variables at the final time $T_{\mathrm{end}}$.  It can be noticed that while the $L^2$ norm of the electric field two form converges with optimal order (cf. Fig \ref{fig:err_L2_E2}), its $H^{\div}$ norm loses one order (cf. Fig \ref{fig:err_Hdiv_E2}). This is due to the presence of the forcing term and the fact that in \firedrake the degree of freedom for Raviart Thomas are not computed via moments but instead via point normal evaluations (cf.  \url{https://github.com/FEniCS/fiat/issues/40}). Therefore, projecting on $H^{\div}$ space does not exactly commute with the divergence. This fact entails the observed loss of convergence. All other variables converge with optimal order. The convergence results for the dual formulation are shown in Fig.~\ref{fig:conv_var_maxwell_dual}. All variables converge with optimal as there is no forcing term in the magnetic field two form equation. In Fig. \ref{fig:diff_dual_maxwell}, the $L^2$ norm of the difference between the dual representation of variables. As in the dual field continuous Galerkin formulation, the dual representation of the variables converges under h-p refinement with order $h^s$. \\

The size reduction between the continuous and hybrid formulation is the same for the primal and dual system as the two formulation uses forms of the same degree and is reported in Table \ref{tab:dofs_cont_hyb_maxwell}. in The hybrid system has $45\%$ of degrees of freedom compared with the continuous formulation in the worst case. The size decreases rapidly to $30\%$ when the number of elements and polynomial degree is increased, showing a clear computational advantage with respect to the continuous formulation presented in Sec. \ref{sec:conGal_dis}.

\begin{table}[tbhp]
    \centering
    \begin{tabular}{c|c|c|c|c}
        Pol. Degree $s$ & $N_{\text{elem}}$ & $N^\circ$ dofs. continuous & $N^\circ$ dofs. hybrid & $N_{\text{hyb}}/N_{\text{cont}}$ \\ \hline
        \multirow{5}{*}{1} & 1 & 43 & 19 & 44\% \\
         & 2 & 290 & 98 & 38\% \\
         & 4 & 2140 & 604 & 28\% \\
         & 8 & 16472 & 4184 & 25\% \\
         & 16 & 129328 & 31024 & 24\% \\ \hline
         \multirow{4}{*}{2} & 1 & 164 & 74 & 45\% \\
          & 2 & 1156 & 436 & 37\% \\
          & 4 & 8696 & 2936 & 33\% \\
          & 8 & 67504 & 21424 & 32\% \\ \hline
          \multirow{3}{*}{3} & 1 & 399 & 165 & 41\% \\   
           & 2 & 2886 & 1014 & 35\% \\ 
           & 4 & 21972 & 6996 & 32\% \\ \hline
    \end{tabular}
    \caption{Size of the primal and dual system for the Maxwell equations: continuous and hybrid formulation.}
    \label{tab:dofs_cont_hyb_maxwell}
\end{table}

\begin{figure}[tbhp]%
\centering
\subfloat[][{$||\dual{E}^2||$ at $t=\pi/4 \; [\mathrm{s}]$}]{%
\label{fig:E2_1_fichera}%
\includegraphics[width=0.4\columnwidth]{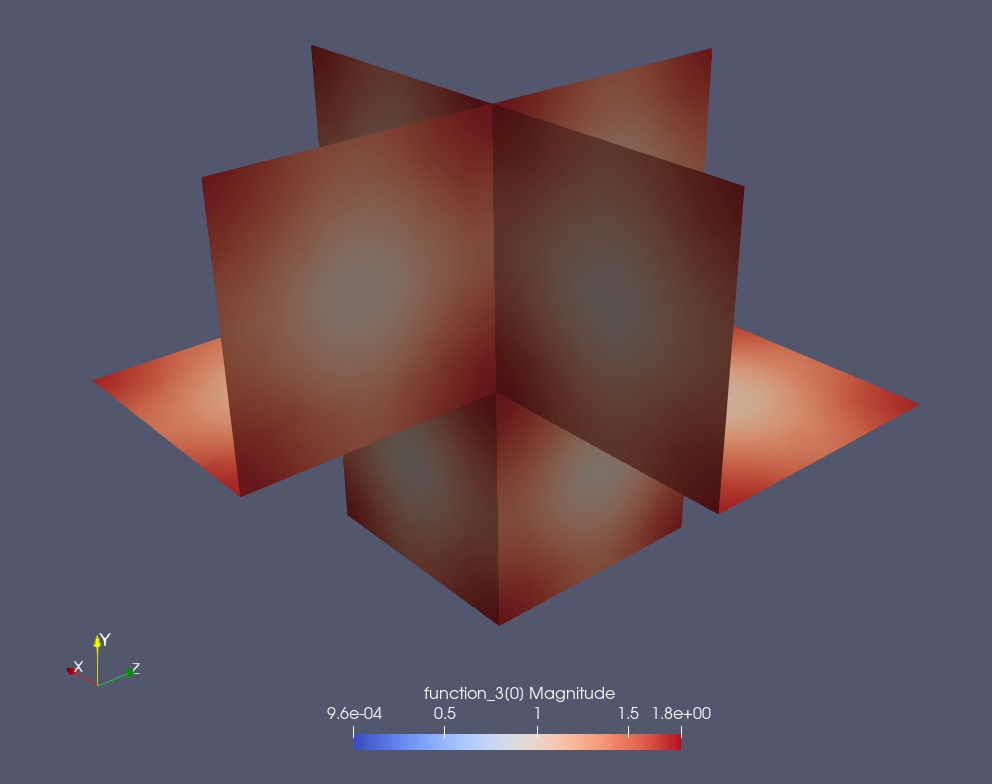}}%
\hspace{8pt}%
\subfloat[][{$||E^1||$ at $t=\pi/4 \; [\mathrm{s}]$}]{%
	\label{fig:E1_1_fichera}%
\includegraphics[width=0.4\columnwidth]{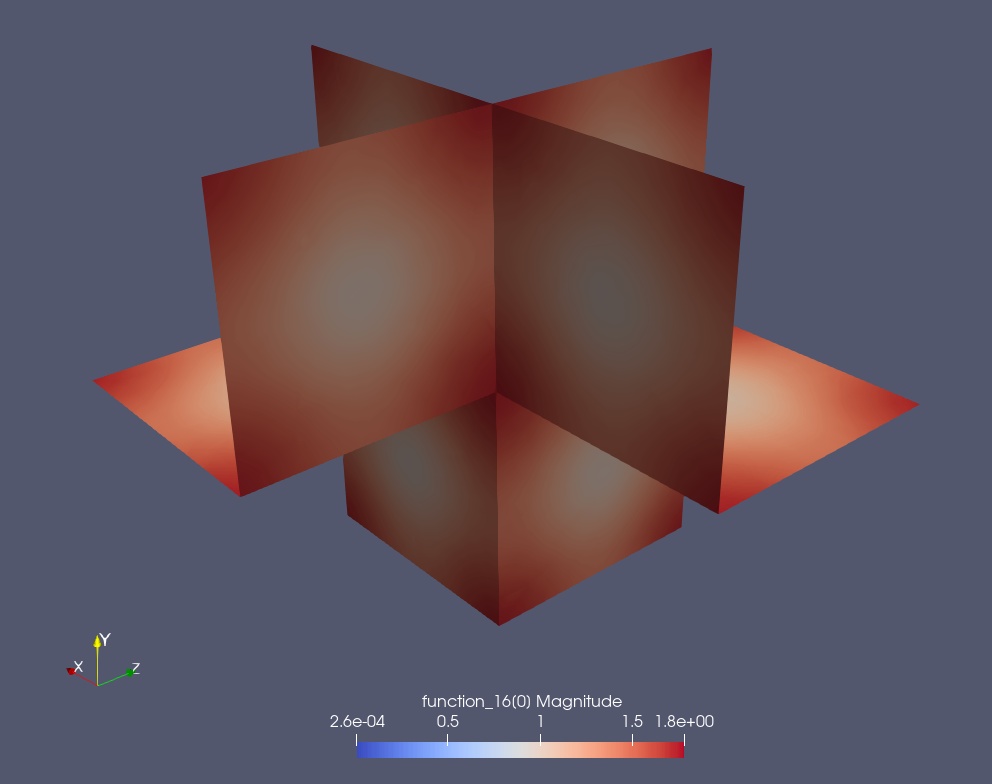}}\\
\subfloat[][{$||\dual{E}^2||$ at $t=\pi/2 \; [\mathrm{s}]$}]{%
\label{fig:E2_2_fichera}%
\includegraphics[width=0.4\columnwidth]{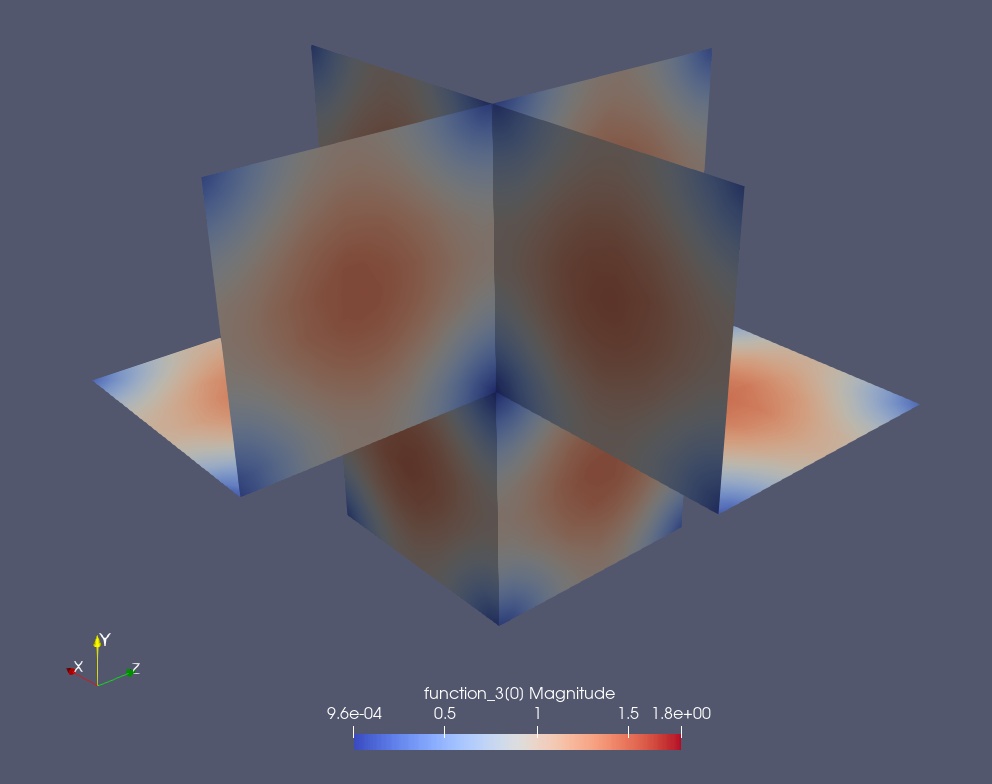}}%
\hspace{8pt}%
\subfloat[][{$||E^1||$ at $t=\pi/2 \; [\mathrm{s}]$}]{%
	\label{fig:E1_2_fichera}%
\includegraphics[width=0.4\columnwidth]{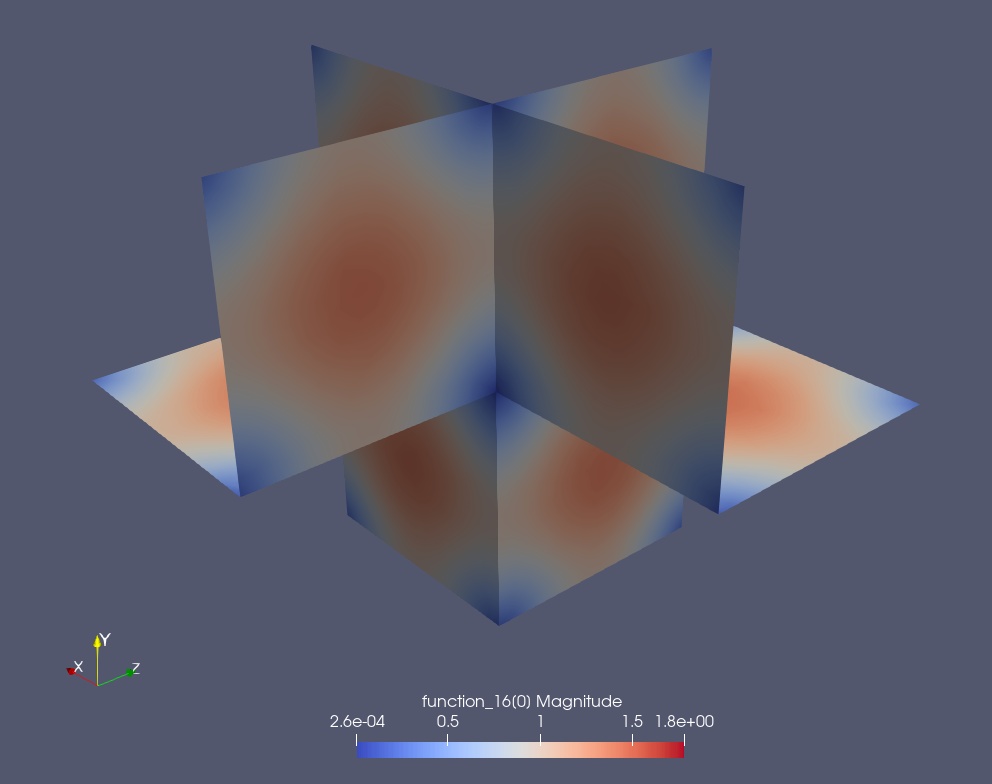}}\\
\subfloat[][{$||\dual{E}^2||$ at $t=3\pi/4 \; [\mathrm{s}]$}]{%
\label{fig:E2_3_fichera}%
\includegraphics[width=0.4\columnwidth]{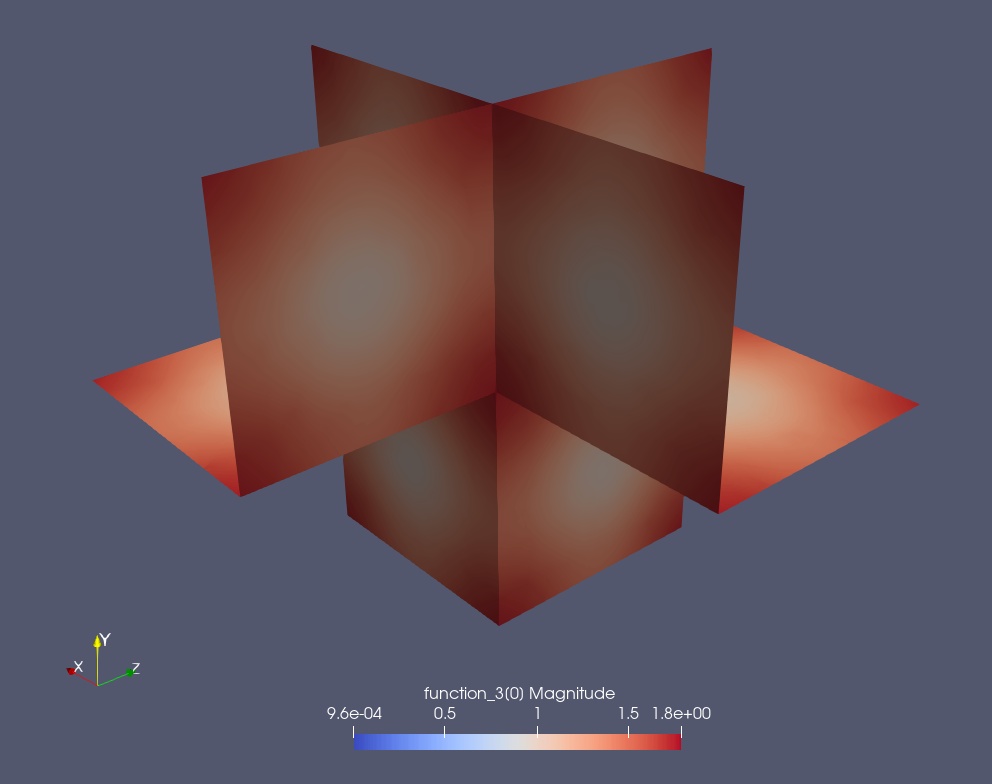}}%
\hspace{8pt}%
\subfloat[][{$||E^1||$ at $t=3\pi/4 \; [\mathrm{s}]$}]{%
	\label{fig:E1_3_fichera}%
\includegraphics[width=0.4\columnwidth]{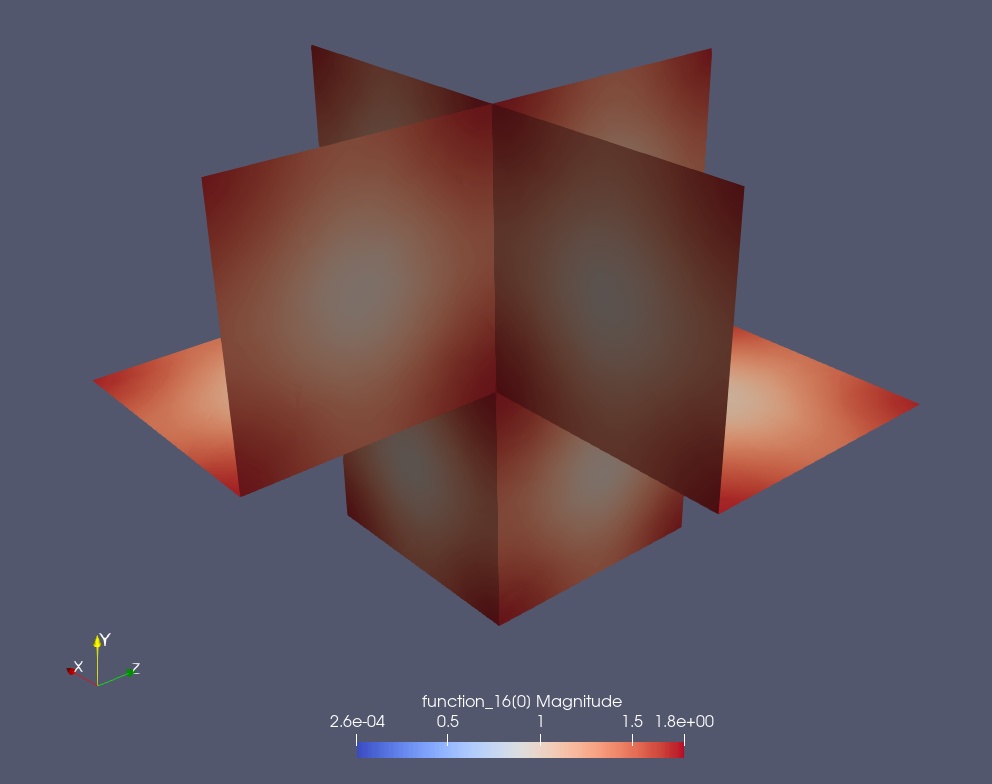}}\\
\subfloat[][{$||\dual{E}^2||$ at $t=\pi \; [\mathrm{s}]$}]{%
\label{fig:E2_4_fichera}%
\includegraphics[width=0.4\columnwidth]{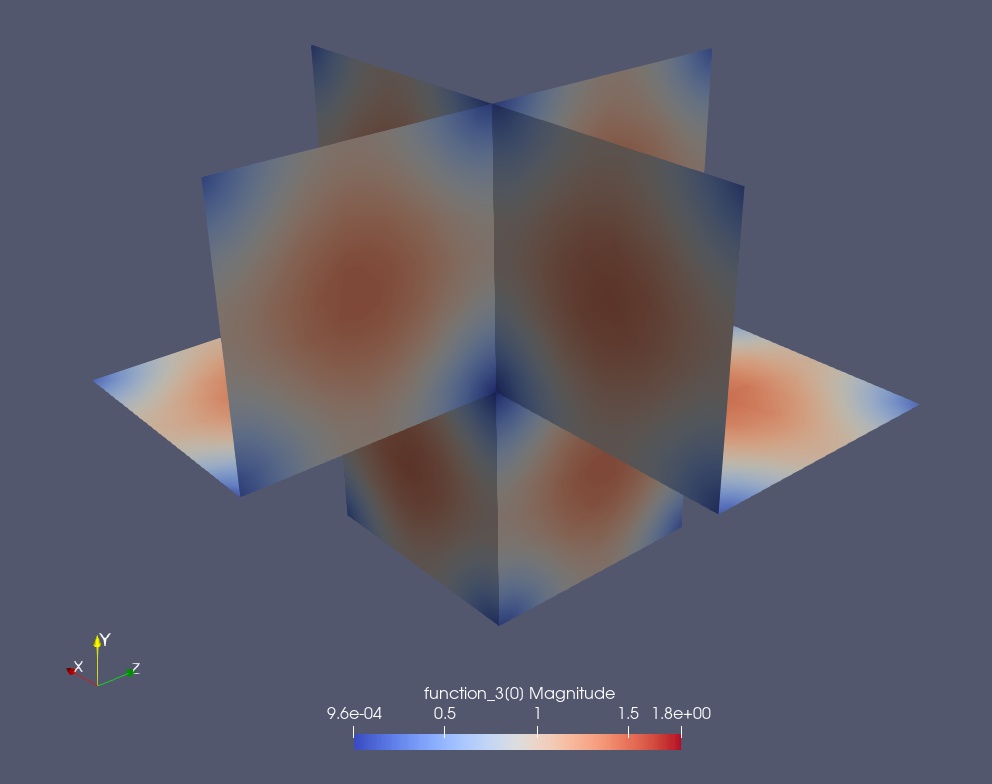}}%
\hspace{8pt}%
\subfloat[][{$||E^1||$ at $t=\pi \; [\mathrm{s}]$}]{%
	\label{fig:E1_4_fichera}%
\includegraphics[width=0.4\columnwidth]{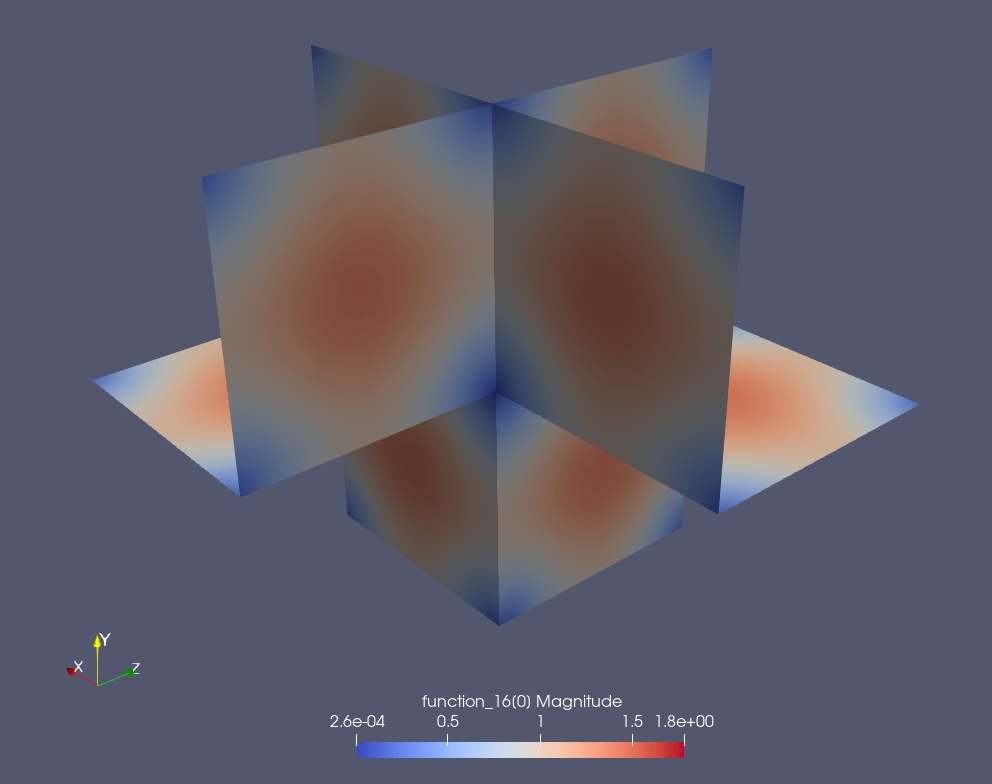}}
\caption{Snapshots of the magnitude of the electric calculated for the Fichera corner using the primal and dual system at different time instants.}
\label{fig:fichera_maxwell}%
\end{figure}

{
\subsection{A non convex domain: the Fichera corner}
As a last example, the Fichera corner geometry is considered. The domain is a cube with one octant removed $M = [-1, 1]^3/[-1,0]^3$. An unstructured mesh with size $h=1/8$  and Second order Raviart-Thomas and Nédélec finite elements are used. The electric permeability and magnetic permeability $\varepsilon=2, \; \mu=3/2$. The following analytical solution is considered
\begin{align*}
    E = \begin{pmatrix}
        \sin(2 t - 3 z)\\
        \sin(2 t - 3 x)\\
        \sin(2 t - 3 y)
        \end{pmatrix}, 
        \qquad
    H = \begin{pmatrix}
        \sin(2 t - 3 y)\\
        \sin(2 t - 3 z)\\
        \sin(2 t - 3 x)
        \end{pmatrix}, 
        \qquad
     j = \begin{pmatrix}
        \cos(2 t - 3 z)\\
        \cos(2 t - 3 x)\\
        \cos(2 t - 3 y)
        \end{pmatrix},
\end{align*}
where a current has been introduced as a forcing to make the given solution the true one. Electric boundary conditions are considered at the boundary. The total simulation is take to be one period of the given sinusoidal solution, i.e. $T_{\mathrm{end}}=\pi$. The time step is taken to be $\Delta T = \pi/100$. Snapshots of the solution at the time instant $\pi/4, \; \pi/2, \; 3\pi/4, \; \pi$ are reported in Fig. \ref{fig:fichera_maxwell}. The solution is the same for the first and third snapshots as for the second and fourth, as expected from the analytical solution. A perfect agreement between primal and dual formulation is again observed.
}

\section{Conclusion}\label{sec:conclusion}
In this work, the mixed dual field formulation is modified by taking the variable that do not undergo the exterior derivative to live in a broken finite element space. The resulting formulation is directly amenable to hybridization by introducing appropriate multipliers enforcing the continuity of the regular variable. The properties of the dual field scheme are left untouched and the hybrid formulation is highly advantageous since the broken variable, being local, is completely discarded from the global system, resulting in a huge computational gain obtained.\\

A development on this work would be the employment of algebraic dual polynomials \cite{jain2021algdual} to improve the condition number of the resulting matrix system.  Another important topic is the use of higher order time integration scheme and how to apply to those the static condensation procedure. The presented framework may be extended by introducing post-processing schemes. This may be achieved exploiting the primal-dual structure of the equations. Furthermore,  non-conforming and hybridizable discontinuous Galerkin (HDG) methods may be devised considering a more general local problem where the exterior derivative and the codifferential are taken weakly. Another important example of port-Hamiltonian system is the elastodynamics problem. A general exterior calculus framework for continuum mechanics is given in \cite{rashad2023intrinsic,rashad2025port}. In this case, a more general Hilbert complex is needed and more involved finite elements need to be used. The dual field method may also be applicable to this case. 

\section*{Code Availability}
The code used for the present work is hosted at: \\
\url{https://github.com/a-brugnoli/ph_hybridization}.

\section*{Acknowledgements}
The first author would like to thank Ari Stern for sharing the code used in \cite{stern2023hyb}.
 
\section*{Funding}

This work was supported by the PortWings project funded by the European Research Council [Grant Agreement No. 787675]

\bibliography{biblio}

\revone{
\appendix

\section{Algebraic realization of the weak formulation}\label{app:alg_mimetic}

We detail the algebraic computation of the different terms in the weak formulation. For a discussion of the algebraic realization arising in conforming finite element discretization, the reader can consult \cite[Sec. 5.1]{brugnoli2022df}. A different notation is used for inner and outer oriented forms. This is important when considering the natural duality pairing of forms, as the combination of an inner and outer oriented form leads to an orientation invariant operation. For the interconnection of different elements inner and outer oriented forms undergo a different treatment. Inner oriented forms do not carry information about the normal orientation, so they keep their sign from one element to the other Outer oriented forms carry information about the ambient normal orientation, so they change sign from one element to the other. This will be addressed when introducing the jump and average operators. At the discrete level, this distinction is yet to be implemented in standard finite element libraries like \firedrake. So the same basis and degrees of freedom are used for both.

\paragraph{General notation} Given a computational mesh $\mathcal{T}_h$, the set of facets of the mesh (the codimension 1 skeleton) is the union of the cells boundary $\mathcal{F}_h = \bigcup_{T \in \mathcal{T}_h} \partial T$, whereas the disjoint union is denoted by $\partial\mathcal{T}_h = \bigsqcup_{T \in \mathcal{T}_h} \partial T$. The notation $[\mathbf{A}]_{ij}$ will be used to denoted the element of matrix $\mathbf{A}$ corresponding to the $i$-th row and the $j$-th column.

\subsection{Local finite element forms and their facet version}

\paragraph{Basis expression for local forms}
Let $W_h^k(T) \subset H\Omega^k(T)$ be a finite (local) subcomplex.
 For a generic inner or outer oriented $k$-discrete conforming forms one has the following finite element expansion
\begin{equation}\label{eq:localbases}
     \mu_{h, T}^k = \sum_{i=1}^{N^k_T} \phi_{i, T}^{k}(\xi) \mu_{i, T}^k, \qquad \mu_{h, T}^k \in W_{h}^k(T),
\end{equation}
where $N_T^k$ is the number of degrees of freedom for $W_{h}^k(T)$, $\mu_{i, T}^k \in \bbR$ is the degree of freedom, and $\phi_{i, T}^{k} : M \rightarrow  W_{h}^k(T) \subset H\Omega^k(T)$ is a section of $W_{h}^k(T)$, corresponding to a finite element basis function. 

\paragraph{Inner product}

Given two discrete local forms $(\nu_{h, T}^k, \mu_{h, T}^k) \in W_{h, T}^k$, the inner product reads
\begin{equation}\label{eq:alg_inner}
    \inprDom[T]{\nu_{h, T}^k}{\mu_{h, T}^k} = (\bm{\nu}^k_T)^\top \mathbf{M}^k_T \bm{\mu}^k_T,
\end{equation}
where $\bm{\nu}^k_T, \; \bm{\mu}^k_T \in \bbR^{N_T^{k}}$ are the vectors collecting the degrees of freedom $\nu_{i, T}^k, \; \mu_{i, T}^k$ respectively and the mass matrix $\mathbf{M}^k_T \in \mathbb{R}^{N_T^{k}\times N_T^{k}}$ of order $k$ (symmetric and positive definite) is computed as  
$$[\mathbf{M}^k_T]_{ij} = \inprDom[T]{\phi_{i, T}^{k}}{\phi_{j, T}^{k}}.$$

\paragraph{Exterior derivative}
The expression of the exterior derivative is here specialized for the inner product. Given a form $\nu_{h, T}^{k+1} \in W_{h}^{k+1}(T)$ and $\mu_{h, T}^k \in {W}_{h}^{k}(T)$ ($k\le n-1$), the inner product of $\nu_{h, T}^{k+1}$ and $\d \mu_{h, T}^k$ is expressed by
\begin{equation}\label{eq:alg_local_d}
    \inprDom[T]{\nu_{h, T}^{k+1}}{\d\mu_{h, T}^k} = (\bm{\nu}^{k+1}_T)^\top \mathbf{D}^{k}_T \bm{\mu}^{k}_T
\end{equation}
where $\mathbf{D}^{k}_T \in \mathbb{R}^{N_T^{k+1} \times N_T^{k}}$ is computed as  $[\mathbf{D}^{k}_T]_{ij} = \inprDom[T]{\phi_{i, T}^{k+1}}{\d\phi_{j, T}^{k}}.$

\paragraph{Trace}
For a discrete differential form of order $k \le n-1$ the trace operator is also defined. When the trace of a local form is considered the simplices that do not lie on the boundary can be discarded in the expansion. The trace of a discrete $k$-form $\mu_{h, T} \in {W}_{h}^k(T)$ reads
\begin{equation}
    \mu_{h, \partial T}^{k, \bm{t}} =\tr \mu_{h, T}^k = \sum_{i=1}^{N_T^{k}} \tr(\phi_{i, T}^{k}(\xi)) \mu_{i, T}^k = \sum_{l=1}^{N_{\partial T}^{k}} \psi_{l, \partial T}^{k, \bm{t}}(\xi) \mu_{l, \partial T}^{k, \bm{t}}.
\end{equation}
where $N_{\partial T}^{k}$ is the number of degrees of freedom (with $k\le j\le n-1$) along the boundary for a polynomial differential form of order $k$. The degrees of freedom along the boundary $\mu_{l, T}^k$ are associated to mesh entities lying on the boundary. The trace matrix simply collects them
\begin{equation}\label{eq:alg_local_trace}
    \bm{\mu}_{\partial T}^{k, \bm{t}} = \mathbf{T}^k_{\partial T} \bm{\mu}^k_T, \qquad
    [\mathbf{T}^k_{\partial T}]_{li} = 
\begin{cases}
    1, \quad \text{if } \psi_{l, \partial T}^{k, \bm{t}}(\xi) \equiv \tr(\phi_{i, T}^{k}(\xi))\\
    0, \quad \text{otherwise},
\end{cases} \quad 
\begin{aligned}
    \forall l = 1, \dots, N_{\partial T}^{k}, \\
    \forall i = 1, \dots, N_T^{k}.
\end{aligned}  
\end{equation}

\paragraph{Inner product of facets elements} 
The finite element expansion of a facet element $\mu_{h, \partial T}^{k, \bm{t}} \in W_h^{k, \bm{t}}(T)$ reads
\begin{equation*}
    \mu_{h, \partial T}^{k, \bm{t}} = \sum_{l=1}^{N_{\partial T}^k} \psi_{l, \partial T}^{k, \bm{t}}(\xi) \mu_{l, \partial T}^{k, \bm{t}}.
\end{equation*}

Since $W_h^{k, \bm{t}}(T) \in L^2(\partial T)$,  the dual space  $W_h^{k, \bm{n}}(T)$ is identified with $W_h^{k, \bm{t}}(T)$. So given $\nu_{h, \partial T}^{k, \bm{n}} \in W_h^{k, \bm{n}}(T) = W_h^{k, \bm{t}}(T)$ its pairing with $\mu_{h, \partial T}^{k, \bm{t}}$ is the $L^2$ inner product over the boundary $\partial T$
\begin{equation*}
    \inprBd[\partial T]{\nu_{h, \partial T}^{k, \bm{n}}}{\mu_{h, \partial T}^{k, \bm{t}}} = (\bm{\nu}_{\partial T}^{k, \bm{n}})^\top \mathbf{M}_{\partial T}^k \bm{\mu}_{\partial T}^{k, \bm{t}}.
\end{equation*}
where the mass matrix at the boundary is computed as 
\begin{equation*}
    [\mathbf{M}_{\partial T}^k]_{ij} = \inprBd[\partial T]{\psi_{i, \partial T}^{k, \bm{t}}(\xi)}{\psi_{j, \partial T}^{k, \bm{t}}(\xi)}.
\end{equation*}
If one of he two facet forms is obtained via the trace operator then 
\begin{equation*}
    \inprBd[\partial T]{\nu_{h, \partial T}^{k, \bm{n}}}{\tr \mu_{h, T}^{k}} = (\bm{\nu}_{\partial T}^{k, \bm{n}})^\top \mathbf{M}_{\partial T}^k  \mathbf{T}^k_{\partial T} \bm{\mu}_{T}^{k}.
\end{equation*}

\subsection{Broken spaces of differential forms}

\paragraph{Basis expression for broken forms}
Since $W_h^k = \prod_{T\in \mathcal{T}} W_h^k(T)$, broken differential forms $\mu_{h}^k \in W_h^k$ are obtained by summing the contribution of all cells in the mesh
\begin{equation}\label{eq:brokenbases}
     \mu_{h}^k = \sum_{T \in \mathcal{T}_h}\mu_{h, T}^k, \qquad \mu_{h, T}^k \in W_{h}^k(T),
\end{equation}
All the local operation translates into block diagonal structures. 
\begin{itemize}
    \item Inner product: given two discrete broken forms $(\nu_{h}^k, \mu_{h}^k) \in W_{h}^k$, the inner product reads
\begin{equation*}
    \inprDom[\mathcal{T}_h]{\nu_{h}^k}{\mu_{h}^k} = (\bm{\nu}^k)^\top \mathbf{M}^k_{\mathcal{T}_h}\bm{\mu}^k, \qquad 
\end{equation*}
where $\bm{\nu}^k, \; \bm{\mu}^k \in \bbR^{N^k}$ with $N^k = \sum_{T\in \mathcal{T}_h}N_T^{k}$ are the vectors collecting the degrees of freedom $\mu_i^k, \; \mu_i^k$ respectively. The mass matrix is given by the block diagonal matrix  
\begin{equation*}
    \mathbf{M}^k_{\mathcal{T}_h} = \mathrm{Diag}_{T\in \mathcal{T}_h}(\mathbf{M}^k_T).
\end{equation*}
\item Exterior derivative: given broken forms $\nu_{h}^{k+1} \in W_{h}^{k+1}$ and $\mu_{h, T}^k \in {W}_{h}^{k}$ ($k\le n-1$), the inner product of $\nu_{h}^{k+1}$ and $\d \mu_{h}^k$ is expressed by
\begin{equation*}
    \inprDom[\mathcal{T}_h]{\nu_{h}^{k+1}}{\d\mu_{h}^k} = (\bm{\nu}^{k+1})^\top \mathbf{D}^{k}_{\mathcal{T}_h} \bm{\mu}^{k}
\end{equation*}
where $\mathbf{D}^{k}_{\mathcal{T}_h}$ is again block diagonal $\mathbf{D}^k_{\mathcal{T}_h} = \mathrm{Diag}_{T\in \mathcal{T}_h}(\mathbf{D}^k_T)$.
\item Trace: the trace of a broken form $\mu_{h, T}^k \in {W}_{h}^{k}$ ($k\le n-1$) over $\partial \mathcal{T}_h$
\begin{equation*}
    \tr_{} \mu_{h}^k = \tr \sum_{T \in \mathcal{T}_h} \mu_{h, T}^k  = \sum_{\partial T \in \partial \mathcal{T}_h} \mu_{h, \partial T}^{k, \bm{t}} = \sum_{\partial T \in \partial \mathcal{T}_h}\sum_{l=1}^{N_{\partial T}^{k}} \psi_{l, \partial T}^{k}(\xi) \mu_{l, \partial T}^{k, \bm{t}},
\end{equation*}
gives rises to a block diagonal operator
\begin{equation*}
    \bm{\mu}^{k, \bm{t}} = \mathbf{T}^k_{\partial \mathcal{T}_h} \bm{\mu}^k, \qquad
    \mathbf{T}^k_{\partial \mathcal{T}_h} = \mathrm{Diag}_{\partial T \in \partial \mathcal{T}_h} (\mathbf{T}^k_{\partial T}).
\end{equation*}
\end{itemize}

\paragraph{Inner product of broken facets elements}
Broken facet forms $\mu_h^{k, \bm{t}} \in W_h^{k, \bm{t}}$ are expressed in a basis summing the contribution of the boundary of each cell
\begin{equation*}
    \mu_{h}^{k, \bm{t}} = \sum_{\partial T \in \partial \mathcal{T}_h}\mu_{h, \partial T}^{k, \bm{t}}, \qquad \mu_{h, \partial T}^{k, \bm{t}} \in W_{h}^{k, \bm{t}}(T).
\end{equation*}
Broken facet forms can be paired together, giving rise to a block diagonal mass facets matrix
\begin{equation*}
    \inprBd[\partial \mathcal{T}_h]{\nu_{h}^{k, \bm{n}}}{\mu_{h}^{k, \bm{t}}} = (\bm{\nu}^{k, \bm{n}})^\top \mathbf{M}_{\partial \mathcal{T}_h}^k \bm{\mu}^{k, \bm{t}}.
\end{equation*}
where the mass matrix at the boundary is computed as 
\begin{equation*}
    \mathbf{M}_{\partial \mathcal{T}_h}^k = \mathrm{Diag}_{\partial T\in \partial\mathcal{T}_h} \mathbf{M}_{\partial T}.
\end{equation*}
If one of he two facet forms is obtained via the trace, then
\begin{equation*}
    \inprBd[\partial \mathcal{T}_h]{\nu_{h}^{k, \bm{n}}}{\tr \mu_{h}^{k}} = (\bm{\nu}^{k, \bm{n}})^\top \mathbf{M}_{\partial \mathcal{T}_h}^k  \mathbf{T}^k_{\partial \mathcal{T}_h} \bm{\mu}^{k}.
\end{equation*}

\subsection{Inner product of broken and unbroken facets forms}

Unbroken discrete forms over the facets are obtained by making the broken version single valued, i.e. $V_h^{k, \bm{t}} = W_h^{k, \bm{t}} \cap V^k$. When an unbroken facet form $v_h^{k, \bm{t}} \in V_h^{k, \bm{t}}$ is paired with a broken one $\mu_h^{k, \bm{n}} \in W_h^{k, \bm{n}}$, the integral over the disjoint union $\partial T_h$ can be converted to an integral over the set of facets $\mathcal{F}_h$ by introducing the jump operator. As anticipated at the beginning of the appendix, inner and outer oriented forms behave differently under interconnection. When the inner product of broken and unbroken forms is considered an average appear. 

\paragraph{Average operator} 
The average of a form is simply given by
\begin{equation*}
\begin{aligned}
&\llbracket \cdot\rrbracket :  W_h^{k, \bm{n}} \longrightarrow V_h^{k, \bm{n}}, \\
& \llbracket \mu_h^{k, \bm{n}}\rrbracket = \begin{cases}
    \frac{1}{2} (\mu_{h, \partial T^+}^{k, \bm{n}} + \mu_{h, \partial T^-}^{k, \bm{n}}), \\
    \mu_{h, \partial T}^{k, \bm{n}},
    \end{cases} \qquad
    \begin{aligned}
    &\text{on the face } f \in \partial T^+ \cap \partial T^-, \\
    &\text{on the face } f \in \partial T \cap \partial M. 
    \end{aligned}
\end{aligned}
\end{equation*}

Consequently the duality pairing of $v_h^{k, \bm{t}} \in V_h^{k, \bm{t}}$ and $\mu_h^{k, \bm{n}} \in W_h^{k, \bm{n}} = W_h^{k, \bm{t}}$ over $\partial \mathcal{T}_h$ reads
\begin{equation*}
    \inprBd[\partial \mathcal{T}_h]{v_h^{k, \bm{t}}}{\mu_h^{k, \bm{n}}} = \inprBd[\mathcal{F}_h]{v_h^{k, \bm{t}}}{\llbracket \mu_h^{k, \bm{n}}\rrbracket}.
\end{equation*}
Computing the expression in a basis gives the following matrix expression
\begin{equation*}
    \inprBd[\mathcal{F}_h]{v_h^{k, \bm{t}}}{\llbracket \mu_h^{k, \bm{n}}\rrbracket} = \mathbf{v}^{k, \bm{t}} \mathbf{M}^k_{\mathcal{F}_h} \mathbf{\Xi}^{k}_{\mathcal{F}_h} \bm{\mu}^{k, \bm{n}}.
\end{equation*}
where $\bm{\Xi}^{k}_{\mathcal{F}_h}$ is the matrix expression of the jump operator.

\subsection{Duality pairing of inner and outer oriented forms over the domain boundary}

An unbroken facet form $v^{k, \bm{t}}_h \in V_h^{k, \bm{t}}$ can be restricted to the boundary of the domain $\partial M$ by taking its trace
\begin{equation*}
    v^{k, \bm{t}}_{h, \partial} = \tr v^{k, \bm{t}}_h, \qquad \text{or algebraically} \quad  \mathbf{v}^{k, \bm{t}}_\partial = \mathbf{T}^{k, \bm{t}}_{\partial M} \mathbf{v}^{k, \bm{t}}.
\end{equation*}

Given an unbroken facet form $v^{k, \bm{t}}_h \in V_h^{k, \bm{t}}$ and an outer oriented control variables defined on the domain boundary $\dual{u}^{n-k-1}_h \in \tr \dual{V}^{n-k-1}_h|_{\partial M}$ their duality pairing over $\partial M$ gives the following algebraic expression
\begin{equation*}
    \dualprBd[\partial M]{\tr v^{k, \bm{t}}_h}{\dual{u}^{n-k-1}_h} = (\mathbf{T}^{k, \bm{t}}_{\partial M}\mathbf{v}^{k, \bm{t}})^\top \mathbf{\Psi}_{\partial M}^{n-k-1} \dual{\mathbf{u}}^{n-k-1}.
\end{equation*}
The matrix $\mathbf{\Psi}_{\partial M}^{n-k-1}$ is in general rectangular and rank deficient. Analogously, an outer oriented from can be paired with an outer oriented control. Furthermore the boundary integral can be restricted to the subpartitions $\Gamma_1,\; \Gamma_2$.

}

\end{document}